\documentclass[11pt]{article}
\usepackage[T1]{fontenc}
\usepackage[utf8]{inputenc}
\usepackage{lmodern}
\usepackage[a4paper,margin=23mm]{geometry}
\usepackage{microtype}
\usepackage{setspace}
\usepackage{amsmath,amssymb,amsthm,mathtools}
\numberwithin{equation}{section}
\allowdisplaybreaks[1]
\usepackage{enumitem}
\usepackage{booktabs}
\usepackage{float}
\usepackage{array}
\usepackage{longtable}
\usepackage{xcolor}
\usepackage[hypertexnames=false]{hyperref}
\usepackage[nameinlink,capitalise,noabbrev]{cleveref}
\usepackage{graphicx}
\usepackage{tikz}
\usetikzlibrary{arrows.meta,positioning,shapes.geometric}

\definecolor{gridgrey}{HTML}{F1F5F9}

\newcommand{\drawGridBackground}[3]{
    \foreach \relX in {0,...,#1} {
        \pgfmathparse{int(mod(\relX + #3, 2))}
        \ifnum\pgfmathresult=1
            \fill[gridgrey] (\relX, 0) rectangle (\relX+1, #2+1);
        \fi
    }
    \ifnum#1>0
        \foreach \x in {1,...,#1} {
            \draw[black, line width=0.5pt] (\x, 0) -- (\x, #2+1);
        }
    \fi
    \draw[dashed, gray!60, line width=0.4pt] (0, 0) rectangle (#1+1, #2+1);
}

\newcommand{\standardBox}[5][2]{
    \pgfmathparse{int(mod(#2 + #5, 2))}
    \edef\isGrey{\pgfmathresult}

    \node[draw=black, fill=white, thick,
          minimum width={0.45cm + (#1 - 2) * 0.3cm},
          minimum height=0.45cm,
          text width={0.45cm + (#1 - 2) * 0.3cm},
          align=center, inner sep=0.5pt,
          font=\tiny\sffamily\bfseries, transform shape]
          at ({#2 + 0.5 + (#1 - 2) * 0.5}, #3 + 0.5) {
            \ifnum\isGrey=1 #4 \else \rotatebox{90}{{\centering #4}} \fi
          };
}

\newcommand{\specialBox}[3]{
    \begin{scope}[shift={(#1+0.5, #2+0.5)}, transform shape]
        \pgfmathparse{int(mod(#1 + #3, 2))}
        \ifnum\pgfmathresult=1
            \def\canvasColor{gridgrey}\def\inverseColor{white}
        \else
            \def\canvasColor{white}\def\inverseColor{gridgrey}
        \fi
        \fill[\inverseColor] (-0.6, -0.357) rectangle (0.6, 0.357);
        \draw[line width=0.5pt, fill=\canvasColor] (-0.5, 0.357)
              arc[start angle=180, end angle=360, x radius=0.5, y radius=0.214];
        \draw[line width=0.5pt, fill=\canvasColor] (-0.5, -0.357)
              arc[start angle=180, end angle=0, x radius=0.5, y radius=0.214];
    \end{scope}
}

\newcommand{\OnPicture}[2]{%
  \vcenter{\hbox{%
    \begin{tikzpicture}[
      x=0.23cm,
      y=0.29cm,
      line width=0.4pt,
      line cap=round,
      line join=round
    ]
      \path[use as bounding box]
        (-#1,-1.08) rectangle (#1,1.08);
      #2
    \end{tikzpicture}%
  }}%
}

\newcommand{\OnOne}{%
  \OnPicture{0.12}{
    \draw (0,1)--(0,-1);
  }%
}

\newcommand{\OnTwo}[1]{%
  \OnPicture{0.55}{
    \ifcase#1\relax
    \or % 1
      \draw
        (-0.5,1)--(-0.5,-1)
        (0.5,1)--(0.5,-1);
    \or % 2
      \draw (-0.5,1)
        .. controls (-0.5,0) and (0.5,0) .. (0.5,1);
      \draw (-0.5,-1)
        .. controls (-0.5,0) and (0.5,0) .. (0.5,-1);
    \or % 3
      \draw
        (-0.5,1)--(0.5,-1)
        (0.5,1)--(-0.5,-1);
    \fi
  }%
}

\newcommand{\OnThree}[1]{%
  \OnPicture{1.05}{
    \ifcase#1\relax

    \or % 1
      \draw
        (-1,1)--(-1,-1)
        (0,1)--(0,-1)
        (1,1)--(1,-1);

    \or % 2
      \draw (-1,1)--(-1,-1);
      \draw (0,1)
        .. controls (0,0) and (1,0) .. (1,1);
      \draw (0,-1)
        .. controls (0,0) and (1,0) .. (1,-1);

    \or % 3
      \draw (1,1)--(1,-1);
      \draw (-1,1)
        .. controls (-1,0) and (0,0) .. (0,1);
      \draw (-1,-1)
        .. controls (-1,0) and (0,0) .. (0,-1);

    \or % 4
      \draw (-1,-1)--(1,1);
      \draw (-1,1)
        .. controls (-1,0) and (0,0) .. (0,1);
      \draw (0,-1)
        .. controls (0,0) and (1,0) .. (1,-1);

    \or % 5
      \draw (-1,1)--(1,-1);
      \draw (0,1)
        .. controls (0,0) and (1,0) .. (1,1);
      \draw (-1,-1)
        .. controls (-1,0) and (0,0) .. (0,-1);

    \or % 6
      \draw
        (-1,1)--(0,-1)
        (0,1)--(-1,-1)
        (1,1)--(1,-1);

    \or % 7
      \draw
        (-1,1)--(-1,-1)
        (0,1)--(1,-1)
        (1,1)--(0,-1);

    \or % 8
      \draw (-1,1)
        .. controls (-1,0.1) and (1,0.1) .. (1,1);
      \draw (-1,-1)
        .. controls (-1,0) and (0,0) .. (0,-1);
      \draw (0,1)--(1,-1);

    \or % 9
      \draw (-1,1)
        .. controls (-1,0) and (0,0) .. (0,1);
      \draw (-1,-1)
        .. controls (-1,-0.1) and (1,-0.1) .. (1,-1);
      \draw (1,1)--(0,-1);

    \or % 10
      \draw (-1,1)
        .. controls (-1,0.1) and (1,0.1) .. (1,1);
      \draw (0,-1)
        .. controls (0,0) and (1,0) .. (1,-1);
      \draw (0,1)--(-1,-1);

    \or % 11
      \draw (0,1)
        .. controls (0,0) and (1,0) .. (1,1);
      \draw (-1,-1)
        .. controls (-1,-0.1) and (1,-0.1) .. (1,-1);
      \draw (-1,1)--(0,-1);

    \or % 12
      \draw
        (-1,1)--(1,-1)
        (0,1)--(-1,-1)
        (1,1)--(0,-1);

    \or % 13
      \draw (-1,1)
        .. controls (-1,0.1) and (1,0.1) .. (1,1);
      \draw (-1,-1)
        .. controls (-1,-0.1) and (1,-0.1) .. (1,-1);
      \draw (0,1)--(0,-1);

    \or % 14
      \draw
        (-1,1)--(0,-1)
        (0,1)--(1,-1)
        (1,1)--(-1,-1);

    \fi
  }%
}

\hypersetup{
  colorlinks=true,
  linkcolor=blue!55!black,
  citecolor=green!40!black,
  urlcolor=blue!65!black,
  pdftitle={Triangle Relation Planar Algebras: Cohomology, Classification, and New Families},
  pdfauthor={Jiaqi Li and Zhengwei Liu}
}

\newtheorem{theorem}{Theorem}[section]
\newtheorem{proposition}[theorem]{Proposition}
\newtheorem{lemma}[theorem]{Lemma}
\newtheorem{corollary}[theorem]{Corollary}

\theoremstyle{definition}
\newtheorem{definition}[theorem]{Definition}

\theoremstyle{remark}
\newtheorem{remark}[theorem]{Remark}
\newtheorem{question}[theorem]{Question}

\newcommand{\C}{\mathbb C}
\newcommand{\Q}{\mathbb Q}
\newcommand{\Z}{\mathbb Z}

\newcommand{\supp}{\operatorname{Supp}}

\newcommand{\im}{\operatorname{im}}

\newcommand{\dual}[1]{\overline{#1}}
\newcommand{\F}{\mathcal F}
\newcommand{\Espan}{\mathcal E}

\title{Triangle Relation Planar Algebras: Cohomology, Classification, and New Families}

\author{
    Jiaqi Li$^{1}$\thanks{\texttt{lijiaqi25@mails.tsinghua.edu.cn}}
  \qquad
  Zhengwei Liu$^{2,3,4}$\thanks{\texttt{liuzhengwei@mail.tsinghua.edu.cn}}\\[0.8em]
  {\small $^{1}$Qiuzhen College, Tsinghua University, Beijing, China}\\[0.3em]
  {\small $^{2}$Yau Mathematical Sciences Center,}\\
  {\small Tsinghua University, Beijing, China}\\[0.3em]
  {\small $^{3}$Beijing Institute of Mathematical Sciences and Applications,}\\
  {\small Beijing, China}\\[0.3em]
  {\small $^{4}$Department of Mathematical Sciences,}\\
  {\small Tsinghua University, Beijing, China}
}

\date{}

\usepackage{needspace}
\newcommand{\tc}{\widetilde c}

\newcommand{\BMW}{\mathcal W}
\newcommand{\PB}{\mathcal P^B}
\definecolor{authorblue}{RGB}{30,76,145}

\makeatletter
\renewcommand{\@seccntformat}[1]{\csname the#1\endcsname\hspace{1em}}
\renewcommand\subsection{\@startsection{subsection}{2}{\z@}%
  {-3.25ex \@plus -1ex \@minus -.2ex}%
  {1.25ex \@plus .2ex}%
  {\normalfont\large\bfseries}}
\makeatother
\newcommand{\opPic}[1]{\vcenter{\hbox{\begin{tikzpicture}[
 x=0.8cm,y=0.8cm,line width=0.5pt,line cap=round,line join=round,
 opbox/.style={draw,fill=white,minimum width=0.8cm,minimum height=0.55cm,
 inner sep=0pt,font=\normalsize}]#1\end{tikzpicture}}}}
\newcommand{\opEditor}[1]{\vcenter{\hbox{\begin{tikzpicture}[
 x=0.3cm,y=0.7cm,scale=1.15]#1\end{tikzpicture}}}}
\newcommand{\opP}[1]{\opEditor{\drawGridBackground{2}{0}{0}\standardBox{1}{0}{$p_{#1}$}{0}}}
\newcommand{\opLoop}[1]{\opPic{\ifnum#1=0
 \fill[gridgrey] (0,0) circle[radius=0.40];
 \else \fill[gridgrey] (-0.55,-0.55) rectangle (0.55,0.55);
 \fill[white] (0,0) circle[radius=0.40];\fi
 \draw (0,0) circle[radius=0.40];}}
\newcommand{\opCap}[1]{\opPic{
 \fill[gridgrey] (-0.24,-0.95)--(-0.24,0.60)
 .. controls (-0.24,1.02) and (0.24,1.02) .. (0.24,0.60)--(0.24,-0.95)--cycle;
 \draw (-0.24,-0.95)--(-0.24,0.60)
 .. controls (-0.24,1.02) and (0.24,1.02) .. (0.24,0.60)--(0.24,-0.95);
 \ifnum#1=1 \node[opbox] at (0,0) {$p_i$};\fi}}
\newcommand{\opPartial}{\opPic{
 \fill[gridgrey] (-0.24,-0.95) rectangle (1.13,0.95);
 \fill[white] (0.24,0.30)--(0.24,0.55)
 .. controls (0.24,1.02) and (0.95,1.02) .. (0.95,0)
 .. controls (0.95,-1.02) and (0.24,-1.02) .. (0.24,-0.55)--(0.24,-0.30)--cycle;
 \draw (-0.24,0.95)--(-0.24,-0.95);
 \draw (0.24,0.30)--(0.24,0.55)
 .. controls (0.24,1.02) and (0.95,1.02) .. (0.95,0)
 .. controls (0.95,-1.02) and (0.24,-1.02) .. (0.24,-0.55)--(0.24,-0.30);
 \node[opbox] at (0,0) {$p_i$};}}
\newcommand{\opStrand}{\opPic{\fill[gridgrey] (0,-0.95) rectangle (0.50,0.95);
 \draw (0,-0.95)--(0,0.95);}}
\newcommand{\opMult}{\opEditor{\drawGridBackground{2}{1}{0}
 \standardBox{1}{1}{$p_i$}{0}\standardBox{1}{0}{$p_j$}{0}}}
\newcommand{\opConv}{\opPic{
 \fill[gridgrey] (0.25,-1.1) rectangle (2.25,1.1);
 \fill[white] (0.85,0.30)--(0.85,0.60)
 .. controls (0.85,1.20) and (1.65,1.20) .. (1.65,0.60)
 --(1.65,-0.60)
 .. controls (1.65,-1.20) and (0.85,-1.20) .. (0.85,-0.60)--cycle;
 \draw (0.25,-1.1)--(0.25,1.1) (2.25,-1.1)--(2.25,1.1);
 \draw (0.85,0.30)--(0.85,0.60)
 .. controls (0.85,1.20) and (1.65,1.20) .. (1.65,0.60)--(1.65,0.30);
 \draw (0.85,-0.30)--(0.85,-0.60)
 .. controls (0.85,-1.20) and (1.65,-1.20) .. (1.65,-0.60)--(1.65,-0.30);
 \node[opbox] at (0.55,0) {$p_i$};\node[opbox] at (1.95,0) {$p_j$};}}
\newcommand{\opIdentity}{\opEditor{\drawGridBackground{2}{0}{0}}}
\newcommand{\opJones}{\opEditor{\drawGridBackground{2}{0}{0}\specialBox{1}{0}{0}}}
\newcommand{\TRpicture}[1]{%
  \vcenter{\hbox{%
    \begin{tikzpicture}[x=0.42cm,y=0.68cm]
      #1
    \end{tikzpicture}%
  }}%
}

\newcommand{\TRbackground}{%
  \fill[gridgrey] (1,0) rectangle (2,3);
  \fill[gridgrey] (3,0) rectangle (4,3);
  \foreach \x in {1,2,3}{
    \draw[black,line width=0.5pt] (\x,0)--(\x,3);
  }
}

\newcommand{\opT}[1]{%
  \TRpicture{%
    \TRbackground
    \ifnum#1=0
      \standardBox{1}{2}{$p_i$}{0}
      \standardBox{2}{1}{$p_j$}{0}
      \standardBox{1}{0}{$p_k$}{0}
    \else
      \standardBox{2}{2}{$p_i$}{0}
      \standardBox{1}{1}{$p_j$}{0}
      \standardBox{2}{0}{$p_k$}{0}
    \fi
  }%
}

\newcommand{\TRfamily}[3]{%
  \TRpicture{%
    \TRbackground
    \ifcase#1\relax
    \or
      \standardBox{1}{2}{$p_{#2}$}{0}
      \standardBox{2}{0}{$p_{#3}$}{0}
    \or
      \standardBox{1}{2}{$p_{#2}$}{0}
      \standardBox{1}{0}{$p_{#3}$}{0}
      \specialBox{2}{1}{0}
    \or
      \standardBox{2}{2}{$p_{#2}$}{0}
      \standardBox{1}{0}{$p_{#3}$}{0}
    \fi
  }%
}

\newcommand{\opE}[1]{\TRfamily{#1}{a}{b}}

\newcommand{\TRminusbackground}{%
  \fill[gridgrey] (0,0) rectangle (1,3);
  \fill[gridgrey] (2,0) rectangle (3,3);
  \foreach \x in {1,2,3}{
    \draw[black,line width=0.5pt] (\x,0)--(\x,3);
  }
}

\newcommand{\opTprime}{%
  \TRpicture{%
    \TRminusbackground
    \standardBox{1}{2}{$p'_a$}{0}
    \standardBox{2}{1}{$p'_b$}{0}
    \standardBox{1}{0}{$p'_c$}{0}
  }%
}

\newcommand{\opEprime}[1]{%
  \TRpicture{%
    \TRminusbackground
    \ifcase#1\relax
    \or
      \standardBox{1}{2}{$p'_d$}{0}
      \standardBox{2}{0}{$p'_e$}{0}
    \or
      \standardBox{1}{2}{$p'_d$}{0}
      \standardBox{1}{0}{$p'_e$}{0}
      \specialBox{2}{1}{1}
    \or
      \standardBox{2}{2}{$p'_d$}{0}
      \standardBox{1}{0}{$p'_e$}{0}
    \fi
  }%
}

\newcommand{\opBare}[1]{\vcenter{\hbox{\begin{tikzpicture}[
 x=0.55cm,y=0.55cm,line width=0.45pt,line cap=round,line join=round]#1\end{tikzpicture}}}}
\newcommand{\opCross}{\opBare{\draw (-0.55,-0.55)--(0.55,0.55) (-0.55,0.55)--(0.55,-0.55);}}
\newcommand{\opMono}{\opBare{\draw (-0.55,-0.5)--(0,0)
 .. controls (0.7,0.65) and (0.6,1.1) .. (0,1.1)
 .. controls (-0.6,1.1) and (-0.7,0.65) .. (0,0)--(0.55,-0.5);}}
\newcommand{\opBigon}{\opBare{
 \draw (-1,-0.4)--(-0.65,0).. controls (-0.25,0.55) and (0.25,0.55) .. (0.65,0)--(1,-0.4);
 \draw (-1,0.4)--(-0.65,0).. controls (-0.25,-0.55) and (0.25,-0.55) .. (0.65,0)--(1,0.4);}}
\newcommand{\opOnPicture}[2]{\vcenter{\hbox{\begin{tikzpicture}[
 x=0.42cm,y=0.42cm,line width=0.45pt,line cap=round,line join=round]
 \path[use as bounding box] (-#1,-1.08) rectangle (#1,1.08);
 #2\end{tikzpicture}}}}
\newcommand{\opSchematicTriangle}{\opOnPicture{1.05}{
 \draw (-1,1)--(1,-1) (-1,-1)--(1,1) (-0.5,-1)--(-0.5,1);
}}
\newcommand{\opOnThree}[1]{\opOnPicture{1.05}{
  \ifcase#1\relax
  \or % 1: three vertical strands
    \draw (-1,1)--(-1,-1) (0,1)--(0,-1) (1,1)--(1,-1);
  \or % 2
    \draw (-1,1)--(-1,-1);
    \draw (0,1) .. controls (0,0) and (1,0) .. (1,1);
    \draw (0,-1) .. controls (0,0) and (1,0) .. (1,-1);
  \or % 3
    \draw (1,1)--(1,-1);
    \draw (-1,1) .. controls (-1,0) and (0,0) .. (0,1);
    \draw (-1,-1) .. controls (-1,0) and (0,0) .. (0,-1);
  \or % 4
    \draw (-1,-1)--(1,1);
    \draw (-1,1) .. controls (-1,0) and (0,0) .. (0,1);
    \draw (0,-1) .. controls (0,0) and (1,0) .. (1,-1);
  \or % 5
    \draw (-1,1)--(1,-1);
    \draw (0,1) .. controls (0,0) and (1,0) .. (1,1);
    \draw (-1,-1) .. controls (-1,0) and (0,0) .. (0,-1);
  \or % 6
    \draw (-1,1)--(0,-1) (0,1)--(-1,-1) (1,1)--(1,-1);
  \or % 7
    \draw (-1,1)--(-1,-1) (0,1)--(1,-1) (1,1)--(0,-1);
  \or % 8
    \draw (-1,1) .. controls (-1,0.1) and (1,0.1) .. (1,1);
    \draw (-1,-1) .. controls (-1,0) and (0,0) .. (0,-1);
    \draw (0,1)--(1,-1);
  \or % 9
    \draw (-1,1) .. controls (-1,0) and (0,0) .. (0,1);
    \draw (-1,-1) .. controls (-1,-0.1) and (1,-0.1) .. (1,-1);
    \draw (1,1)--(0,-1);
  \or % 10
    \draw (-1,1) .. controls (-1,0.1) and (1,0.1) .. (1,1);
    \draw (0,-1) .. controls (0,0) and (1,0) .. (1,-1);
    \draw (0,1)--(-1,-1);
  \or % 11
    \draw (0,1) .. controls (0,0) and (1,0) .. (1,1);
    \draw (-1,-1) .. controls (-1,-0.1) and (1,-0.1) .. (1,-1);
    \draw (-1,1)--(0,-1);
  \or % 12
    \draw (-1,1)--(1,-1) (0,1)--(-1,-1) (1,1)--(0,-1);
  \or % 13
    \draw (-1,1) .. controls (-1,0.1) and (1,0.1) .. (1,1);
    \draw (-1,-1) .. controls (-1,-0.1) and (1,-0.1) .. (1,-1);
    \draw (0,1)--(0,-1);
  \or % 14
    \draw (-1,1)--(0,-1) (0,1)--(1,-1) (1,1)--(-1,-1);
  \fi
}}

\newcommand{\qPanel}[1]{%
  \vcenter{\hbox{%
    \begin{tikzpicture}[
      x=.48cm,
      y=.48cm,
      line width=.55pt
    ]
      \path[use as bounding box]
        (-.25,-.5) rectangle (3.75,3.5);

      \fill[gridgrey]
        (-.25,-.5) rectangle (3.75,3.5);

      #1
    \end{tikzpicture}%
  }}%
}

\newcommand{\qBox}[4]{%
  \node[
    draw,
    fill=white,
    minimum width=.56cm,
    minimum height=.40cm,
    inner sep=1pt
  ] at (#1,#2) {%
    \rotatebox{#4}{$\scriptstyle #3$}%
  };
}

\newcommand{\qLeft}{%
  \draw[fill=white]
    (1,0) -- (1,3)
    .. controls (1,3.55) and (.2,3.55) .. (.2,3)
    -- (.2,0)
    .. controls (.2,-.55) and (1,-.55) .. (1,0)
    -- cycle;
}

\newcommand{\qRight}{%
  \draw[fill=white]
    (2,0) -- (2,3)
    .. controls (2,3.55) and (3,3.55) .. (3,3)
    -- (3,0)
    .. controls (3,-.55) and (2,-.55) .. (2,0)
    -- cycle;
}

\newcommand{\qClosedT}{%
  \qPanel{%
    \qLeft
    \qRight

    \qBox{1.5}{2.5}{p_i}{0}
    \qBox{2.5}{1.5}{p_j}{90}
    \qBox{1.5}{.5}{p_k}{0}
  }%
}

\newcommand{\qClosedE}[1]{%
  \qPanel{%
    \qLeft

    \ifnum#1=2
      \draw[fill=white]
        (2,2) -- (2,3)
        .. controls (2,3.55) and (3,3.55) .. (3,3)
        -- (3,2)
        .. controls (3,1.6) and (2,1.6) .. (2,2)
        -- cycle;

      \draw[fill=white]
        (2,1) -- (2,0)
        .. controls (2,-.55) and (3,-.55) .. (3,0)
        -- (3,1)
        .. controls (3,1.4) and (2,1.4) .. (2,1)
        -- cycle;

      \qBox{1.5}{2.5}{p_a}{0}
      \qBox{1.5}{.5}{p_b}{0}

    \else
      \qRight

      \ifnum#1=1
        \qBox{1.5}{2.5}{p_a}{0}
        \qBox{2.5}{.5}{p_b}{90}
      \else
        \qBox{2.5}{2.5}{p_a}{90}
        \qBox{1.5}{.5}{p_b}{0}
      \fi
    \fi
  }%
}

\newcommand{\qProductTrace}{%
  \qPanel{%
    \qLeft
    \qRight

    \qBox{1.5}{2.5}{p_s}{0}
    \qBox{1.5}{.5}{p_k}{0}
  }%
}

\newcommand{\eBox}[4]{%
  \node[
    draw,
    fill=white,
    minimum width=.43cm,
    minimum height=.34cm,
    inner sep=.4pt
  ] at (#1,#2) {%
    \rotatebox{#4}{$\scriptstyle #3$}%
  };
}

\newcommand{\eLoop}[3]{%
  \draw[fill=#3]
    (#1,.45)--(#1,3.55)
    .. controls (#1,3.95) and (#2,3.95) .. (#2,3.55)
    --(#2,.45)
    .. controls (#2,.05) and (#1,.05) .. (#1,.45)
    --cycle;
}

\newcommand{\eInput}{%
  \rPanel{%
    \fill[gridgrey] (.4,0) rectangle (1.2,4);
    \fill[gridgrey] (3.6,0) rectangle (4,4);

    \draw (.4,0)--(.4,4);
    \draw (1.2,0)--(1.2,4);
    \draw (3.6,0)--(3.6,4);

    \eLoop{2}{2.8}{gridgrey}

    \eBox{1.6}{3}{p_i}{90}
    \eBox{.8}{2}{p_\ell}{0}
    \eBox{3.2}{2}{p_j}{90}
    \eBox{1.6}{1}{p_k}{90}
  }%
}

\newcommand{\eE}[3]{%
  \rPanel{%
    \rBase{0}
    \ifnum#1=1
      \rBox{1.5}{2.9}{#2}{0}
      \rBox{2.5}{.65}{#3}{90}
    \else
      \ifnum#1=2
        \rSpecial
        \rBox{1.5}{2.9}{#2}{0}
        \rBox{1.5}{.65}{#3}{0}
      \else
        \rBox{2.5}{2.9}{#2}{90}
        \rBox{1.5}{.65}{#3}{0}
      \fi
    \fi
  }%
}

\newcommand{\eT}[3]{%
  \rPanel{%
    \rBase{0}
    \rBox{1.5}{2.9}{#1}{0}
    \rBox{2.5}{1.95}{#2}{90}
    \rBox{1.5}{.65}{#3}{0}
  }%
}

\newcommand{\eW}[3]{%
  \rPanel{%
    \rBase{0}
    \rBox{2.5}{2.9}{#1}{90}
    \rBox{1.5}{1.95}{#2}{0}
    \rBox{2.5}{.65}{#3}{90}
  }%
}

\newcommand{\eWraw}[1]{%
  \rPanel{%
    \fill[gridgrey] (.4,0) rectangle (1.2,4);
    \fill[gridgrey] (3.6,0) rectangle (4,4);

    \draw (.4,0)--(.4,4);
    \draw (1.2,0)--(1.2,4);
    \draw (3.6,0)--(3.6,4);

    \ifnum#1=1
      \draw[fill=gridgrey]
        (2,1)
        .. controls (2,.4) and (2.8,.4) .. (2.8,1)
        --(2.8,2.1)
        .. controls (2.8,2.7) and (2,2.7) .. (2,2.1)
        --cycle;

      \eBox{.8}{3}{p_a}{0}
      \eBox{1.6}{1.55}{p_b}{90}
      \eBox{3.2}{1.55}{p_j}{90}
    \else
      \draw[fill=gridgrey]
        (2,2)
        .. controls (2,1.4) and (2.8,1.4) .. (2.8,2)
        --(2.8,3.1)
        .. controls (2.8,3.7) and (2,3.7) .. (2,3.1)
        --cycle;

      \eBox{.8}{1}{p_b}{0}
      \eBox{1.6}{2.55}{p_a}{90}
      \eBox{3.2}{2.55}{p_j}{90}
    \fi
  }%
}

\newcommand{\eTraw}[1]{%
  \rPanel{%
    \fill[gridgrey] (.4,0) rectangle (2.8,4);
    \fill[gridgrey] (3.6,0) rectangle (4,4);

    \draw (.4,0)--(.4,4);
    \draw (2.8,0)--(2.8,4);
    \draw (3.6,0)--(3.6,4);

    \ifnum#1=1
      \draw[fill=white]
        (1.2,2)
        .. controls (1.2,1.4) and (2,1.4) .. (2,2)
        --(2,3)
        .. controls (2,3.6) and (1.2,3.6) .. (1.2,3)
        --cycle;

      \eBox{.8}{2.5}{p_\ell}{0}
      \eBox{2.4}{2.5}{p_a}{0}
      \eBox{3.2}{.75}{p_b}{90}
    \else
      \draw[fill=white]
        (1.2,1)
        .. controls (1.2,.4) and (2,.4) .. (2,1)
        --(2,2)
        .. controls (2,2.6) and (1.2,2.6) .. (1.2,2)
        --cycle;

      \eBox{.8}{1.5}{p_\ell}{0}
      \eBox{2.4}{1.5}{p_b}{0}
      \eBox{3.2}{3.1}{p_a}{90}
    \fi
  }%
}

\newcommand{\eTmiddle}{%
  \rPanel{%
    \fill[gridgrey] (.4,0) rectangle (2.8,4);
    \fill[gridgrey] (3.6,0) rectangle (4,4);

    \draw (.4,0)--(.4,4);
    \draw (2.8,0)--(2.8,4);
    \draw (3.6,0)--(3.6,4);

    \draw[fill=white]
      (1.2,.6)--(1.2,3.4)
      .. controls (1.2,4) and (2,4) .. (2,3.4)
      --(2,2.7)
      .. controls (2,2.2) and (2.8,2.2) .. (2.8,2.7)
      --(2.8,1.3)
      .. controls (2.8,1.8) and (2,1.8) .. (2,1.3)
      --(2,.6)
      .. controls (2,0) and (1.2,0) .. (1.2,.6)
      --cycle;

    \eBox{.8}{2}{p_\ell}{0}
    \eBox{2.4}{3.2}{p_p}{0}
    \eBox{2.4}{.8}{p_q}{0}
  }%
}

\newcommand{\opMarkedBox}[3]{%
  \opEditor{%
    \drawGridBackground{2}{0}{#1}
    \standardBox{1}{0}{\rotatebox{#3}{$#2$}}{0}
  }%
}
\newcommand{\rPanel}[1]{%
  \vcenter{\hbox{%
    \begin{tikzpicture}[
      x=.43cm,
      y=.43cm,
      line width=.5pt
    ]
      \path[use as bounding box]
        (0,0) rectangle (4,4);
      \fill[white]
        (0,0) rectangle (4,4);
      #1
    \end{tikzpicture}%
  }}%
}

\newcommand{\rBox}[4]{%
  \node[
    draw,
    fill=white,
    minimum width=.50cm,
    minimum height=.34cm,
    inner sep=.5pt
  ] at (#1,#2) {%
    \rotatebox{#4}{$\scriptstyle #3$}%
  };
}

\newcommand{\rBase}[1]{%
  \ifnum#1=0
    \fill[gridgrey] (1,0) rectangle (2,4);
    \fill[gridgrey] (3,0) rectangle (4,4);

    \draw (1,0)--(1,4);
    \draw (2,0)--(2,4);
    \draw (3,0)--(3,4);
  \else
    \ifnum#1=1
      \fill[gridgrey] (1,0) rectangle (2,3.4);
      \fill[gridgrey] (3,0) rectangle (4,4);

      \fill[gridgrey]
        (1,3.4)
        .. controls (1,4.05) and (2,4.05) .. (2,3.4)
        -- cycle;

      \draw
        (1,0)--(1,3.4)
        .. controls (1,4.05) and (2,4.05) .. (2,3.4)
        --(2,0);

      \draw (3,0)--(3,4);
    \else
      \fill[gridgrey] (1,0) rectangle (2,4);
      \fill[gridgrey] (3,0) rectangle (4,4);
      \fill[gridgrey] (2,3.4) rectangle (3,4);

      \fill[white]
        (2,3.4)
        .. controls (2,4.05) and (3,4.05) .. (3,3.4)
        -- cycle;

      \draw (1,0)--(1,4);

      \draw
        (2,0)--(2,3.4)
        .. controls (2,4.05) and (3,4.05) .. (3,3.4)
        --(3,0);
    \fi
  \fi
}

\newcommand{\rSpecial}{%
  \fill[gridgrey]
    (1.99,1.2) rectangle (3.01,1.9);

  \fill[white]
    (2,1.9)
    .. controls (2,1.53) and (3,1.53) .. (3,1.9)
    -- cycle;

  \fill[white]
    (2,1.2)
    .. controls (2,1.57) and (3,1.57) .. (3,1.2)
    -- cycle;

  \draw
    (2,1.9)
    .. controls (2,1.53) and (3,1.53) .. (3,1.9);

  \draw
    (2,1.2)
    .. controls (2,1.57) and (3,1.57) .. (3,1.2);
}

\newcommand{\rT}[1]{%
  \rPanel{%
    \rBase{#1}
    \rBox{1.5}{2.9}{p_i}{0}
    \rBox{2.5}{1.95}{p_j}{90}
    \rBox{1.5}{.65}{p_k}{0}
  }%
}

\newcommand{\rE}[2]{%
  \rPanel{%
    \rBase{#1}
    \ifnum#2=1
      \rBox{1.5}{2.9}{p_a}{0}
      \rBox{2.5}{.65}{p_b}{90}
    \else
      \ifnum#2=2
        \rSpecial
        \rBox{1.5}{2.9}{p_a}{0}
        \rBox{1.5}{.65}{p_b}{0}
      \else
        \rBox{2.5}{2.9}{p_a}{90}
        \rBox{1.5}{.65}{p_b}{0}
      \fi
    \fi
  }%
}

\newcommand{\rProduct}[2]{%
  \rPanel{%
    \rBase{2}
    \rBox{1.5}{2.65}{#1}{0}
    \rBox{1.5}{.9}{#2}{0}
  }%
}

\newcommand{\rOut}[1]{%
  \rPanel{%
    \rBase{2}
    \rBox{1.5}{1.6}{#1}{0}
  }%
}

\newcommand{\rV}{%
  \rPanel{%
    \rBase{0}
    \rBox{1.5}{3.5}{p_\ell}{0}
    \rBox{1.5}{2.5}{p_i}{0}
    \rBox{2.5}{1.5}{p_j}{90}
    \rBox{1.5}{.5}{p_k}{0}
  }%
}

\newcommand{\rVE}[1]{%
  \rPanel{%
    \rBase{0}
    \rBox{1.5}{3.5}{p_\ell}{0}
    \ifnum#1=1
      \rBox{1.5}{2.5}{p_a}{0}
      \rBox{2.5}{.5}{p_b}{90}
    \else
      \ifnum#1=2
        \rSpecial
        \rBox{1.5}{2.5}{p_a}{0}
        \rBox{1.5}{.5}{p_b}{0}
      \else
        \rBox{2.5}{2}{p_a}{90}
        \rBox{1.5}{.5}{p_b}{0}
      \fi
    \fi
  }%
}

\Needspace{110mm}

\newcommand{\OnFour}[1]{%
  \OnPicture{1.58}{%
    \ifcase#1\relax
    \or % 1: 0 vertices. Strand endpoints ((1, 8), (2, 7), (3, 6), (4, 5))
      \draw (-1.50000,1.00000) -- (-1.50000,-1.00000);
      \draw (-0.50000,1.00000) -- (-0.50000,-1.00000);
      \draw (0.50000,1.00000) -- (0.50000,-1.00000);
      \draw (1.50000,1.00000) -- (1.50000,-1.00000);
    \or % 2: 0 vertices. Strand endpoints ((1, 2), (3, 6), (4, 5), (7, 8))
      \draw (-1.50000,1.00000) .. controls (-1.50000,0.57333) and (-0.50000,0.57333) .. (-0.50000,1.00000);
      \draw (0.50000,1.00000) -- (0.50000,-1.00000);
      \draw (1.50000,1.00000) -- (1.50000,-1.00000);
      \draw (-0.50000,-1.00000) .. controls (-0.50000,-0.57333) and (-1.50000,-0.57333) .. (-1.50000,-1.00000);
    \or % 3: 0 vertices. Strand endpoints ((1, 2), (3, 8), (4, 5), (6, 7))
      \draw (-1.50000,1.00000) .. controls (-1.50000,0.57333) and (-0.50000,0.57333) .. (-0.50000,1.00000);
      \draw (0.50000,1.00000) -- (-1.50000,-1.00000);
      \draw (1.50000,1.00000) -- (1.50000,-1.00000);
      \draw (0.50000,-1.00000) .. controls (0.50000,-0.57333) and (-0.50000,-0.57333) .. (-0.50000,-1.00000);
    \or % 4: 0 vertices. Strand endpoints ((1, 2), (3, 8), (4, 7), (5, 6))
      \draw (-1.50000,1.00000) .. controls (-1.50000,0.57333) and (-0.50000,0.57333) .. (-0.50000,1.00000);
      \draw (0.50000,1.00000) -- (-1.50000,-1.00000);
      \draw (1.50000,1.00000) -- (-0.50000,-1.00000);
      \draw (1.50000,-1.00000) .. controls (1.50000,-0.57333) and (0.50000,-0.57333) .. (0.50000,-1.00000);
    \or % 5: 0 vertices. Strand endpoints ((1, 6), (2, 3), (4, 5), (7, 8))
      \draw (-1.50000,1.00000) -- (0.50000,-1.00000);
      \draw (-0.50000,1.00000) .. controls (-0.50000,0.57333) and (0.50000,0.57333) .. (0.50000,1.00000);
      \draw (1.50000,1.00000) -- (1.50000,-1.00000);
      \draw (-0.50000,-1.00000) .. controls (-0.50000,-0.57333) and (-1.50000,-0.57333) .. (-1.50000,-1.00000);
    \or % 6: 0 vertices. Strand endpoints ((1, 6), (2, 5), (3, 4), (7, 8))
      \draw (-1.50000,1.00000) -- (0.50000,-1.00000);
      \draw (-0.50000,1.00000) -- (1.50000,-1.00000);
      \draw (0.50000,1.00000) .. controls (0.50000,0.57333) and (1.50000,0.57333) .. (1.50000,1.00000);
      \draw (-0.50000,-1.00000) .. controls (-0.50000,-0.57333) and (-1.50000,-0.57333) .. (-1.50000,-1.00000);
    \or % 7: 0 vertices. Strand endpoints ((1, 8), (2, 3), (4, 5), (6, 7))
      \draw (-1.50000,1.00000) -- (-1.50000,-1.00000);
      \draw (-0.50000,1.00000) .. controls (-0.50000,0.57333) and (0.50000,0.57333) .. (0.50000,1.00000);
      \draw (1.50000,1.00000) -- (1.50000,-1.00000);
      \draw (0.50000,-1.00000) .. controls (0.50000,-0.57333) and (-0.50000,-0.57333) .. (-0.50000,-1.00000);
    \or % 8: 0 vertices. Strand endpoints ((1, 8), (2, 3), (4, 7), (5, 6))
      \draw (-1.50000,1.00000) -- (-1.50000,-1.00000);
      \draw (-0.50000,1.00000) .. controls (-0.50000,0.57333) and (0.50000,0.57333) .. (0.50000,1.00000);
      \draw (1.50000,1.00000) -- (-0.50000,-1.00000);
      \draw (1.50000,-1.00000) .. controls (1.50000,-0.57333) and (0.50000,-0.57333) .. (0.50000,-1.00000);
    \or % 9: 0 vertices. Strand endpoints ((1, 8), (2, 5), (3, 4), (6, 7))
      \draw (-1.50000,1.00000) -- (-1.50000,-1.00000);
      \draw (-0.50000,1.00000) -- (1.50000,-1.00000);
      \draw (0.50000,1.00000) .. controls (0.50000,0.57333) and (1.50000,0.57333) .. (1.50000,1.00000);
      \draw (0.50000,-1.00000) .. controls (0.50000,-0.57333) and (-0.50000,-0.57333) .. (-0.50000,-1.00000);
    \or % 10: 0 vertices. Strand endpoints ((1, 8), (2, 7), (3, 4), (5, 6))
      \draw (-1.50000,1.00000) -- (-1.50000,-1.00000);
      \draw (-0.50000,1.00000) -- (-0.50000,-1.00000);
      \draw (0.50000,1.00000) .. controls (0.50000,0.57333) and (1.50000,0.57333) .. (1.50000,1.00000);
      \draw (1.50000,-1.00000) .. controls (1.50000,-0.57333) and (0.50000,-0.57333) .. (0.50000,-1.00000);
    \or % 11: 0 vertices. Strand endpoints ((1, 2), (3, 4), (5, 6), (7, 8))
      \draw (-1.50000,1.00000) .. controls (-1.50000,0.57333) and (-0.50000,0.57333) .. (-0.50000,1.00000);
      \draw (0.50000,1.00000) .. controls (0.50000,0.57333) and (1.50000,0.57333) .. (1.50000,1.00000);
      \draw (1.50000,-1.00000) .. controls (1.50000,-0.57333) and (0.50000,-0.57333) .. (0.50000,-1.00000);
      \draw (-0.50000,-1.00000) .. controls (-0.50000,-0.57333) and (-1.50000,-0.57333) .. (-1.50000,-1.00000);
    \or % 12: 0 vertices. Strand endpoints ((1, 2), (3, 4), (5, 8), (6, 7))
      \draw (-1.50000,1.00000) .. controls (-1.50000,0.57333) and (-0.50000,0.57333) .. (-0.50000,1.00000);
      \draw (0.50000,1.00000) .. controls (0.50000,0.57333) and (1.50000,0.57333) .. (1.50000,1.00000);
      \draw (1.50000,-1.00000) .. controls (1.50000,-0.26099) and (-1.50000,-0.26099) .. (-1.50000,-1.00000);
      \draw (0.50000,-1.00000) .. controls (0.50000,-0.57333) and (-0.50000,-0.57333) .. (-0.50000,-1.00000);
    \or % 13: 0 vertices. Strand endpoints ((1, 4), (2, 3), (5, 6), (7, 8))
      \draw (-1.50000,1.00000) .. controls (-1.50000,0.26099) and (1.50000,0.26099) .. (1.50000,1.00000);
      \draw (-0.50000,1.00000) .. controls (-0.50000,0.57333) and (0.50000,0.57333) .. (0.50000,1.00000);
      \draw (1.50000,-1.00000) .. controls (1.50000,-0.57333) and (0.50000,-0.57333) .. (0.50000,-1.00000);
      \draw (-0.50000,-1.00000) .. controls (-0.50000,-0.57333) and (-1.50000,-0.57333) .. (-1.50000,-1.00000);
    \or % 14: 0 vertices. Strand endpoints ((1, 4), (2, 3), (5, 8), (6, 7))
      \draw (-1.50000,1.00000) .. controls (-1.50000,0.26099) and (1.50000,0.26099) .. (1.50000,1.00000);
      \draw (-0.50000,1.00000) .. controls (-0.50000,0.57333) and (0.50000,0.57333) .. (0.50000,1.00000);
      \draw (1.50000,-1.00000) .. controls (1.50000,-0.26099) and (-1.50000,-0.26099) .. (-1.50000,-1.00000);
      \draw (0.50000,-1.00000) .. controls (0.50000,-0.57333) and (-0.50000,-0.57333) .. (-0.50000,-1.00000);
    \or % 15: 1 vertices. Strand endpoints ((1, 7), (2, 8), (3, 6), (4, 5))
      \draw (-1.50000,1.00000) -- (-0.50000,-1.00000);
      \draw (-0.50000,1.00000) -- (-1.50000,-1.00000);
      \draw (0.50000,1.00000) -- (0.50000,-1.00000);
      \draw (1.50000,1.00000) -- (1.50000,-1.00000);
    \or % 16: 1 vertices. Strand endpoints ((1, 8), (2, 6), (3, 7), (4, 5))
      \draw (-1.50000,1.00000) -- (-1.50000,-1.00000);
      \draw (-0.50000,1.00000) -- (0.50000,-1.00000);
      \draw (0.50000,1.00000) -- (-0.50000,-1.00000);
      \draw (1.50000,1.00000) -- (1.50000,-1.00000);
    \or % 17: 1 vertices. Strand endpoints ((1, 8), (2, 7), (3, 5), (4, 6))
      \draw (-1.50000,1.00000) -- (-1.50000,-1.00000);
      \draw (-0.50000,1.00000) -- (-0.50000,-1.00000);
      \draw (0.50000,1.00000) -- (1.50000,-1.00000);
      \draw (1.50000,1.00000) -- (0.50000,-1.00000);
    \or % 18: 1 vertices. Strand endpoints ((1, 2), (3, 5), (4, 6), (7, 8))
      \draw (-1.50000,1.00000) .. controls (-1.50000,0.57333) and (-0.50000,0.57333) .. (-0.50000,1.00000);
      \draw (0.50000,1.00000) -- (1.50000,-1.00000);
      \draw (1.50000,1.00000) -- (0.50000,-1.00000);
      \draw (-0.50000,-1.00000) .. controls (-0.50000,-0.57333) and (-1.50000,-0.57333) .. (-1.50000,-1.00000);
    \or % 19: 1 vertices. Strand endpoints ((1, 2), (3, 5), (4, 8), (6, 7))
      \draw (-1.50000,1.00000) .. controls (-1.50000,0.57333) and (-0.50000,0.57333) .. (-0.50000,1.00000);
      \draw (0.50000,1.00000) -- (1.50000,-1.00000);
      \draw (1.50000,1.00000) -- (-1.50000,-1.00000);
      \draw (0.50000,-1.00000) .. controls (0.50000,-0.57333) and (-0.50000,-0.57333) .. (-0.50000,-1.00000);
    \or % 20: 1 vertices. Strand endpoints ((1, 2), (3, 7), (4, 5), (6, 8))
      \draw (-1.50000,1.00000) .. controls (-1.50000,0.57333) and (-0.50000,0.57333) .. (-0.50000,1.00000);
      \draw (0.50000,1.00000) -- (-0.50000,-1.00000);
      \draw (1.50000,1.00000) -- (1.50000,-1.00000);
      \draw (0.50000,-1.00000) .. controls (0.50000,-0.39660) and (-1.50000,-0.39660) .. (-1.50000,-1.00000);
    \or % 21: 1 vertices. Strand endpoints ((1, 2), (3, 7), (4, 8), (5, 6))
      \draw (-1.50000,1.00000) .. controls (-1.50000,0.57333) and (-0.50000,0.57333) .. (-0.50000,1.00000);
      \draw (0.50000,1.00000) -- (-0.50000,-1.00000);
      \draw (1.50000,1.00000) -- (-1.50000,-1.00000);
      \draw (1.50000,-1.00000) .. controls (1.50000,-0.57333) and (0.50000,-0.57333) .. (0.50000,-1.00000);
    \or % 22: 1 vertices. Strand endpoints ((1, 2), (3, 8), (4, 6), (5, 7))
      \draw (-1.50000,1.00000) .. controls (-1.50000,0.57333) and (-0.50000,0.57333) .. (-0.50000,1.00000);
      \draw (0.50000,1.00000) -- (-1.50000,-1.00000);
      \draw (1.50000,1.00000) -- (0.50000,-1.00000);
      \draw (1.50000,-1.00000) .. controls (1.50000,-0.39660) and (-0.50000,-0.39660) .. (-0.50000,-1.00000);
    \or % 23: 1 vertices. Strand endpoints ((1, 3), (2, 6), (4, 5), (7, 8))
      \draw (-1.50000,1.00000) .. controls (-1.50000,0.39660) and (0.50000,0.39660) .. (0.50000,1.00000);
      \draw (-0.50000,1.00000) -- (0.50000,-1.00000);
      \draw (1.50000,1.00000) -- (1.50000,-1.00000);
      \draw (-0.50000,-1.00000) .. controls (-0.50000,-0.57333) and (-1.50000,-0.57333) .. (-1.50000,-1.00000);
    \or % 24: 1 vertices. Strand endpoints ((1, 3), (2, 8), (4, 5), (6, 7))
      \draw (-1.50000,1.00000) .. controls (-1.50000,0.39660) and (0.50000,0.39660) .. (0.50000,1.00000);
      \draw (-0.50000,1.00000) -- (-1.50000,-1.00000);
      \draw (1.50000,1.00000) -- (1.50000,-1.00000);
      \draw (0.50000,-1.00000) .. controls (0.50000,-0.57333) and (-0.50000,-0.57333) .. (-0.50000,-1.00000);
    \or % 25: 1 vertices. Strand endpoints ((1, 3), (2, 8), (4, 7), (5, 6))
      \draw (-1.50000,1.00000) .. controls (-1.50000,0.39660) and (0.50000,0.39660) .. (0.50000,1.00000);
      \draw (-0.50000,1.00000) -- (-1.50000,-1.00000);
      \draw (1.50000,1.00000) -- (-0.50000,-1.00000);
      \draw (1.50000,-1.00000) .. controls (1.50000,-0.57333) and (0.50000,-0.57333) .. (0.50000,-1.00000);
    \or % 26: 1 vertices. Strand endpoints ((1, 5), (2, 3), (4, 6), (7, 8))
      \draw (-1.50000,1.00000) -- (1.50000,-1.00000);
      \draw (-0.50000,1.00000) .. controls (-0.50000,0.57333) and (0.50000,0.57333) .. (0.50000,1.00000);
      \draw (1.50000,1.00000) -- (0.50000,-1.00000);
      \draw (-0.50000,-1.00000) .. controls (-0.50000,-0.57333) and (-1.50000,-0.57333) .. (-1.50000,-1.00000);
    \or % 27: 1 vertices. Strand endpoints ((1, 5), (2, 3), (4, 8), (6, 7))
      \draw (-1.50000,1.00000) -- (1.50000,-1.00000);
      \draw (-0.50000,1.00000) .. controls (-0.50000,0.57333) and (0.50000,0.57333) .. (0.50000,1.00000);
      \draw (1.50000,1.00000) -- (-1.50000,-1.00000);
      \draw (0.50000,-1.00000) .. controls (0.50000,-0.57333) and (-0.50000,-0.57333) .. (-0.50000,-1.00000);
    \or % 28: 1 vertices. Strand endpoints ((1, 5), (2, 6), (3, 4), (7, 8))
      \draw (-1.50000,1.00000) -- (1.50000,-1.00000);
      \draw (-0.50000,1.00000) -- (0.50000,-1.00000);
      \draw (0.50000,1.00000) .. controls (0.50000,0.57333) and (1.50000,0.57333) .. (1.50000,1.00000);
      \draw (-0.50000,-1.00000) .. controls (-0.50000,-0.57333) and (-1.50000,-0.57333) .. (-1.50000,-1.00000);
    \or % 29: 1 vertices. Strand endpoints ((1, 5), (2, 8), (3, 4), (6, 7))
      \draw (-1.50000,1.00000) -- (1.50000,-1.00000);
      \draw (-0.50000,1.00000) -- (-1.50000,-1.00000);
      \draw (0.50000,1.00000) .. controls (0.50000,0.57333) and (1.50000,0.57333) .. (1.50000,1.00000);
      \draw (0.50000,-1.00000) .. controls (0.50000,-0.57333) and (-0.50000,-0.57333) .. (-0.50000,-1.00000);
    \or % 30: 1 vertices. Strand endpoints ((1, 6), (2, 4), (3, 5), (7, 8))
      \draw (-1.50000,1.00000) -- (0.50000,-1.00000);
      \draw (-0.50000,1.00000) .. controls (-0.50000,0.39660) and (1.50000,0.39660) .. (1.50000,1.00000);
      \draw (0.50000,1.00000) -- (1.50000,-1.00000);
      \draw (-0.50000,-1.00000) .. controls (-0.50000,-0.57333) and (-1.50000,-0.57333) .. (-1.50000,-1.00000);
    \or % 31: 1 vertices. Strand endpoints ((1, 7), (2, 3), (4, 5), (6, 8))
      \draw (-1.50000,1.00000) -- (-0.50000,-1.00000);
      \draw (-0.50000,1.00000) .. controls (-0.50000,0.57333) and (0.50000,0.57333) .. (0.50000,1.00000);
      \draw (1.50000,1.00000) -- (1.50000,-1.00000);
      \draw (0.50000,-1.00000) .. controls (0.50000,-0.39660) and (-1.50000,-0.39660) .. (-1.50000,-1.00000);
    \or % 32: 1 vertices. Strand endpoints ((1, 7), (2, 3), (4, 8), (5, 6))
      \draw (-1.50000,1.00000) -- (-0.50000,-1.00000);
      \draw (-0.50000,1.00000) .. controls (-0.50000,0.57333) and (0.50000,0.57333) .. (0.50000,1.00000);
      \draw (1.50000,1.00000) -- (-1.50000,-1.00000);
      \draw (1.50000,-1.00000) .. controls (1.50000,-0.57333) and (0.50000,-0.57333) .. (0.50000,-1.00000);
    \or % 33: 1 vertices. Strand endpoints ((1, 7), (2, 5), (3, 4), (6, 8))
      \draw (-1.50000,1.00000) -- (-0.50000,-1.00000);
      \draw (-0.50000,1.00000) -- (1.50000,-1.00000);
      \draw (0.50000,1.00000) .. controls (0.50000,0.57333) and (1.50000,0.57333) .. (1.50000,1.00000);
      \draw (0.50000,-1.00000) .. controls (0.50000,-0.39660) and (-1.50000,-0.39660) .. (-1.50000,-1.00000);
    \or % 34: 1 vertices. Strand endpoints ((1, 7), (2, 8), (3, 4), (5, 6))
      \draw (-1.50000,1.00000) -- (-0.50000,-1.00000);
      \draw (-0.50000,1.00000) -- (-1.50000,-1.00000);
      \draw (0.50000,1.00000) .. controls (0.50000,0.57333) and (1.50000,0.57333) .. (1.50000,1.00000);
      \draw (1.50000,-1.00000) .. controls (1.50000,-0.57333) and (0.50000,-0.57333) .. (0.50000,-1.00000);
    \or % 35: 1 vertices. Strand endpoints ((1, 8), (2, 3), (4, 6), (5, 7))
      \draw (-1.50000,1.00000) -- (-1.50000,-1.00000);
      \draw (-0.50000,1.00000) .. controls (-0.50000,0.57333) and (0.50000,0.57333) .. (0.50000,1.00000);
      \draw (1.50000,1.00000) -- (0.50000,-1.00000);
      \draw (1.50000,-1.00000) .. controls (1.50000,-0.39660) and (-0.50000,-0.39660) .. (-0.50000,-1.00000);
    \or % 36: 1 vertices. Strand endpoints ((1, 8), (2, 4), (3, 5), (6, 7))
      \draw (-1.50000,1.00000) -- (-1.50000,-1.00000);
      \draw (-0.50000,1.00000) .. controls (-0.50000,0.39660) and (1.50000,0.39660) .. (1.50000,1.00000);
      \draw (0.50000,1.00000) -- (1.50000,-1.00000);
      \draw (0.50000,-1.00000) .. controls (0.50000,-0.57333) and (-0.50000,-0.57333) .. (-0.50000,-1.00000);
    \or % 37: 1 vertices. Strand endpoints ((1, 8), (2, 4), (3, 7), (5, 6))
      \draw (-1.50000,1.00000) -- (-1.50000,-1.00000);
      \draw (-0.50000,1.00000) .. controls (-0.50000,0.39660) and (1.50000,0.39660) .. (1.50000,1.00000);
      \draw (0.50000,1.00000) -- (-0.50000,-1.00000);
      \draw (1.50000,-1.00000) .. controls (1.50000,-0.57333) and (0.50000,-0.57333) .. (0.50000,-1.00000);
    \or % 38: 1 vertices. Strand endpoints ((1, 8), (2, 6), (3, 4), (5, 7))
      \draw (-1.50000,1.00000) -- (-1.50000,-1.00000);
      \draw (-0.50000,1.00000) -- (0.50000,-1.00000);
      \draw (0.50000,1.00000) .. controls (0.50000,0.57333) and (1.50000,0.57333) .. (1.50000,1.00000);
      \draw (1.50000,-1.00000) .. controls (1.50000,-0.39660) and (-0.50000,-0.39660) .. (-0.50000,-1.00000);
    \or % 39: 1 vertices. Strand endpoints ((1, 2), (3, 4), (5, 7), (6, 8))
      \draw (-1.50000,1.00000) .. controls (-1.50000,0.57333) and (-0.50000,0.57333) .. (-0.50000,1.00000);
      \draw (0.50000,1.00000) .. controls (0.50000,0.57333) and (1.50000,0.57333) .. (1.50000,1.00000);
      \draw (1.50000,-1.00000) .. controls (1.50000,-0.39660) and (-0.50000,-0.39660) .. (-0.50000,-1.00000);
      \draw (0.50000,-1.00000) .. controls (0.50000,-0.39660) and (-1.50000,-0.39660) .. (-1.50000,-1.00000);
    \or % 40: 1 vertices. Strand endpoints ((1, 3), (2, 4), (5, 6), (7, 8))
      \draw (-1.50000,1.00000) .. controls (-1.50000,0.39660) and (0.50000,0.39660) .. (0.50000,1.00000);
      \draw (-0.50000,1.00000) .. controls (-0.50000,0.39660) and (1.50000,0.39660) .. (1.50000,1.00000);
      \draw (1.50000,-1.00000) .. controls (1.50000,-0.57333) and (0.50000,-0.57333) .. (0.50000,-1.00000);
      \draw (-0.50000,-1.00000) .. controls (-0.50000,-0.57333) and (-1.50000,-0.57333) .. (-1.50000,-1.00000);
    \or % 41: 1 vertices. Strand endpoints ((1, 3), (2, 4), (5, 8), (6, 7))
      \draw (-1.50000,1.00000) .. controls (-1.50000,0.39660) and (0.50000,0.39660) .. (0.50000,1.00000);
      \draw (-0.50000,1.00000) .. controls (-0.50000,0.39660) and (1.50000,0.39660) .. (1.50000,1.00000);
      \draw (1.50000,-1.00000) .. controls (1.50000,-0.26099) and (-1.50000,-0.26099) .. (-1.50000,-1.00000);
      \draw (0.50000,-1.00000) .. controls (0.50000,-0.57333) and (-0.50000,-0.57333) .. (-0.50000,-1.00000);
    \or % 42: 1 vertices. Strand endpoints ((1, 4), (2, 3), (5, 7), (6, 8))
      \draw (-1.50000,1.00000) .. controls (-1.50000,0.26099) and (1.50000,0.26099) .. (1.50000,1.00000);
      \draw (-0.50000,1.00000) .. controls (-0.50000,0.57333) and (0.50000,0.57333) .. (0.50000,1.00000);
      \draw (1.50000,-1.00000) .. controls (1.50000,-0.39660) and (-0.50000,-0.39660) .. (-0.50000,-1.00000);
      \draw (0.50000,-1.00000) .. controls (0.50000,-0.39660) and (-1.50000,-0.39660) .. (-1.50000,-1.00000);
    \or % 43: 2 vertices. Strand endpoints ((1, 6), (2, 8), (3, 7), (4, 5))
      \draw (-1.50000,1.00000) -- (0.50000,-1.00000);
      \draw (-0.50000,1.00000) -- (-1.50000,-1.00000);
      \draw (0.50000,1.00000) -- (-0.50000,-1.00000);
      \draw (1.50000,1.00000) -- (1.50000,-1.00000);
    \or % 44: 2 vertices. Strand endpoints ((1, 7), (2, 6), (3, 8), (4, 5))
      \draw (-1.50000,1.00000) -- (-0.50000,-1.00000);
      \draw (-0.50000,1.00000) -- (0.50000,-1.00000);
      \draw (0.50000,1.00000) -- (-1.50000,-1.00000);
      \draw (1.50000,1.00000) -- (1.50000,-1.00000);
    \or % 45: 2 vertices. Strand endpoints ((1, 7), (2, 8), (3, 5), (4, 6))
      \draw (-1.50000,1.00000) -- (-0.50000,-1.00000);
      \draw (-0.50000,1.00000) -- (-1.50000,-1.00000);
      \draw (0.50000,1.00000) -- (1.50000,-1.00000);
      \draw (1.50000,1.00000) -- (0.50000,-1.00000);
    \or % 46: 2 vertices. Strand endpoints ((1, 8), (2, 5), (3, 7), (4, 6))
      \draw (-1.50000,1.00000) -- (-1.50000,-1.00000);
      \draw (-0.50000,1.00000) -- (1.50000,-1.00000);
      \draw (0.50000,1.00000) -- (-0.50000,-1.00000);
      \draw (1.50000,1.00000) -- (0.50000,-1.00000);
    \or % 47: 2 vertices. Strand endpoints ((1, 8), (2, 6), (3, 5), (4, 7))
      \draw (-1.50000,1.00000) -- (-1.50000,-1.00000);
      \draw (-0.50000,1.00000) -- (0.50000,-1.00000);
      \draw (0.50000,1.00000) -- (1.50000,-1.00000);
      \draw (1.50000,1.00000) -- (-0.50000,-1.00000);
    \or % 48: 2 vertices. Strand endpoints ((1, 2), (3, 5), (4, 7), (6, 8))
      \draw (-1.50000,1.00000) .. controls (-1.50000,0.57333) and (-0.50000,0.57333) .. (-0.50000,1.00000);
      \draw (0.50000,1.00000) -- (1.50000,-1.00000);
      \draw (1.50000,1.00000) -- (-0.50000,-1.00000);
      \draw (0.50000,-1.00000) .. controls (0.50000,-0.39660) and (-1.50000,-0.39660) .. (-1.50000,-1.00000);
    \or % 49: 2 vertices. Strand endpoints ((1, 2), (3, 6), (4, 8), (5, 7))
      \draw (-1.50000,1.00000) .. controls (-1.50000,0.57333) and (-0.50000,0.57333) .. (-0.50000,1.00000);
      \draw (0.50000,1.00000) -- (0.50000,-1.00000);
      \draw (1.50000,1.00000) -- (-1.50000,-1.00000);
      \draw (1.50000,-1.00000) .. controls (1.50000,-0.39660) and (-0.50000,-0.39660) .. (-0.50000,-1.00000);
    \or % 50: 2 vertices. Strand endpoints ((1, 2), (3, 7), (4, 6), (5, 8))
      \draw (-1.50000,1.00000) .. controls (-1.50000,0.57333) and (-0.50000,0.57333) .. (-0.50000,1.00000);
      \draw (0.50000,1.00000) -- (-0.50000,-1.00000);
      \draw (1.50000,1.00000) -- (0.50000,-1.00000);
      \draw (1.50000,-1.00000) .. controls (1.50000,-0.26099) and (-1.50000,-0.26099) .. (-1.50000,-1.00000);
    \or % 51: 2 vertices. Strand endpoints ((1, 3), (2, 5), (4, 6), (7, 8))
      \draw (-1.50000,1.00000) .. controls (-1.50000,0.39660) and (0.50000,0.39660) .. (0.50000,1.00000);
      \draw (-0.50000,1.00000) -- (1.50000,-1.00000);
      \draw (1.50000,1.00000) -- (0.50000,-1.00000);
      \draw (-0.50000,-1.00000) .. controls (-0.50000,-0.57333) and (-1.50000,-0.57333) .. (-1.50000,-1.00000);
    \or % 52: 2 vertices. Strand endpoints ((1, 3), (2, 5), (4, 8), (6, 7))
      \draw (-1.50000,1.00000) .. controls (-1.50000,0.39660) and (0.50000,0.39660) .. (0.50000,1.00000);
      \draw (-0.50000,1.00000) -- (1.50000,-1.00000);
      \draw (1.50000,1.00000) -- (-1.50000,-1.00000);
      \draw (0.50000,-1.00000) .. controls (0.50000,-0.57333) and (-0.50000,-0.57333) .. (-0.50000,-1.00000);
    \or % 53: 2 vertices. Strand endpoints ((1, 3), (2, 7), (4, 5), (6, 8))
      \draw (-1.50000,1.00000) .. controls (-1.50000,0.39660) and (0.50000,0.39660) .. (0.50000,1.00000);
      \draw (-0.50000,1.00000) -- (-0.50000,-1.00000);
      \draw (1.50000,1.00000) -- (1.50000,-1.00000);
      \draw (0.50000,-1.00000) .. controls (0.50000,-0.39660) and (-1.50000,-0.39660) .. (-1.50000,-1.00000);
    \or % 54: 2 vertices. Strand endpoints ((1, 3), (2, 7), (4, 8), (5, 6))
      \draw (-1.50000,1.00000) .. controls (-1.50000,0.39660) and (0.50000,0.39660) .. (0.50000,1.00000);
      \draw (-0.50000,1.00000) -- (-0.50000,-1.00000);
      \draw (1.50000,1.00000) -- (-1.50000,-1.00000);
      \draw (1.50000,-1.00000) .. controls (1.50000,-0.57333) and (0.50000,-0.57333) .. (0.50000,-1.00000);
    \or % 55: 2 vertices. Strand endpoints ((1, 3), (2, 8), (4, 6), (5, 7))
      \draw (-1.50000,1.00000) .. controls (-1.50000,0.39660) and (0.50000,0.39660) .. (0.50000,1.00000);
      \draw (-0.50000,1.00000) -- (-1.50000,-1.00000);
      \draw (1.50000,1.00000) -- (0.50000,-1.00000);
      \draw (1.50000,-1.00000) .. controls (1.50000,-0.39660) and (-0.50000,-0.39660) .. (-0.50000,-1.00000);
    \or % 56: 2 vertices. Strand endpoints ((1, 4), (2, 6), (3, 5), (7, 8))
      \draw (-1.50000,1.00000) .. controls (-1.50000,0.26099) and (1.50000,0.26099) .. (1.50000,1.00000);
      \draw (-0.50000,1.00000) -- (0.50000,-1.00000);
      \draw (0.50000,1.00000) -- (1.50000,-1.00000);
      \draw (-0.50000,-1.00000) .. controls (-0.50000,-0.57333) and (-1.50000,-0.57333) .. (-1.50000,-1.00000);
    \or % 57: 2 vertices. Strand endpoints ((1, 4), (2, 8), (3, 5), (6, 7))
      \draw (-1.50000,1.00000) .. controls (-1.50000,0.26099) and (1.50000,0.26099) .. (1.50000,1.00000);
      \draw (-0.50000,1.00000) -- (-1.50000,-1.00000);
      \draw (0.50000,1.00000) -- (1.50000,-1.00000);
      \draw (0.50000,-1.00000) .. controls (0.50000,-0.57333) and (-0.50000,-0.57333) .. (-0.50000,-1.00000);
    \or % 58: 2 vertices. Strand endpoints ((1, 4), (2, 8), (3, 7), (5, 6))
      \draw (-1.50000,1.00000) .. controls (-1.50000,0.26099) and (1.50000,0.26099) .. (1.50000,1.00000);
      \draw (-0.50000,1.00000) -- (-1.50000,-1.00000);
      \draw (0.50000,1.00000) -- (-0.50000,-1.00000);
      \draw (1.50000,-1.00000) .. controls (1.50000,-0.57333) and (0.50000,-0.57333) .. (0.50000,-1.00000);
    \or % 59: 2 vertices. Strand endpoints ((1, 5), (2, 3), (4, 7), (6, 8))
      \draw (-1.50000,1.00000) -- (1.50000,-1.00000);
      \draw (-0.50000,1.00000) .. controls (-0.50000,0.57333) and (0.50000,0.57333) .. (0.50000,1.00000);
      \draw (1.50000,1.00000) -- (-0.50000,-1.00000);
      \draw (0.50000,-1.00000) .. controls (0.50000,-0.39660) and (-1.50000,-0.39660) .. (-1.50000,-1.00000);
    \or % 60: 2 vertices. Strand endpoints ((1, 5), (2, 4), (3, 6), (7, 8))
      \draw (-1.50000,1.00000) -- (1.50000,-1.00000);
      \draw (-0.50000,1.00000) .. controls (-0.50000,0.39660) and (1.50000,0.39660) .. (1.50000,1.00000);
      \draw (0.50000,1.00000) -- (0.50000,-1.00000);
      \draw (-0.50000,-1.00000) .. controls (-0.50000,-0.57333) and (-1.50000,-0.57333) .. (-1.50000,-1.00000);
    \or % 61: 2 vertices. Strand endpoints ((1, 5), (2, 4), (3, 8), (6, 7))
      \draw (-1.50000,1.00000) -- (1.50000,-1.00000);
      \draw (-0.50000,1.00000) .. controls (-0.50000,0.39660) and (1.50000,0.39660) .. (1.50000,1.00000);
      \draw (0.50000,1.00000) -- (-1.50000,-1.00000);
      \draw (0.50000,-1.00000) .. controls (0.50000,-0.57333) and (-0.50000,-0.57333) .. (-0.50000,-1.00000);
    \or % 62: 2 vertices. Strand endpoints ((1, 5), (2, 7), (3, 4), (6, 8))
      \draw (-1.50000,1.00000) -- (1.50000,-1.00000);
      \draw (-0.50000,1.00000) -- (-0.50000,-1.00000);
      \draw (0.50000,1.00000) .. controls (0.50000,0.57333) and (1.50000,0.57333) .. (1.50000,1.00000);
      \draw (0.50000,-1.00000) .. controls (0.50000,-0.39660) and (-1.50000,-0.39660) .. (-1.50000,-1.00000);
    \or % 63: 2 vertices. Strand endpoints ((1, 6), (2, 3), (4, 8), (5, 7))
      \draw (-1.50000,1.00000) -- (0.50000,-1.00000);
      \draw (-0.50000,1.00000) .. controls (-0.50000,0.57333) and (0.50000,0.57333) .. (0.50000,1.00000);
      \draw (1.50000,1.00000) -- (-1.50000,-1.00000);
      \draw (1.50000,-1.00000) .. controls (1.50000,-0.39660) and (-0.50000,-0.39660) .. (-0.50000,-1.00000);
    \or % 64: 2 vertices. Strand endpoints ((1, 6), (2, 8), (3, 4), (5, 7))
      \draw (-1.50000,1.00000) -- (0.50000,-1.00000);
      \draw (-0.50000,1.00000) -- (-1.50000,-1.00000);
      \draw (0.50000,1.00000) .. controls (0.50000,0.57333) and (1.50000,0.57333) .. (1.50000,1.00000);
      \draw (1.50000,-1.00000) .. controls (1.50000,-0.39660) and (-0.50000,-0.39660) .. (-0.50000,-1.00000);
    \or % 65: 2 vertices. Strand endpoints ((1, 7), (2, 3), (4, 6), (5, 8))
      \draw (-1.50000,1.00000) -- (-0.50000,-1.00000);
      \draw (-0.50000,1.00000) .. controls (-0.50000,0.57333) and (0.50000,0.57333) .. (0.50000,1.00000);
      \draw (1.50000,1.00000) -- (0.50000,-1.00000);
      \draw (1.50000,-1.00000) .. controls (1.50000,-0.26099) and (-1.50000,-0.26099) .. (-1.50000,-1.00000);
    \or % 66: 2 vertices. Strand endpoints ((1, 7), (2, 4), (3, 5), (6, 8))
      \draw (-1.50000,1.00000) -- (-0.50000,-1.00000);
      \draw (-0.50000,1.00000) .. controls (-0.50000,0.39660) and (1.50000,0.39660) .. (1.50000,1.00000);
      \draw (0.50000,1.00000) -- (1.50000,-1.00000);
      \draw (0.50000,-1.00000) .. controls (0.50000,-0.39660) and (-1.50000,-0.39660) .. (-1.50000,-1.00000);
    \or % 67: 2 vertices. Strand endpoints ((1, 7), (2, 4), (3, 8), (5, 6))
      \draw (-1.50000,1.00000) -- (-0.50000,-1.00000);
      \draw (-0.50000,1.00000) .. controls (-0.50000,0.39660) and (1.50000,0.39660) .. (1.50000,1.00000);
      \draw (0.50000,1.00000) -- (-1.50000,-1.00000);
      \draw (1.50000,-1.00000) .. controls (1.50000,-0.57333) and (0.50000,-0.57333) .. (0.50000,-1.00000);
    \or % 68: 2 vertices. Strand endpoints ((1, 7), (2, 6), (3, 4), (5, 8))
      \draw (-1.50000,1.00000) -- (-0.50000,-1.00000);
      \draw (-0.50000,1.00000) -- (0.50000,-1.00000);
      \draw (0.50000,1.00000) .. controls (0.50000,0.57333) and (1.50000,0.57333) .. (1.50000,1.00000);
      \draw (1.50000,-1.00000) .. controls (1.50000,-0.26099) and (-1.50000,-0.26099) .. (-1.50000,-1.00000);
    \or % 69: 2 vertices. Strand endpoints ((1, 8), (2, 4), (3, 6), (5, 7))
      \draw (-1.50000,1.00000) -- (-1.50000,-1.00000);
      \draw (-0.50000,1.00000) .. controls (-0.50000,0.39660) and (1.50000,0.39660) .. (1.50000,1.00000);
      \draw (0.50000,1.00000) -- (0.50000,-1.00000);
      \draw (1.50000,-1.00000) .. controls (1.50000,-0.39660) and (-0.50000,-0.39660) .. (-0.50000,-1.00000);
    \or % 70: 2 vertices. Strand endpoints ((1, 3), (2, 4), (5, 7), (6, 8))
      \draw (-1.50000,1.00000) .. controls (-1.50000,0.39660) and (0.50000,0.39660) .. (0.50000,1.00000);
      \draw (-0.50000,1.00000) .. controls (-0.50000,0.39660) and (1.50000,0.39660) .. (1.50000,1.00000);
      \draw (1.50000,-1.00000) .. controls (1.50000,-0.39660) and (-0.50000,-0.39660) .. (-0.50000,-1.00000);
      \draw (0.50000,-1.00000) .. controls (0.50000,-0.39660) and (-1.50000,-0.39660) .. (-1.50000,-1.00000);
    \or % 71: 3 vertices. Strand endpoints ((1, 5), (2, 8), (3, 7), (4, 6))
      \draw (-1.50000,1.00000) -- (1.50000,-1.00000);
      \draw (-0.50000,1.00000) -- (-1.50000,-1.00000);
      \draw (0.50000,1.00000) -- (-0.50000,-1.00000);
      \draw (1.50000,1.00000) -- (0.50000,-1.00000);
    \or % 72: 3 vertices. Strand endpoints ((1, 6), (2, 8), (3, 5), (4, 7))
      \draw (-1.50000,1.00000) -- (0.50000,-1.00000);
      \draw (-0.50000,1.00000) -- (-1.50000,-1.00000);
      \draw (0.50000,1.00000) -- (1.50000,-1.00000);
      \draw (1.50000,1.00000) -- (-0.50000,-1.00000);
    \or % 73: 3 vertices. Strand endpoints ((1, 7), (2, 5), (3, 8), (4, 6))
      \draw (-1.50000,1.00000) -- (-0.50000,-1.00000);
      \draw (-0.50000,1.00000) -- (1.50000,-1.00000);
      \draw (0.50000,1.00000) -- (-1.50000,-1.00000);
      \draw (1.50000,1.00000) -- (0.50000,-1.00000);
    \or % 74: 3 vertices. Strand endpoints ((1, 7), (2, 6), (3, 5), (4, 8))
      \draw (-1.50000,1.00000) -- (-0.50000,-1.00000);
      \draw (-0.50000,1.00000) -- (0.50000,-1.00000);
      \draw (0.50000,1.00000) -- (1.50000,-1.00000);
      \draw (1.50000,1.00000) -- (-1.50000,-1.00000);
    \or % 75: 3 vertices. Strand endpoints ((1, 3), (2, 5), (4, 7), (6, 8))
      \draw (-1.50000,1.00000) .. controls (-1.50000,0.39660) and (0.50000,0.39660) .. (0.50000,1.00000);
      \draw (-0.50000,1.00000) -- (1.50000,-1.00000);
      \draw (1.50000,1.00000) -- (-0.50000,-1.00000);
      \draw (0.50000,-1.00000) .. controls (0.50000,-0.39660) and (-1.50000,-0.39660) .. (-1.50000,-1.00000);
    \or % 76: 3 vertices. Strand endpoints ((1, 3), (2, 6), (4, 8), (5, 7))
      \draw (-1.50000,1.00000) .. controls (-1.50000,0.39660) and (0.50000,0.39660) .. (0.50000,1.00000);
      \draw (-0.50000,1.00000) -- (0.50000,-1.00000);
      \draw (1.50000,1.00000) -- (-1.50000,-1.00000);
      \draw (1.50000,-1.00000) .. controls (1.50000,-0.39660) and (-0.50000,-0.39660) .. (-0.50000,-1.00000);
    \or % 77: 3 vertices. Strand endpoints ((1, 3), (2, 7), (4, 6), (5, 8))
      \draw (-1.50000,1.00000) .. controls (-1.50000,0.39660) and (0.50000,0.39660) .. (0.50000,1.00000);
      \draw (-0.50000,1.00000) -- (-0.50000,-1.00000);
      \draw (1.50000,1.00000) -- (0.50000,-1.00000);
      \draw (1.50000,-1.00000) .. controls (1.50000,-0.26099) and (-1.50000,-0.26099) .. (-1.50000,-1.00000);
    \or % 78: 3 vertices. Strand endpoints ((1, 4), (2, 7), (3, 5), (6, 8))
      \draw (-1.50000,1.00000) .. controls (-1.50000,0.26099) and (1.50000,0.26099) .. (1.50000,1.00000);
      \draw (-0.50000,1.00000) -- (-0.50000,-1.00000);
      \draw (0.50000,1.00000) -- (1.50000,-1.00000);
      \draw (0.50000,-1.00000) .. controls (0.50000,-0.39660) and (-1.50000,-0.39660) .. (-1.50000,-1.00000);
    \or % 79: 3 vertices. Strand endpoints ((1, 4), (2, 8), (3, 6), (5, 7))
      \draw (-1.50000,1.00000) .. controls (-1.50000,0.26099) and (1.50000,0.26099) .. (1.50000,1.00000);
      \draw (-0.50000,1.00000) -- (-1.50000,-1.00000);
      \draw (0.50000,1.00000) -- (0.50000,-1.00000);
      \draw (1.50000,-1.00000) .. controls (1.50000,-0.39660) and (-0.50000,-0.39660) .. (-0.50000,-1.00000);
    \or % 80: 3 vertices. Strand endpoints ((1, 5), (2, 4), (3, 7), (6, 8))
      \draw (-1.50000,1.00000) -- (1.50000,-1.00000);
      \draw (-0.50000,1.00000) .. controls (-0.50000,0.39660) and (1.50000,0.39660) .. (1.50000,1.00000);
      \draw (0.50000,1.00000) -- (-0.50000,-1.00000);
      \draw (0.50000,-1.00000) .. controls (0.50000,-0.39660) and (-1.50000,-0.39660) .. (-1.50000,-1.00000);
    \or % 81: 3 vertices. Strand endpoints ((1, 6), (2, 4), (3, 8), (5, 7))
      \draw (-1.50000,1.00000) -- (0.50000,-1.00000);
      \draw (-0.50000,1.00000) .. controls (-0.50000,0.39660) and (1.50000,0.39660) .. (1.50000,1.00000);
      \draw (0.50000,1.00000) -- (-1.50000,-1.00000);
      \draw (1.50000,-1.00000) .. controls (1.50000,-0.39660) and (-0.50000,-0.39660) .. (-0.50000,-1.00000);
    \or % 82: 3 vertices. Strand endpoints ((1, 7), (2, 4), (3, 6), (5, 8))
      \draw (-1.50000,1.00000) -- (-0.50000,-1.00000);
      \draw (-0.50000,1.00000) .. controls (-0.50000,0.39660) and (1.50000,0.39660) .. (1.50000,1.00000);
      \draw (0.50000,1.00000) -- (0.50000,-1.00000);
      \draw (1.50000,-1.00000) .. controls (1.50000,-0.26099) and (-1.50000,-0.26099) .. (-1.50000,-1.00000);
    \or % 83: quadrilateral. Boundary pairs (1,2),(3,4),(5,6),(7,8).
      \draw (-.9,.42)--(.9,.42)--(.9,-.42)--(-.9,-.42)--cycle;
      \draw (-1.5,1)--(-.9,.42)--(-.5,1);
      \draw (.5,1)--(.9,.42)--(1.5,1);
      \draw (1.5,-1)--(.9,-.42)--(.5,-1);
      \draw (-.5,-1)--(-.9,-.42)--(-1.5,-1);
    \or % 84: quadrilateral. Boundary pairs (2,3),(4,5),(6,7),(8,1).
      \draw (-1,0)--(0,.55)--(1,0)--(0,-.55)--cycle;
      \draw (-1.5,1)--(-1,0)--(-1.5,-1);
      \draw (-.5,1)--(0,.55)--(.5,1);
      \draw (1.5,1)--(1,0)--(1.5,-1);
      \draw (.5,-1)--(0,-.55)--(-.5,-1);
    \fi
  }%
}

\definecolor{TRfacelabel}{RGB}{183,25,32}

\newcommand{\TRVBowtie}{%
  \draw (-1,-1.4)--(-1,1.4) (1,-1.4)--(1,1.4);
  \draw (-1.4,1.05)--(1.4,-1.05);
  \draw (-1.4,-1.05)--(1.4,1.05);
  \node[tr face] at (-.66,0) {$A$};
  \node[tr face] at (.66,0) {$B$};
}
\newcommand{\TRVFigure}[1]{%
\begin{tikzpicture}[
  x=.70cm,y=.70cm,line width=.85pt,line cap=round,line join=round,
  tr face/.style={text=TRfacelabel,font=\normalsize,inner sep=1pt},
  tr leg/.style={draw,circle,inner sep=0pt,minimum size=4mm,font=\small}
]
\path[use as bounding box] (-2.35,-2.65) rectangle (3.5,1.75);
\TRVBowtie
\ifcase#1\relax
\or % 1: the labelled (3,3)-v configuration.
  \node[tr leg] at (-1.82,1.10) {1};
  \node[tr leg] at (-1.82,-1.10) {2};
  \node[tr leg] at (-1,-1.79) {3};
  \node[tr leg] at (1,-1.79) {4};
\or % 2: connect external legs 1 and 2; C is a bigon.
  \draw (-1.4,1.05)
    .. controls (-2.20,1.65) and (-2.20,-1.65) .. (-1.4,-1.05);
  \node[tr face] at (-1.59,0) {$C$};
\or % 3: connect external legs 2 and 3; C is a monogon.
  \draw (-1.4,-1.05)
    .. controls (-1.75,-1.3125) and (-1.87,-1.92) .. (-1.48,-1.98)
    .. controls (-1.11,-2.04) and (-1,-1.75) .. (-1,-1.4);
  \node[tr face] at (-1.39,-1.52) {$C$};
\or % 4: connect external legs 3 and 4; C is a triangle.
  \draw (-1,-1.4)
    .. controls (-1,-2.40) and (1,-2.40) .. (1,-1.4);
  \node[tr face] at (0,-1.39) {$C$};
\or % 5: fork (a); C is a bigon meeting B at one vertex.
  \draw (1,-1.4)
    .. controls (1,-1.74) and (1.56,-1.99) .. (1.95,-2.05)
    --(2.30,-2.104);
  \draw (1.4,-1.05)
    .. controls (1.90,-1.425) and (1.969,-1.74) .. (1.95,-2.05)
    --(1.928,-2.41);
  \node[tr face] at (1.40,-1.52) {$C$};
\or % 6: fork (b); C is a triangle sharing an edge with B.
  \draw (1.4,1.05)
    .. controls (2.05,1.5375) and (2.85,.80) .. (3,0)
    --(3.10,-.5333);
  \draw (1.4,-1.05)
    .. controls (2.05,-1.5375) and (2.85,-.80) .. (3,0)
    --(3.10,.5333);
  \node[tr face] at (2.02,0) {$C$};
\fi
\end{tikzpicture}%
}

\newcommand{\bBase}[1]{\ifnum#1<3 \rBase{#1}\else
 \begin{scope}[shift={(0,4)},yscale=-1]
 \ifnum#1=3\rBase{1}\else\rBase{2}\fi\end{scope}\fi}
\newcommand{\gE}[4]{\rPanel{\bBase{#1}
 \ifnum#2=1\rBox{1.5}{2.9}{#3}{0}\rBox{2.5}{.85}{#4}{90}
 \else\ifnum#2=2\rSpecial\rBox{1.5}{2.9}{#3}{0}\rBox{1.5}{.85}{#4}{0}
 \else\rBox{2.5}{2.9}{#3}{90}\rBox{1.5}{.85}{#4}{0}\fi\fi}}
\newcommand{\gW}[4]{\rPanel{\bBase{#1}
 \rBox{2.5}{3.05}{#2}{90}\rBox{1.5}{2}{#3}{0}\rBox{2.5}{.95}{#4}{90}}}
\newcommand{\bOut}[1]{\rPanel{\bBase{4}\rBox{1.5}{2}{#1}{0}}}
\newcommand{\wGone}{\rPanel{\bBase{4}
 \rBox{2.5}{3}{p_i}{90}\rBox{1.5}{1.8}{p_j}{0}}}
\newcommand{\wProduct}{\rPanel{\bBase{4}
 \rBox{2.5}{3}{p_i}{90}\rBox{2.5}{1}{p_k}{90}}}
\newcommand{\wV}{\rPanel{\rBase{0}
 \rBox{2.5}{3.5}{p_i}{90}\rBox{1.5}{2.5}{p_j}{0}
 \rBox{2.5}{1.5}{p_k}{90}\rBox{2.5}{.5}{p_\ell}{90}}}
\newcommand{\wVE}[1]{\rPanel{\rBase{0}\rBox{2.5}{.5}{p_\ell}{90}
 \ifnum#1=1\rBox{1.5}{3.1}{p_a}{0}\rBox{2.5}{2}{p_b}{90}
 \else\ifnum#1=2
 \begin{scope}[shift={(0,.9)}]\rSpecial\end{scope}
 \rBox{1.5}{3.4}{p_a}{0}\rBox{1.5}{1.8}{p_b}{0}
 \else\rBox{2.5}{3.1}{p_a}{90}\rBox{1.5}{2}{p_b}{0}\fi\fi}}

\newcommand{\localstep}[1]{\par\Needspace{6\baselineskip}\smallskip\noindent\textit{#1}\par}

\definecolor{facecol}{RGB}{164,42,53}
\definecolor{good}{RGB}{31,116,101}
\definecolor{bad}{RGB}{190,116,26}
\tikzset{face/.style={text=facecol,font=\small,inner sep=1pt},leg/.style={font=\scriptsize,fill=white,inner sep=1.5pt},vtx/.style={circle,fill=black,inner sep=1.5pt}}
\newcommand{\RVHbase}{%
\draw (-1,-1.4)--(-1,1.4) (1,-1.4)--(1,1.4);
\draw (-1.4,1.05)--(1.4,-1.05) (-1.4,-1.05)--(1.4,1.05);
\node[face] at (-.65,0) {$A$};\node[face] at (.65,0) {$B$};
}
\newcommand{\RVHlabels}{%
\node[leg,left] at (-1.4,1.05) {$1$};\node[leg,left] at (-1.4,-1.05) {$2$};
\node[leg,below] at (-1,-1.4) {$3$};\node[leg,below] at (1,-1.4) {$4$};
\node[leg,right] at (1.4,-1.05) {$5$};\node[leg,above] at (-1,1.4) {$6$};
\node[leg,above] at (1,1.4) {$7$};\node[leg,right] at (1.4,1.05) {$8$};}
\newcommand{\RVLocal}[1]{%
\begin{tikzpicture}[x=.74cm,y=.74cm,line width=.75pt,line cap=round,line join=round]
\path[use as bounding box] (-2.7,-2.55) rectangle (1.85,1.7);
\RVHbase
\ifcase#1
\RVHlabels
\or
\draw[good,thick] (-1.4,1.05)..controls(-2.4,1.65)and(-2.4,-1.65)..(-1.4,-1.05);
\node[face] at (-1.65,0) {$C$};
\or
\draw[good,thick] (-1.4,-1.05)..controls(-2,-1.5)and(-1.8,-2.2)..(-1.35,-2.0)..controls(-1,-1.85)and(-1,-1.7)..(-1,-1.4);
\node[face] at (-1.4,-1.5) {$C$};
\or
\draw[good,thick] (-1,-1.4)..controls(-1,-2.4)and(1,-2.4)..(1,-1.4);
\node[face] at (0,-1.38) {$C$};
\or
\draw[good,thick] (-1.4,1.05)..controls(-1.9,1.43)and(-2.1,.65)..(-2.3,0)--(-2.5,-.45);
\draw[good,thick] (-1.4,-1.05)..controls(-1.9,-1.43)and(-2.1,-.65)..(-2.3,0)--(-2.5,.45);
\node[vtx] at (-2.3,0) {};\node[face] at (-1.6,0) {$C$};
\or
\draw[good,thick] (-1.4,-1.05)..controls(-1.9,-1.43)and(-1.95,-1.7)..(-1.8,-2.05)--(-1.63,-2.42);
\draw[good,thick] (-1,-1.4)..controls(-1,-1.8)and(-1.45,-2)..(-1.8,-2.05)--(-2.2,-2.1);
\node[vtx] at (-1.8,-2.05) {};\node[face] at (-1.45,-1.55) {$C$};
\or
\draw[bad,thick] (-1,-1.4)..controls(-1,-1.9)and(-.5,-2.05)..(0,-2.1)--(.45,-2.22);
\draw[bad,thick] (1,-1.4)..controls(1,-1.9)and(.5,-2.05)..(0,-2.1)--(-.45,-2.22);
\node[vtx] at (0,-2.1) {};\node[face] at (0,-1.35) {$C$};
\fi
\end{tikzpicture}}
\newcommand{\RVPanel}[3]{\begin{minipage}[t]{#1\linewidth}\centering #2\par\small #3\end{minipage}}
\newcommand{\RVQuad}[1]{%
\begin{tikzpicture}[scale=.85,line width=.8pt]
\draw[gray!50,dashed] (0,0) circle(1.65);
\foreach \i/\lab in {0/1,1/2,2/3,3/4,4/5,5/8,6/7,7/6}{
\coordinate (p\i) at ({180+45*\i}:1.65);
\node[leg] at ({180+45*\i}:1.91) {$\lab$};}
\ifnum#1=0
\foreach \k/\a/\b/\col in {0/0/1/good,1/2/3/bad,2/4/5/good,3/6/7/bad}{
\coordinate (v\k) at ({202.5+90*\k}:.84);
\draw[\col,thick] (p\a)--(v\k)--(p\b);\node[vtx] at(v\k){};}
\else
\foreach \k/\a/\b in {0/1/2,1/3/4,2/5/6,3/7/0}{
\coordinate(v\k)at({247.5+90*\k}:.84);
\draw[good,thick](p\a)--(v\k)--(p\b);\node[vtx]at(v\k){};}
\fi
\draw(v0)--(v1)--(v2)--(v3)--cycle;
\end{tikzpicture}}

\newcommand{\TransportProjection}[4]{%
  \pgfmathparse{int(mod(#1+#4,2))}%
  \edef\transportParity{\pgfmathresult}%
  \node[draw=black,fill=white,thick,minimum width=.60cm,
    minimum height=.45cm,text width=.56cm,align=center,inner sep=0pt,
    font=\tiny\sffamily\bfseries,transform shape]
    at (#1+.5,#2+.5) {\ifnum\transportParity=1 $#3$%
    \else\rotatebox{90}{$#3$}\fi};
}
\newcommand{\NativeTransport}[5]{%
\begin{tikzpicture}[x=.60cm,y=.60cm,baseline={([yshift=-.5ex]current bounding box.center)}]
\drawGridBackground{3}{2}{#2}
\ifcase#1
 \TransportProjection{1}{2}{#3}{#2}
 \TransportProjection{2}{1}{#4}{#2}
 \TransportProjection{1}{0}{#5}{#2}
\or
 \TransportProjection{1}{2}{#3}{#2}
 \TransportProjection{2}{0}{#4}{#2}
\or
 \TransportProjection{1}{2}{#3}{#2}
 \specialBox{2}{1}{#2}
 \TransportProjection{1}{0}{#4}{#2}
\or
 \TransportProjection{2}{2}{#3}{#2}
 \TransportProjection{1}{0}{#4}{#2}
\fi
\end{tikzpicture}}

\newcommand{\TransportGridBackground}[3]{%
 \foreach \relX in {0,...,#1} {
  \pgfmathparse{int(mod(\relX+#3,2))}
  \ifnum\pgfmathresult=1
   \fill[gridgrey] (\relX,0) rectangle (\relX+1,#2+1);
  \fi
 }
 \foreach \x in {1,...,#1} {\draw[line width=.5pt] (\x,0)--(\x,#2+1);}
}
\renewcommand{\TransportProjection}[4]{%
 \pgfmathparse{int(mod(#1+#4,2))}\edef\transportParity{\pgfmathresult}
 \node[draw=black,fill=white,line width=.5pt,
 minimum width=.64cm,minimum height=.64cm,inner sep=0pt,outer sep=0pt,
 font=\normalsize,transform shape]
 at (#1+.5,#2+.5)
 {\ifnum\transportParity=1 $#3$\else\rotatebox{90}{$#3$}\fi};
}
\renewcommand{\NativeTransport}[5]{%
\begin{tikzpicture}[x=.50cm,y=1.05cm,baseline={([yshift=-.5ex]current bounding box.center)}]
\TransportGridBackground{3}{2}{#2}
\ifcase#1
 \TransportProjection{1}{2}{#3}{#2}
 \TransportProjection{2}{1}{#4}{#2}
 \TransportProjection{1}{0}{#5}{#2}
\or
 \TransportProjection{1}{2}{#3}{#2}
 \TransportProjection{2}{0}{#4}{#2}
\or
 \TransportProjection{1}{2}{#3}{#2}
 \specialBox{2}{1}{#2}
 \TransportProjection{1}{0}{#4}{#2}
\or
 \TransportProjection{2}{2}{#3}{#2}
 \TransportProjection{1}{0}{#4}{#2}
\fi
\end{tikzpicture}}

\begin{document}

\maketitle

\begin{abstract}
We establish a finite system of consistency equations for
triangle relation planar algebras. Within the Yang--Baxter
relation class, we characterize triangle relations by
$H^2(K;\mathbb C)=0$, where $K$ is the fusion complex
recording the support of the fusion coefficients.
Together with a reconstruction theorem, this gives an
algebraic classification procedure.

We classify bi-commutative triangle relation planar algebras
with $\dim P_{2,\pm}\leq5$. In ranks two, three, and four,
the classification yields respectively $2$, $7$, and $19$
realizable fusion-support pair strata, up to relabelling
and shading reversal. These results recover the low-rank
examples of Bisch, Jones, Liu, and Lu, and extend the
rank-four classification beyond exchange relations.
In rank five, we obtain two one-parameter families
with different duality types whose generic members
do not arise from representations of groups or
quantum groups.
\end{abstract}
\clearpage
\begingroup
\small
\setcounter{tocdepth}{2}
\tableofcontents
\endgroup
\newpage

% ============================================================
\section{Introduction}
\label{sec:introduction}

The theory of subfactors studies inclusions of von Neumann factors
and the invariants associated with them. Jones' index theorem
revealed a striking rigidity in the possible indices of
$\mathrm{II}_1$ subfactors \cite{JonesIndex}. Beyond the index,
the higher relative commutants and their operations form the
standard invariant. Constructing and classifying these invariants
is a central problem in subfactor theory.

Planar algebras, introduced by Jones, describe the standard
invariant through planar tangles \cite{JonesPA}. Multiplication,
trace, inclusion, and rotation become operations within a single
diagrammatic structure. This language also provides a setting
for presentations by generators and local skein relations.
For planar algebras generated by their $2$-box spaces, generators
can be represented by labelled four-valent vertices, and suitable
local relations allow closed diagrams to be evaluated.

Skein relations have played an important role
in classification. Bisch and Jones studied singly 2-box generated
planar algebras up to 13-dimensional 3-boxes \cite{BischJonesPA,BischJonesPAII}. 
Landau formulated the corresponding skein theory as exchange relation in \cite{Landau2002}. The Bisch-Jones classification is equivalent to classification of singly generated exchange relation planar algebras.
Furthermore,  Liu studied planar algebras with multiple 2-box generators and classified exchange-relation subfactor planar
algebras with 4-dimensional $2$-box spaces in \cite{LiuER}.
More recently, Lu and Liu characterized exchange relations by a forest condition on fusion graphs and classified exchange relation planar algebras up to 6-dimensional generators and discovered a new family of subfactor planar algebras, in addition to algebraic ones\cite{LuThesis,LuLiuER}.

In this paper, we introduce the Triangle Relation (TR) for 2-box generators. The Bisch-Jones-Liu classification of the $14$-dimensional $3$-box case
\cite{BJL2014} is equivalent to singly generated TR planar algebras.
The underlying BMW algebras were introduced by
Birman--Wenzl and Murakami
\cite{BirmanWenzl1989,Murakami1987}.
Their connections with quantum groups and subfactors
were developed by Wenzl \cite{Wenzl1990}.

Our main aim is to classify TR planar algebras with multiple generators, so that we can discover new planar algebras in the classification. 
We achieve this goal in our classification of TR planar algebras up to 5-dim 2-box generators. We discover and construct two new one-parameter families of planar algebras in our classification. 

The first difficulty is consistency. Different reduction
paths may give different evaluations of the same closed
diagram. We reduce path independence for infinitely many
diagrams to finitely many local overlap tests, expressed
as polynomial equations in the fusion-bialgebra setting.
For reconstruction, we also establish compatibility
with the identity resolution.
Even this finite system is computationally demanding.
The direct formulation involves $O(n^5)$ variables
and $O(n^7)$ indexed polynomial equations, where
$n=\dim P_{2,\pm}$.
Reducing the number of unknowns is therefore essential
for effective classification.

The main reduction comes from the topology of fusion
support, encoded by a two-dimensional simplicial complex
$K$, called the \emph{fusion complex}.
Within the YBR class, triangle relations are equivalent
to $H^2(K;\mathbb C)=0$. This criterion sieves supports
by a matrix-rank computation and gives canonical triangle
coefficients by linear algebra over $\mathbb Q$.
Trace closure and the character-table construction
reduce the numerical problem to at most $n-1$
quantum dimensions on one side.
Using root coordinates, the polynomial systems involve
only $O(n)$ variables per chart.
The prior equations make the $(3,1)$ and $(3,2)$
consistency equations redundant, leaving the common-edge
$(3,3)$-e equation.
These reductions yield an algorithm combining support
enumeration, cohomological sieving, and exact polynomial
solving with nonvanishing conditions.
Reconstruction and uniqueness turn the surviving matched
data into planar algebras, which are then compared with
known models and product constructions.

The role of the cohomological sieve parallels that
of the forest criterion in the exchange-relation
classification. Lu and Liu's 6-dimensional investigation
starts from more than $34$ billion candidate fusion
graphs and leads to a new family of subfactors
\cite{LuThesis}.
In our 5-dimensional classification, the duality types
self-dual, $(12)$, and $(12)(34)$ give respectively
$2^{20}$, $2^{13}$, and $2^{10}$ labelled single-sided
support candidates, totaling $1,057,792$ before
relabelling and cohomological sieving.
The condition $H^2=0$ filters these candidates by
linear algebra before numerical solving.
The complete classification yields two new
one-parameter families, with duality types
$\mathrm{id}$ and $(12)$.

A central goal of classification is to reach the
\emph{critical dimension}: the smallest $2$-box dimension
at which examples emerge beyond the chosen known
constructions.
For exchange relations, this threshold is $3$
in the algebraic setting, with the binary-tree example
\cite{LuLiuER}, and $6$ in the subfactor setting
\cite{LuThesis}.
For bi-commutative TR planar algebras, dimension $5$
is the first to produce new families.
Our classification thus reaches the dimension at which
new examples first emerge.

A natural next step is to classify Yang--Baxter relation
planar algebras with multiple $2$-box generators,
extending Liu's classification of the singly generated
case \cite{LiuYBR}.
Our cohomological characterization of triangle relations
suggests seeking analogous structural criteria within
this broader class.
The free-product constructions of Bisch and Jones,
arising from intermediate subfactors
\cite{BischJones1997}, provide basic models for such
classifications and already occur among our TR examples.
Skein theoretical construction and classification beyond 2-box generators have been studied in 
\cite{PetersHaagerup,BMPSExtendedHaagerup,MPSTrivalent,JonesLiuRenThurston}.

\subsection{Cohomological characterization and reconstruction}

A local reduction procedure must satisfy two requirements:
it must terminate, and its value must be independent of the
choices of reductions. We first establish finite consistency
tests for the canonical reductions of closed diagrams.
The graphical argument compares overlapping local reductions
and their complementary diagrams.
When the $2$-box space has a fusion-bialgebra structure,
these graphical tests can be written as polynomial equations
in the fusion data and triangle coefficients.
This gives the algebraic consistency system of
Theorem~\ref{cons:algthm}.

Fusion bialgebras and their Fourier analysis were
introduced by Liu, Palcoux, and Wu in the study of
unitary categorification \cite{LPW21}.
Here we use an algebraic version of this framework,
without imposing positivity, to organize the
$2$-box data.
Let $B$ be a fusion bialgebra with distinguished minimal
idempotents $\{p_i\}_{i\in I}$ and convolution coefficients
defined by
\[
p_k*p_j=\sum_i N_{kj}^{i}p_i.
\]
Introduce three copies $I_X,I_Y,I_Z$ of the label set.
The fusion complex $K(B)$ has triangular faces
\[
[i_X,j_Y,k_Z]\in K(B)^{(2)}
\quad\Longleftrightarrow\quad
N_{kj}^{i}\ne0.
\]
It records which fusion coefficients are nonzero.
Their values and the quantum dimensions remain numerical
data to be determined.

In the YBR setting, the nonzero triangle diagrams form a
basis of the $3$-box space. Resolving the ordinary identity
expresses each internal-face-free diagram as a sum of triangles.
Consequently, the coefficients of a triangle expansion
satisfy a coboundary equation
\[
d_K\alpha_\tau=e_\tau,
\qquad
d_K:C^1(K;\mathbb C)\longrightarrow C^2(K;\mathbb C),
\]
where $e_\tau$ takes value $1$ on the oriented face $\tau$
and $0$ on every other oriented face.
Since $K$ is two-dimensional, these equations are solvable
for every face exactly when $H^2(K;\mathbb C)=0$.
This gives the following characterization, proved in
Theorem~\ref{thm:cohomological-criterion}.

\begin{theorem}
Let $B$ be a YBR fusion bialgebra in the sense of
Section~\ref{sec:algebraic-consistency}, and let $K(B)$
be its fusion complex. Then
\[
B\text{ is of TR type}
\quad\Longleftrightarrow\quad
H^2(K(B);\mathbb C)=0.
\]
\end{theorem}

\begin{remark}
Liu-Lu forest criterion for exchange relations
\cite{LuThesis} admits a cohomological interpretation:
a fusion graph $L_i$ is a forest if and only if
$H^1(L_i;\mathbb C)=0$.
Thus its graph condition is equivalent to
\[
H^1(L_i;\mathbb C)=0
\qquad\text{for every fusion graph }L_i.
\]
These graphs occur as links of the fusion complex.
The exchange-relation criterion concerns first
cohomology of these links, whereas our triangle-relation
criterion concerns second cohomology of the whole complex.
\end{remark}

The coboundary equations also give a way to choose
triangle coefficients before finding a numerical realization.
When $H^2(K;\mathbb Q)=0$, the Hodge gauge selects
a canonical rational solution for each face by
linear algebra over $\mathbb Q$.
Thus a fusion complex pair $(K_+,K_-)$ with
$H^2(K_\pm;\mathbb Q)=0$ determines two native
triangle coefficient tables $c,c'$.

Trace closure then expresses the fusion coefficients in terms
of the quantum dimensions:
\[
N_{ij}^{k}=\frac{Q_{ij}^{k}(d)}{d_k},
\qquad
N'{}_{ab}^{c}=\frac{Q'{}_{ab}^{c}(d')}{d'_c},
\]
where $Q$ and $Q'$ are quadratic forms determined by $c$ and $c'$.
The remaining numerical conditions include the fusion
identities on both sides, the prescribed supports, and
a Fourier matching between the two fusion bialgebras.
The matching matrices $S,S'$ also express the
common-shading triangle coefficients in terms of the
native minus table $c'$.

Together with the fusion structure identities,
the coboundary equations imply the $(3,1)$ and
$(3,2)$-e/v consistency equations, as proved in
Proposition~\ref{prop:prior-redundancy}.
After imposing the fusion and Fourier conditions,
the remaining triangle comparison is the common-edge
$(3,3)$-e equation.
The numerical data $(N,d)$ and $(N',d')$ specify
fusion bialgebras $B_+$ and $B_-$ on the supports
$K_+$ and $K_-$, respectively.
The Fourier matrices $S,S'$ provide their matching.
The following reconstruction theorem realizes
this pair as the two $2$-box spaces of a planar algebra, as shown in Theorem~\ref{thm:construction-after-consistency}.

\begin{theorem}
Let $(K_+,K_-)$ be a fusion complex pair with
$H^2(K_\pm;\mathbb Q)=0$, and let $c,c'$ be its
canonical triangle coefficient tables.
Suppose that numerical data
$(N,N',d,d',S,S',\delta)$ satisfy the reconstruction
equations of Theorem~\ref{thm:construction-after-consistency},
have exactly the supports specified by $K_\pm$,
and have nonzero quantum dimensions and modulus,
with $S,S'$ invertible.
Let $B_+$ and $B_-$ be the fusion bialgebras
determined by $(N,d)$ and $(N',d')$, respectively.
Then there exists a bi-commutative TR planar algebra $P$
realizing these data, with
$
P_{2,+}\cong B_+,
P_{2,-}\cong B_-.
$
\end{theorem}

The reconstruction starts with abstract data.
Its proof establishes a well-defined partition function
and then constructs the planar algebra while preserving
the prescribed $2$-box structures.
Thus the cohomological condition supplies canonical triangle
relations, and the remaining numerical equations determine
whether those relations admit a planar-algebra realization.

\subsection{Classification and new families}

The construction organizes classification into four stages:
enumerating fusion supports, solving the fusion equations,
finding Fourier matchings, and imposing the remaining
consistency equations.
At fixed rank there are finitely many supports to consider.
Their numerical realizations can nevertheless form
positive-dimensional algebraic families.
The algorithm in Section~\ref{sec:algorithm} retains these
families together with the nonvanishing conditions that
specify their supports.

The low-rank calculations recover the known rank-two
and rank-three examples and give complete classifications
in ranks four and five.

\begin{theorem}
\label{thm:intro-classification}
\label{thm:intro-four}
\label{thm:intro-five}
Up to planar-algebra isomorphism and shading reversal,
the bi-commutative TR planar algebras with
$\dim P_{2,\pm}=4$ are precisely the models in
Table~\ref{tab:four}, and those with
$\dim P_{2,\pm}=5$ are precisely the models in
Table~\ref{tab:five}.
Every listed model occurs throughout its stated
parameter domain.
\end{theorem}

The classifications are proved in
Theorems~\ref{thm:four} and~\ref{thm:five}.
A principal feature of the rank-five classification
is the occurrence of \textbf{two new one-parameter families},
$\mathcal H_{\mathrm{sym}}(a)$ and
$\mathcal H_{\mathrm{pair}}(a)$.
Both are defined for $a\in\mathbb C\setminus\{0,1\}$,
with
\[
d=d'=(1,a,a,2a-2,1),
\qquad
\delta^2=4a.
\]
Their duality types are different:
all distinguished labels are self-dual in
$\mathcal H_{\mathrm{sym}}(a)$, whereas duality
exchanges labels $1$ and $2$ in
$\mathcal H_{\mathrm{pair}}(a)$.
Their generic members do not arise from representations
of groups or quantum groups.
Neither family admits a nontrivial tensor-product
or ordered free-product decomposition into lower-rank
planar algebras.
Section~\ref{sec:examples} gives their explicit
fusion rules, Fourier matrices, and specializations.

These families raise a natural question.

\begin{question}
Which members of these two families arise from
representations of Hopf algebras?
\end{question}

\paragraph{Organization of the paper.}
Section~\ref{sec:preliminaries} introduces planar algebras,
fusion bialgebras, and the graphical conventions.
Section~\ref{sec:evaluation} proves finite graphical
consistency and derives its algebraic equations.
Section~\ref{sec:fusioncomplex} develops the cohomological
criterion, canonical triangle coefficients, Fourier matching,
and reconstruction, and then presents the classification
algorithm.
Section~\ref{sec:computations} gives the low-rank
classifications and the new families.
The appendix collects explicit model formulas,
support lists, and the supplementary $4$-box diagrams.

\section{Preliminaries}
\label{sec:preliminaries}
% ============================================================

\subsection{Planar algebras and local relations}

We work with a connected, irreducible, spherical,
nondegenerate shaded planar algebra
$P=\{P_{n,+},P_{n,-}\}_{n\geq0}$ over $\C$,
with every box space finite-dimensional.  The two signs record the shading adjacent to the distinguished
interval of an input or output disc.  Planar tangles act multilinearly on
the box spaces, and isotopic tangles induce the same multilinear map.

Connectedness means $P_{0,+}=P_{0,-}=\C$, with the empty diagram
identified with $1$. Irreducibility means $\dim P_{1,+}=\dim P_{1,-}=1$.
The modulus $\delta\in\C^\times$ is the value of a vertex-free circle bounding an empty disc in
either shading. Sphericality identifies evaluations related by isotopy
on the sphere.

Closing all boundary strings defines the unnormalized Markov trace
$\operatorname{tr}_{m,\varepsilon}:P_{m,\varepsilon}\to\C$.
Our normalization is
\[
\operatorname{tr}_{0,\varepsilon}(1)=1,\qquad
\operatorname{tr}_{m,\varepsilon}(1)=\delta^m.
\]
Nondegeneracy means that the bilinear trace pairing
$(x,y)\mapsto\operatorname{tr}_{m,\varepsilon}(xy)$
is nondegenerate for every $m$ and $\varepsilon$.
For spherical planar algebras, this is equivalent to
detection of every nonzero element by a complementary
closure according to \cite[Proposition~1.29]{JonesPA}.
Under these assumptions, the box algebras are semisimple
by \cite[Corollary~1.31]{JonesPA}.
In the subfactor setting there is, in addition, a conjugate-linear
reflection $x\mapsto x^*$, and
$\operatorname{tr}_{m,\varepsilon}(y^*x)$ is a positive-definite Hermitian
inner product. Positivity and the $*$-structure are not assumed for the
algebraic realizations considered here.

On the 2-box spaces the Fourier maps are fixed as
$\F_+:P_{2,+}\longrightarrow P_{2,-},\; \F_-:P_{2,-}\longrightarrow P_{2,+}$,
both counterclockwise $90^\circ$ rotations. We also write $\F$ when its
source shading is understood. The two compositions are the half-turns on
the respective 2-box spaces and implement the duality involutions.

\begin{figure}[htbp]
\centering
\begin{tikzpicture}[
  scale=.85,
  line cap=round,
  line join=round,
  every node/.style={font=\small}
]
% Left: the original 3-box.
\begin{scope}
  \foreach \j in {0,2,4}{
    \pgfmathsetmacro{\ang}{60*\j}
    \fill[gray!20]
      (\ang:1.2)
      arc[start angle=\ang,end angle=\ang+60,radius=1.2]
      -- (\ang+60:.84)
      arc[start angle=\ang+60,end angle=\ang,radius=.84]
      -- cycle;
  }

  \foreach \j in {0,...,5}{
    \pgfmathsetmacro{\ang}{60*\j}
    \draw[line width=.7pt] (\ang:.84) -- (\ang:1.2);
    \pgfmathtruncatemacro{\lab}{\j+1}
    \node at (\ang:1.48) {$\lab$};
  }

  \draw[line width=.7pt] (0,0) circle (.84);
  \draw[line width=.7pt] (0,0) circle (1.2);
  \node[font=\Large] at (0,0) {$X$};
  \fill (30:.65) circle (1.5pt);
  \node[font=\large] at (30:1.38) {$*$};
  \node[font=\normalsize] at (0,-1.95)
    {$P_{3,\varepsilon}$};
\end{scope}

% Rotation.
\draw[->,line width=.8pt] (1.95,0) -- (4.45,0);
\node[font=\large] at (3.2,.38)
  {$\rho_{3,\varepsilon}$};
\node at (3.2,-.4)
  {$180^\circ$ counterclockwise};

% Right: the rotated 3-box.
\begin{scope}[xshift=6.4cm]
  \foreach \j in {1,3,5}{
    \pgfmathsetmacro{\ang}{60*\j}
    \fill[gray!20]
      (\ang:1.2)
      arc[start angle=\ang,end angle=\ang+60,radius=1.2]
      -- (\ang+60:.84)
      arc[start angle=\ang+60,end angle=\ang,radius=.84]
      -- cycle;
  }

  \foreach \j in {0,...,5}{
    \pgfmathsetmacro{\ang}{60*\j}
    \draw[line width=.7pt] (\ang:.84) -- (\ang:1.2);
    \pgfmathtruncatemacro{\lab}{mod(\j+3,6)+1}
    \node at (\ang:1.48) {$\lab$};
  }

  \draw[line width=.7pt] (0,0) circle (.84);
  \draw[line width=.7pt] (0,0) circle (1.2);
  \node[font=\Large,rotate=180] at (0,0) {$X$};
  \fill (210:.65) circle (1.5pt);
  \node[font=\large] at (30:1.38) {$*$};
  \node[font=\normalsize] at (0,-1.95)
    {$P_{3,-\varepsilon}$};
\end{scope}
\end{tikzpicture}
\caption{The half-turn $\rho_3$.
The numbers track the boundary endpoints, and the black dot
tracks the rotation of the input box.
The outer distinguished interval remains fixed,
while its shading changes.}
\label{fig:annular-half-turn}
\end{figure}
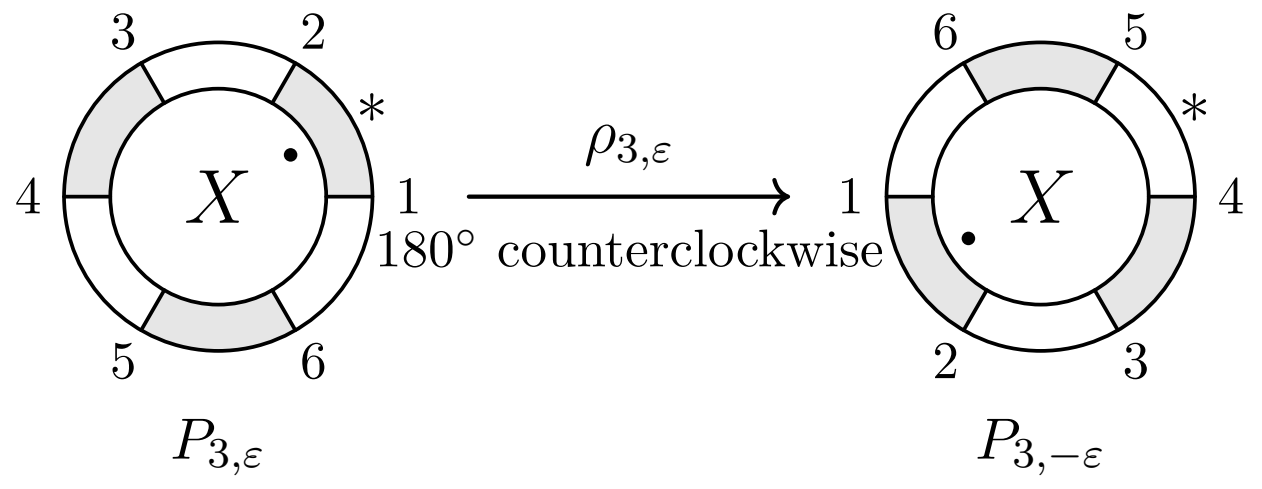

On $3$-boxes, let
$\rho_{3,\varepsilon}:P_{3,\varepsilon}\to P_{3,-\varepsilon}$
denote the counterclockwise $180^\circ$ rotation shown in
Figure~\ref{fig:annular-half-turn}.
It shifts the six boundary positions by three places
and exchanges the two shadings.
We write $\rho_3$ when the source shading is clear.

\begin{remark}
\label{rmk:matrix-convention}
Every matrix in this paper acts on column coordinate vectors.  Thus columns
index source-basis elements and rows index target-basis elements.  In
particular, the matrix of a linear map $T:V\to W$ has entries
$T_{\text{target},\text{source}}$.  This convention is used for the fusion
matrices, the Fourier matrices, the character tables, and the duality
permutation matrices.
\end{remark}

\label{subsec:presentations}
We use the generators-and-relations framework of
\cite{LuThesis}. Let
$L=\coprod_{m,\varepsilon}L_{m,\varepsilon}$ be a set of labels with
specified box numbers and shadings.
\begin{definition}\label{def:planar-presentation}
The universal planar algebra $\mathcal P(L)$ has, in each output type
$(m,\varepsilon)$, the vector space freely spanned by isotopy classes
of compatible $L$-labelled planar tangles. Tangle substitution acts by
gluing, extended multilinearly. Given local relations $R$, the planar
algebra presented by $(L,R)$ is obtained by identifying the two sides
of each relation in every planar context, and extending these
identifications linearly.
\end{definition}
In particular, $\mathcal P(L_2)$ denotes the universal planar algebra
on 2-box labels, before the proposed triangle relations are imposed.
Throughout the construction, these relations first specify computation
rules on free diagrams. They need not yet give consistent evaluations.
We reserve $P$ for an actual planar algebra satisfying the presentation.
A reference to a TR presentation of $\mathcal P(L_2)$ means the universal
object equipped with these proposed rules, not an already established
TR realization.

A planar algebra is 2-box generated if every element is a linear
combination of tangles whose internal labels belong to its 2-box
spaces. In an existing realization the labels can be chosen from a
basis of $P_{2,+}$, with rotated copies supplying the opposite shading.
In reconstruction they are initially symbols for the basis of an
abstract candidate fusion bialgebra. Proving that these symbols become
a basis of the constructed 2-box space is part of the construction.
Contracting labelled two-boxes to four-valent vertices gives the graph
notation below, retaining cyclic orders and distinguished intervals.
These vertices have no over/under-crossing data.

A normalized evaluation is first a linear functional on
$\mathcal P(L_2)_{0,\pm}$. Invariance under every defining replacement
in every closed context makes its value well defined on the presented
0-box space. Computability makes that space at most one-dimensional.
The normalization $Z(\varnothing)=1$ prevents the empty diagram from
vanishing, so it is exactly one-dimensional. This establishes consistency
at degree zero. Preservation of the candidate 2-box data and
nondegeneracy at other degrees are separate assertions proved or stated
explicitly in the construction theorem.

\subsection{Fusion bialgebras}

We use the algebraic definition of fusion bialgebra from
\cite{LuThesis}.  Throughout this paper, a fusion bialgebra
of rank $n$ is indexed by
$I=\{0,1,\ldots,n-1\}$.
Unless another label set is specified, every free label in the one-sided
part of the paper ranges over $I$.  

\begin{definition}
\label{def:fusion-bialgebra}
A fusion bialgebra over $\mathbb K$ is a finite-dimensional vector space
$B$ with associative products $\cdot$ and $*$, a $\mathbb K$-linear involutive anti-automorphism
$x\mapsto\bar x$ of the convolution algebra $(B,*)$, a linear functional $d:B\to\mathbb K$,
and a basis $\{p_0,\ldots,p_{n-1}\}$ with distinguished element $p_0$,
satisfying the following conditions:
\begin{enumerate}[label=\textup{(\roman*)},leftmargin=2.4em]
\item $p_j\cdot p_k=\delta_{jk}p_j$.
\item $p_j*p_k=\sum_sN_{jk}^{s}p_s$ for scalars $N_{jk}^{s}\in\mathbb K$.
\item An involution on $I$ satisfies $\bar p_k=p_{\bar k}$ and $\bar p_0=p_0$.
\item $d(p_i)\ne0$ for all $i$, and
\[
\begin{gathered}
 d(p_j*p_k)=d(p_j)d(p_k),\qquad d(p_{\bar k})=d(p_k),\\
 N_{jk}^{0}=\delta_{j,\bar k}d(p_j),\qquad
 d(p_0)=1,\qquad \sum_s d(p_s)\ne0.
\end{gathered}
\]
\end{enumerate}
\end{definition}
The scalars $N_{jk}^{s}$ are the \emph{fusion coefficients}.
The \emph{quantum dimension} of the minimal idempotent $p_i$ is
$d_i:=d(p_i)$ \cite{LuThesis}. We work over $\mathbb K=\C$.

\begin{remark}
The ordinary unit is
$1_B=\sum_{i=0}^{n-1}p_i$,
while $p_0$ is the unit of the convolution product. If $B$ is realized as the $2$-box algebra of a subfactor
planar algebra, the Schur product theorem implies that
$p_j\star p_k$ is positive. Its coefficients in the
minimal projection basis therefore satisfy
$N_{jk}^{s}\geq0$.
The algebraic definition above does not require this
nonnegativity.
\end{remark}

\begin{definition}
\label{def:based-fusion-isomorphism}
A fusion bialgebra is called semisimple if its convolution algebra
$(B,*)$ is semisimple. Its ordinary algebra is already split semisimple
in the distinguished idempotent basis. A based fusion-bialgebra
isomorphism is a linear bijection preserving both products, duality,
the distinguished element $p_0$, the functional $d$, and the distinguished
basis. In particular it induces a permutation of the minimal idempotents
fixing $0$.
\end{definition}

The linear duality $p_i\mapsto p_{\dual i}$ acts trivially on complex
scalars and is distinct from the conjugate-linear adjoint in a positive
realization. In the algebraic setting, ``minimal projection'' below refers
to the distinguished primitive idempotent. Self-adjointness is asserted
only when a compatible $*$-structure is present.

\begin{definition}
\label{def:fusion-realizability}
A fusion bialgebra $B$ is \emph{planar-algebra realizable} if there are
an irreducible, spherical, nondegenerate planar algebra $P$ and a based fusion-bialgebra isomorphism
$\phi:B\longrightarrow P_{2,+}$.
The isomorphism identifies ordinary multiplication, normalized convolution,
linear duality, $p_0$ and the unnormalized 2-box trace. A modulus compatible
with $\delta^2=\sum_i d_i$ is part of the realization.

\end{definition}

\begin{definition}
A fusion bialgebra is \emph{subfactorizable} if it can be realized as the
$2$-box space $P_{2,+}$ of an irreducible subfactor planar algebra, with the
ordinary and convolution products induced by vertical and normalized
horizontal composition, respectively.
\end{definition}

Suppose that $B$ is realized by an irreducible planar algebra $P$ with
commutative $P_{2,+}$.  We use the same symbols
$\{p_i\}_{i\in I}=\{p_0,p_1,\ldots,p_{n-1}\}$
for the pairwise orthogonal minimal projections of $P_{2,+}$.  In the
normalization of \cite{LuThesis},
$p_0=\frac{1}{\delta}E_{1,+}$,
where $\delta\neq0$ is the modulus, and the convolution is the horizontal
composition multiplied by $\delta$.  We write
\[
p_k*p_j=\sum_{i=0}^{n-1}N_{kj}^{i}p_i,
\qquad
d_i=d(p_i)=\operatorname{tr}_2(p_i),
\]
and use the global-dimension normalization
\begin{equation}
\label{eq:dimension-normalization}
\delta^2=\sum_{i=0}^{n-1}d_i,
\end{equation}
so that $\delta^2$ is the Jones index in a subfactor realization.  In the
purely algebraic classification, $\delta\in\C^\times$ is a distinguished
invertible scalar satisfying this equation. In the subfactor case $\delta>0$ is the positive Jones
modulus.

\label{subsec:2-box-diagrams}
The following conventions agree with those in \cite{LuThesis}.
The modulus is $\delta\ne0$. A contractible circle of either shading has
value $\delta$.

\begin{equation}\opLoop{0}=\delta.\label{pre:loop}\end{equation}
The Jones projection and ordinary unit are normalized by
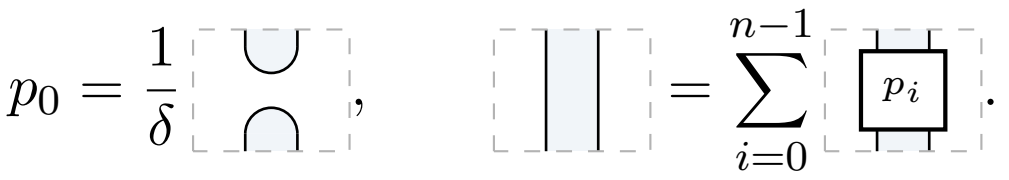
\begin{equation}
 p_0=\frac1\delta\,\opJones,
 \qquad \opIdentity=\sum_{i=0}^{n-1}\opP{i}.
 \label{pre:jones-unit}
\end{equation}
The two elementary one-face evaluations are
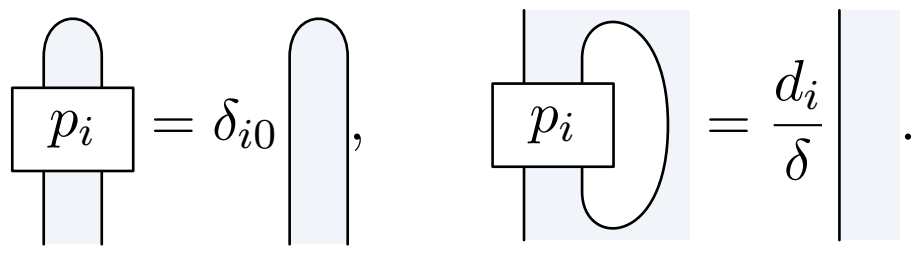
\begin{equation}
 \opCap{1}=\delta_{i0}\,\opCap{0},\qquad
 \opPartial=\frac{d_i}{\delta}\,\opStrand.
 \label{pre:one-face}
\end{equation}
Here $d_i=\operatorname{tr}_{2,+}(p_i)$ is the value obtained
by closing all four boundary points of $p_i$.
The right-hand relation in \eqref{pre:one-face} closes only
two boundary points, giving $d_i/\delta$ times the identity
strand. Closing the remaining strand contributes a factor
of $\delta$, recovering the full trace $d_i$.

Ordinary multiplication is vertical composition:
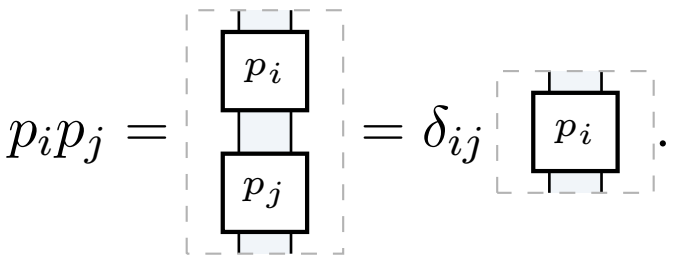
\begin{equation}
 p_i p_j=\opMult=\delta_{ij}\,\opP{i}.
 \label{pre:multiplication}
\end{equation}
Convolution is $\delta$ times horizontal composition:
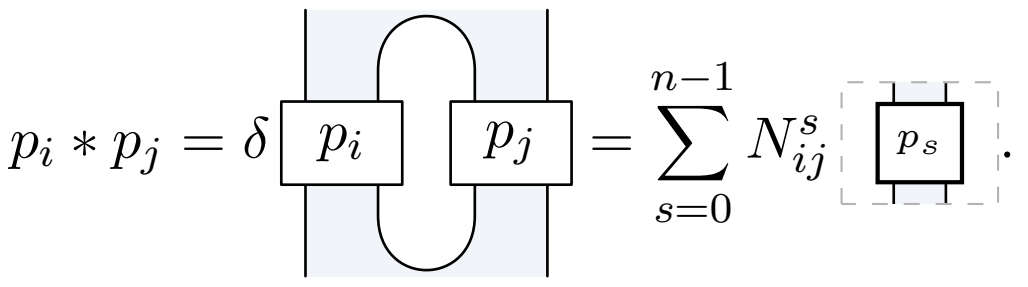
\begin{equation}
 p_i*p_j=\delta\,\opConv=\sum_{s=0}^{n-1}N_{ij}^{s}\,\opP{s}.
 \label{pre:convolution}
\end{equation}
Thus the horizontal diagram itself evaluates to
$\delta^{-1}\sum_sN_{ij}^{s}p_s$. On the opposite side we use minimal idempotents $p'_a$, dimensions
$d'_a$, and fusion coefficients $N'{}_{ab}^{c}$. The local relations
have the opposite checkerboard shading and the corresponding boundary
markings. For example, negative-side multiplication is
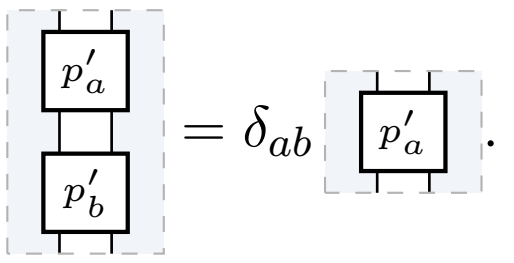
\[
 \opEditor{\drawGridBackground{2}{1}{1}
 \standardBox{1}{1}{$p'_a$}{0}
 \standardBox{1}{0}{$p'_b$}{0}}
 =\delta_{ab}\,
 \opEditor{\drawGridBackground{2}{0}{1}
 \standardBox{1}{0}{$p'_a$}{0}}.
\]
The negative-side convolution relation is
$p'_a*p'_b=\sum_c N'{}_{ab}^{c}p'_c$, with the same factor $\delta$
between convolution and horizontal composition.

The Fourier convention has the graphical form
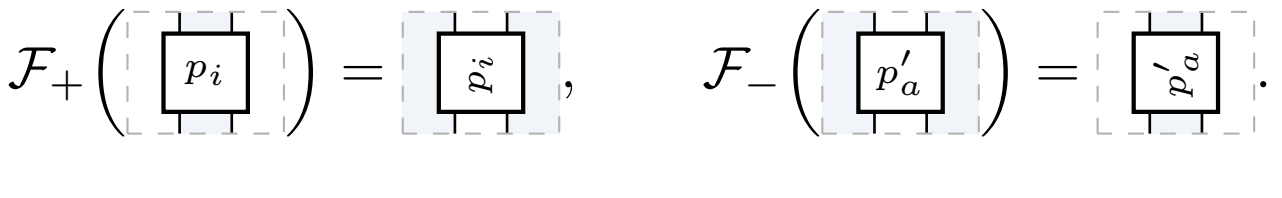
\begin{equation}\label{pre:fourier-pictures}
 \mathcal F_+\!\left(\opMarkedBox{0}{p_i}{0}\right)
 =\opMarkedBox{1}{p_i}{90},
 \qquad
 \mathcal F_-\!\left(\opMarkedBox{1}{p'_a}{0}\right)
 =\opMarkedBox{0}{p'_a}{90}.
\end{equation}
Both maps rotate the marked box counterclockwise through $90^\circ$.
The symbol inside a rotated box retains the original input label. The
right-hand picture of $\mathcal F_+(p_i)$ is not a primitive projection
$p'_a$ in general. The letter orientation follows the editor convention.
The rotation of the distinguished interval is understood. Consequently,
\begin{equation}\label{pre:fourier-duality}
 \mathcal F_-\mathcal F_+(p_i)=p_{\dual i},\qquad
 \mathcal F_+\mathcal F_-(p'_a)=p'_{\dual a}.
\end{equation}

\begin{proposition}
\label{prop:basic-2box-equations}
The fusion coefficients $N_{kj}^{i}$ and quantum dimensions $d_i$
satisfy the Frobenius reciprocity identities
\begin{align}
N_{kj}^{i}d_i
&=
N_{j\dual{i}}^{\dual{k}}d_{\dual{k}}
=
N_{\dual{i}k}^{\dual{j}}d_{\dual{j}}
=
N_{\dual{j}\dual{k}}^{\dual{i}}d_{\dual{i}}
=
N_{\dual{k}i}^{j}d_j
=
N_{i\dual{j}}^{k}d_k,
\label{eq:frobenius-reciprocity}
\end{align}
together with
\begin{align}
N_{kj}^{0}
&=\delta_{k,\dual{j}}d_k,
\label{eq:vacuum-equation}\\
N_{kj}^{i}
&=N_{\dual{j}\,\dual{k}}^{\dual{i}},
\label{eq:dual-fusion-equation}\\
d_i&=d_{\dual{i}},
\label{eq:dim-dual}\\
\delta^2&=\sum_{i\in I}d_i.
\label{eq:global-dimension}
\end{align}
These identities hold for all indicated labels in $I$.
\end{proposition}

\begin{proof}
The graphical identities follow from
\cite{LuThesis}.
The Frobenius reciprocity and duality identities also hold
in an abstract fusion bialgebra.
The modulus identity follows by taking the trace of
$1=\sum_{i\in I}p_i$.
\end{proof}

Following \cite{LuThesis}, we call
\eqref{eq:frobenius-reciprocity}--\eqref{eq:global-dimension}
the \emph{basic consistency equations}.

\begin{proposition}
\label{prop:dimension-associativity}
Let $B$ be a fusion bialgebra with fusion coefficients $N_{kj}^{i}$
and dimensions $d_i$.
The multiplicativity of the dimension functional and the
associativity of convolution give
\begin{align}
d_kd_j&=\sum_{i\in I}N_{kj}^{i}d_i
&& (j,k\in I),
\label{eq:dimension-equation}\\
\sum_{s\in I}N_{kj}^{s}N_{sr}^{\ell}
&=\sum_{s\in I}N_{jr}^{s}N_{ks}^{\ell}
&& (j,k,r,\ell\in I).
\label{eq:associativity-equation}
\end{align}
\end{proposition}

These are the coordinate forms of the defining properties
in Definition~\ref{def:fusion-bialgebra}.

Following \cite{LuThesis}, we call
\eqref{eq:vacuum-equation},
\eqref{eq:dual-fusion-equation},
\eqref{eq:dim-dual},
\eqref{eq:dimension-equation}, and
\eqref{eq:associativity-equation}
the \emph{fusion structure equations}.

\begin{theorem}
\label{prop:2box-fusion-bialgebra}
Let $P$ be a locally finite-dimensional, nondegenerate,
spherical, irreducible planar algebra over $\C$.
If $P_{2,+}$ is commutative, then $P_{2,+}$ is a fusion
bialgebra with its ordinary multiplication, normalized
convolution, half-turn duality, and unnormalized trace.
\end{theorem}

\begin{definition}
Let $P,Q$ be planar algebras. Their \emph{tensor product} $P\otimes Q$
has box spaces $(P\otimes Q)_{n,\pm}=P_{n,\pm}\otimes Q_{n,\pm}$.
For a planar tangle $T$ with $m$ input discs, its action is determined by
\[
 Z_T^{P\otimes Q}(x_1\otimes y_1,\ldots,x_m\otimes y_m)
 =Z_T^P(x_1,\ldots,x_m)\otimes Z_T^Q(y_1,\ldots,y_m).
\]
\end{definition}

\begin{definition}
Let $P,Q$ be planar algebras. Their ordered \emph{free product} $P*Q$
is the planar algebra of diagrams with two types of strings, labelled
$P$ and $Q$. Strings are disjoint, and each internal disc carries an
element of the factor matching all its incident strings. Diagrams are
taken modulo planar isotopy, linearity in the labels, and the evaluation
relations of the two factors. An $n$-box has $4n$ boundary endpoints,
ordered $PQQP\cdots PQQP$ in the plus shading, with compatible factor
shadings. The minus shading uses the reversed factor order. A planar
tangle acts by replacing each string with a parallel pair of the two
types and gluing matching endpoints. Closed strings are evaluated using
the modulus of their factor. This is the construction of
\cite[Definition~2.6]{LiuER}.
\end{definition}

The corresponding constructions for fusion bialgebras are as follows.
\begin{definition}
Let $A,C$ be fusion bialgebras with bases $\{x_i\}$ and $\{y_j\}$.
Their tensor product has basis $\{x_i\otimes y_j\}$, distinguished element
$x_0\otimes y_0$, componentwise products and involution, and
$d(x_i\otimes y_j)=d_A(x_i)d_C(y_j)$.
Their ordered free product is the subspace of $A\otimes C$ with basis
\[
 \{x_i\otimes y_0:i\in I_A\}\ \cup\
 \{1_A\otimes y_j:j\in I_C\setminus\{0\}\},
 \qquad 1_A=\sum_i x_i,
\]
with the inherited products, involution and dimension functional, and
with distinguished element $x_0\otimes y_0$.
\end{definition}
These are the constructions of \cite{LuThesis}.
\subsection{Triangle and Yang--Baxter relations}\label{sec:TR-definitions}

The triangle relation is motivated by a local reduction: replace a
triangular face by a linear combination of diagrams with two vertices.
Suppressing shading, labels and coefficients, its schematic form is
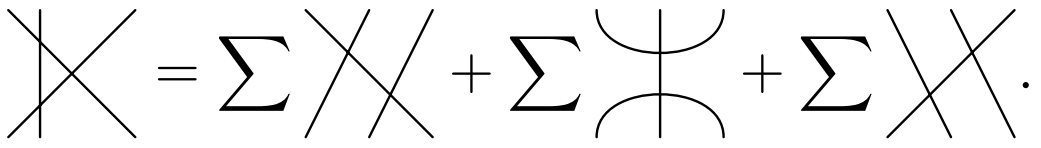
\begin{equation}
 \opSchematicTriangle
 =\sum\opOnThree{12}+\sum\opOnThree{13}+\sum\opOnThree{14}.
 \label{pre:TR-schematic}
\end{equation}
The sums denote weighted sums over vertex labels. The weights are suppressed
only in this schematic equation. Each term preserves the six boundary ports.

\begin{figure}[H]
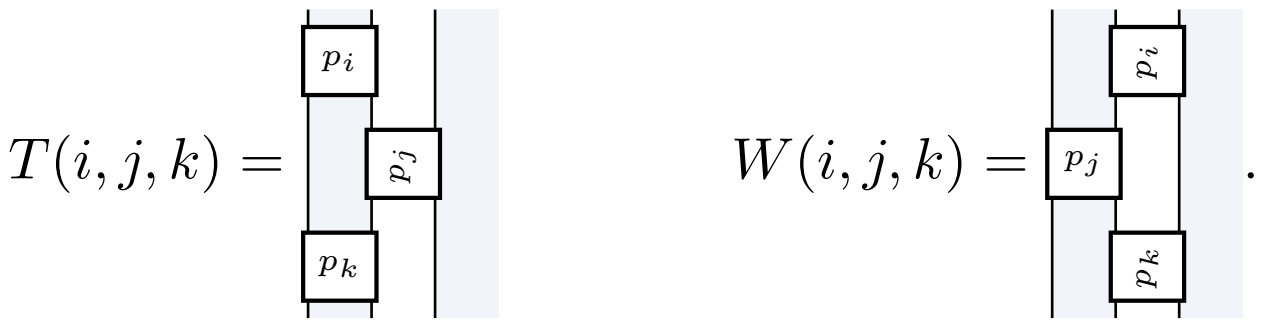

\centering
\[
  T(i,j,k)=\opT{0}
  \qquad\qquad
  W(i,j,k)=\opT{1}.
\]
\caption{The diagrams $T$ and $W$.}
\label{fig:triangle-convention}
\end{figure}

Let $P$ be a locally finite-dimensional, spherical, nondegenerate planar algebra
generated by its commutative $2$-box space $P_{2,+}$, and let
$\{p_i\}_{i\in I}$ be the distinguished minimal projections.  For
$i,j,k\in I$, let
$T(i,j,k),\; W(i,j,k)\in P_{3,+}$
denote the two oriented triangular $3$-face diagrams in the convention of
\Cref{fig:triangle-convention}.  The diagram $W(i,j,k)$ is the opposite
orientation, expressed using the same minimal idempotents. We also use the three internal-face-free
families
$E_1(x,y),\; E_2(x,z),\; E_3(y,z)\in P_{3,+}$.
Define
\begin{equation}
\label{eq:Espan-definition}
\Espan
:=
\operatorname{span}_{\C}
\bigl\{E_1(x,y),E_2(x,z),E_3(y,z):x,y,z\in I\bigr\}
\subseteq P_{3,+}.
\end{equation}

\begin{figure}[H]
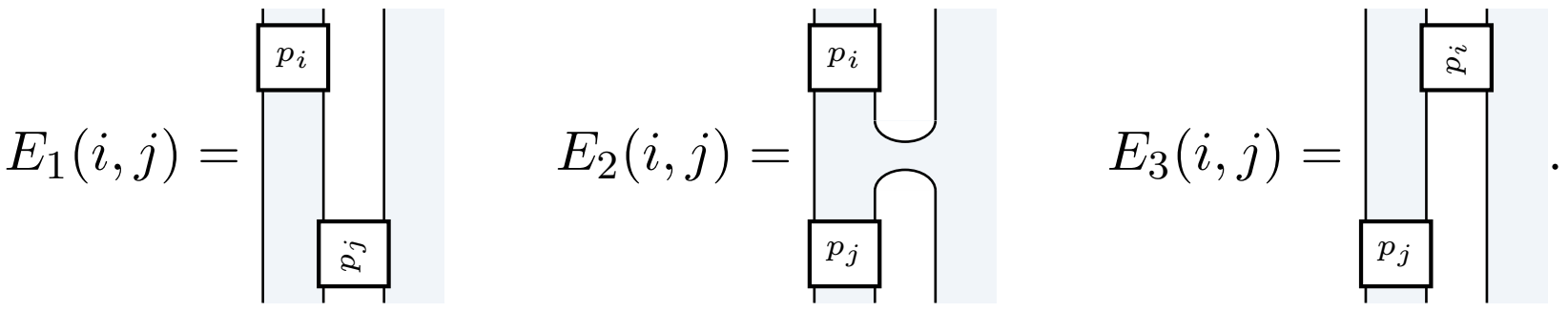

\centering
\[
  E_1(i,j)=\TRfamily{1}{i}{j}
  \qquad
  E_2(i,j)=\TRfamily{2}{i}{j}
  \qquad
  E_3(i,j)=\TRfamily{3}{i}{j}.
\]
\caption{The diagrams $E_1,E_2,E_3$.}
\label{fig:E-diagrams}
\end{figure}

\begin{definition}
\label{def:er-planar-algebra}
In the conventions of \Cref{fig:E-diagrams}, the planar algebra satisfies
an exchange relation if, for every $i,j$, there are scalars
$a_{st}^{ij}$ and $b_{\ell m}^{ij}$ such that
\begin{equation}\label{pre:exchange-relation}
 E_1(i,j)=\sum_{s,t\in I}a_{st}^{ij}E_2(s,t)
       +\sum_{\ell,m\in I}b_{\ell m}^{ij}E_3(\ell,m).
\end{equation}
Together with 2-box generation and the elementary closure relations,
this is the exchange-relation (ER) presentation of
\cite{LuThesis}. Each required external shading uses its
own coefficients and minimal idempotents.
\end{definition}

In the TR relation the input is a triangle with three labelled vertices.
Both local orientations are expanded in the face-free families. The
following formulas use only the minimal idempotents of $P_{2,+}$. Rotated
local occurrences carry the induced boundary markings and labels.

\begin{definition}
\label{def:tr-planar-algebra}
The planar algebra $P$ is called a \emph{triangle relation planar algebra},
or a \emph{TR planar algebra}, if every nonzero triangle of both orientations
belongs to $\Espan$.  Equivalently, there are coefficient tables $c$ and
$\widetilde c$ such that
\begin{align}
T(i,j,k)
&=
\sum_{x,y}c_{1,xy}^{ijk}E_1(x,y)
+
\sum_{x,z}c_{2,xz}^{ijk}E_2(x,z)
+
\sum_{y,z}c_{3,yz}^{ijk}E_3(y,z),
\label{eq:plus-tr-expansion}\\
W(i,j,k)
&=
\sum_{x,y}\widetilde c_{1,xy}^{ijk}E_1(x,y)
+
\sum_{x,z}\widetilde c_{2,xz}^{ijk}E_2(x,z)
+
\sum_{y,z}\widetilde c_{3,yz}^{ijk}E_3(y,z)
\label{eq:W-tr-expansion}
\end{align}
for every $i,j,k\in I$ for which the corresponding triangle is nonzero.
The summation indices $x,y,z$ range over $I$.
\end{definition}

In the diagram convention fixed above, these equations read
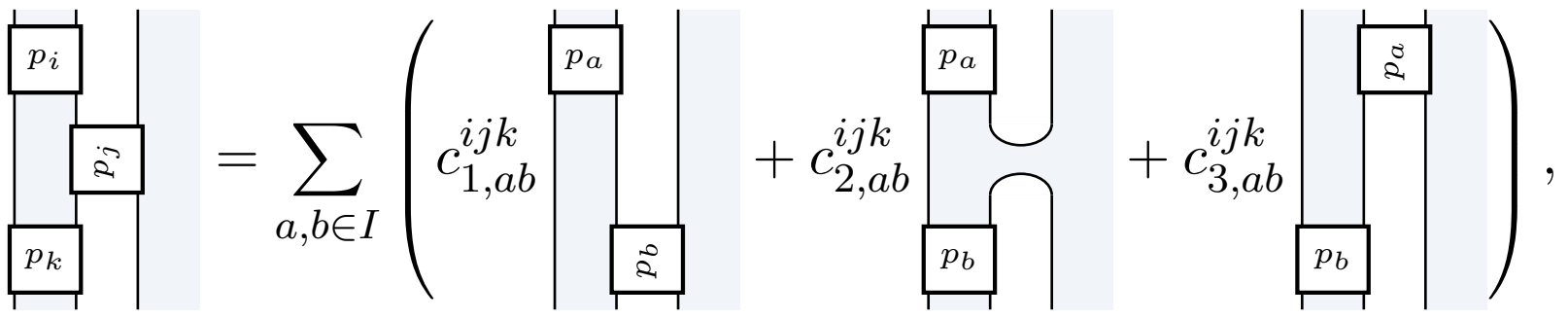
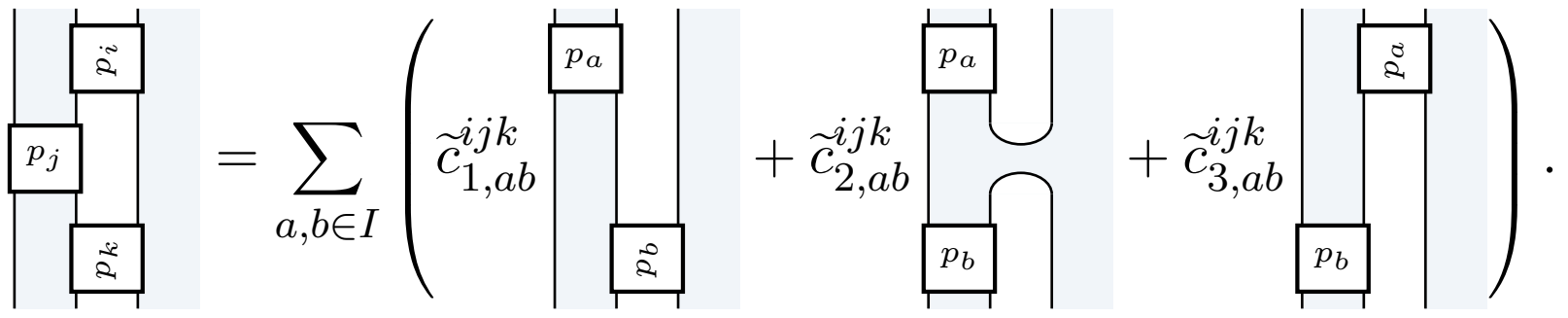
\begin{equation}\label{pre:TR-native-picture}
 \opT{0}=\sum_{a,b\in I}\left(
 c_{1,ab}^{ijk}\,\opE{1}+c_{2,ab}^{ijk}\,\opE{2}
 +c_{3,ab}^{ijk}\,\opE{3}\right),
\end{equation}
\begin{equation}\label{pre:TR-W-picture}
 \opT{1}=\sum_{a,b\in I}\left(
 \widetilde c_{1,ab}^{ijk}\,\opE{1}
 +\widetilde c_{2,ab}^{ijk}\,\opE{2}
 +\widetilde c_{3,ab}^{ijk}\,\opE{3}\right).
\end{equation}
The table $\widetilde c$ expands $W$ in the same families
$E_\alpha$ used for $T$. All labels in these two expansions belong to
$P_{2,+}$. 

\begin{definition}
\label{def:ybr-planar-algebra}
The planar algebra $P$ is called a \emph{Yang--Baxter relation planar
algebra}, or a \emph{YBR planar algebra}, if every $W$-triangle is a linear
combination of $T$-triangles. That is, for every $i,j,k\in I$ there are
coefficients $\lambda_{abc}^{ijk}$ such that
\begin{equation}
\label{eq:ybr-relation}
W(i,j,k)=\sum_{a,b,c\in I}\lambda_{abc}^{ijk}T(a,b,c).
\end{equation}
\end{definition}

\begin{proposition}
\label{prop:triangle-basis}
Let $P$ be a YBR planar algebra.  Then
\[
\operatorname{COB}_{+}
=
\left\{T(i,j,k):N_{kj}^{i}\neq0\right\}
\]
is an orthogonal basis of $P_{3,+}$.
\end{proposition}

\begin{proof}
The YBR relation implies that the triangle diagrams span $P_{3,+}$.
Orthogonality of the minimal idempotents and a direct trace
calculation show that the triangles with $N_{kj}^{i}\ne0$
are mutually orthogonal with nonzero self-pairings, while
the remaining triangles vanish. Hence the displayed family
is a basis. For the detailed calculation
and the YBR spanning argument, see \cite{LuThesis}.
\end{proof}

\begin{proposition}
\label{prop:ybr-one-sided-tr}
Let $P$ be a YBR planar algebra.  Then $P$ is a TR planar algebra if and only
if every nonzero $T(i,j,k)$ belongs to $\Espan$.
\end{proposition}

\begin{proof}
The forward implication follows from the definition.  Conversely, assume
that every $T(a,b,c)$ lies in $\Espan$.  By the YBR relation
\eqref{eq:ybr-relation}, every $W(i,j,k)$ is a linear combination of these
$T$-triangles and hence also lies in $\Espan$.  Thus both oriented
$3$-faces satisfy triangle relations.
\end{proof}

\begin{remark}
Since each $E_\alpha$ is obtained from the ordinary identity resolution and
is therefore a linear combination of $T$-triangles, the general TR
relations also imply the YBR spanning relation in the present ambient
class.
\end{remark}

We call a fusion bialgebra ER-, YBR-, or TR-realizable when its
planar-algebra realization can be chosen with the respective relations.
We also call these ER, YBR, and TR fusion bialgebras,
respectively.

The graphical consistency argument is formulated in $\mathcal P(L_2)$,
where $L_2$ is a finite set of 2-box labels with its prescribed rotated
and shaded copies. A rule replaces a circle, a monogon, a bigon, or a triangle by a
finite linear combination of diagrams with fewer faces. In particular,
a triangle is replaced by internal-face-free diagrams with fewer vertices.
Each rule is specified for every permitted boundary marking and shading.

\section{Consistency of local skein reductions}
\label{sec:evaluation}
We work in the skein setting of Section~\ref{subsec:presentations},
over a field of characteristic zero.
The arguments in this section use only the planar structure of
the diagrams and the local reduction rules.

\subsection{Evaluation and consistency of skein relations}
A closed diagram is regarded as a graph in the sphere, so its face count
includes the unbounded face of a plane drawing. A face whose boundary consists of $k$ edges is called a $k$-face.
A circle with no vertices is treated as a $0$-face.

\begin{lemma}\label{cons:eight}\label{lem:eight-triangles}
Let $D$ be a nonempty connected $0$-box diagram containing no $0$-, $1$-,
or $2$-face. Then it contains at least eight $3$-faces. If $f_r$ denotes
the number of degree-$r$ faces, then
\begin{equation}\label{cons:eightidentity}
f_3=8+\sum_{r\ge5}(r-4)f_r.
\end{equation}
\end{lemma}
\begin{proof}
The graph is four-valent. Writing $V,E,F$ for its vertex, edge and face
counts gives $E=2V$, $V-E+F=2$, and $\sum_r rf_r=2E$. Hence
$\sum_r(4-r)f_r=4F-2E=8$, which is the asserted identity.
\end{proof}

\begin{lemma}\label{cons:lowbox}\label{lem:low-degree-internal-face}
If a connected $r$-box diagram, $r=1,2,3$, has an internal face, then it
has an internal $1$-, $2$-, or $3$-face.
\end{lemma}
\begin{proof}
See \cite{LuThesis} for the planar graph argument.
It is independent of the choice of local skein relations. An internal
face is a complementary region disjoint from the outer boundary.
\end{proof}

\begin{remark}
The restriction to $r\le3$ is essential. A quadrilateral with two external
legs at each vertex is a 4-box with an internal $4$-face and no
internal $0$-, $1$-, $2$-, or $3$-face:
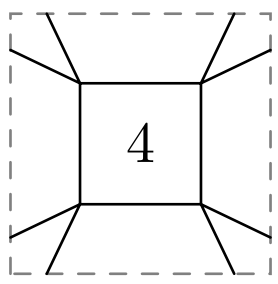
\[
\begin{tikzpicture}[scale=.40,line width=.5pt,line cap=round]
\draw[dashed,gray] (-2.15,-2.15) rectangle (2.15,2.15);
\draw (-1,1)--(1,1)--(1,-1)--(-1,-1)--cycle;
\draw (-1,1)--(-1.55,2.15) (-1,1)--(-2.15,1.55);
\draw (1,1)--(1.55,2.15) (1,1)--(2.15,1.55);
\draw (1,-1)--(1.55,-2.15) (1,-1)--(2.15,-1.55);
\draw (-1,-1)--(-1.55,-2.15) (-1,-1)--(-2.15,-1.55);
\node at (0,0) {$4$};
\end{tikzpicture}
\]
Shading and vertex labels are suppressed. In particular, an internal
face in a 4-box need not admit a canonical reduction.
\end{remark}

\begin{definition}\label{def:computation-operations}
Consider labelled diagrams with prescribed local skein relations.
We use the following moves to apply these relations locally,
leaving the rest of the diagram unchanged and extending
linearly to diagram expansions.
The pictures suppress shading. Each represents the corresponding rule
in either shading. 
\begin{enumerate}[label=\textbf{Move \arabic*:},start=0,leftmargin=2.7cm,itemsep=5pt]
\item Remove a vertex-free circle bounding an empty disc:
\(\opBare{\draw (0,0) circle[radius=0.4];}=\delta\).
\item Reduce a monogon:
\(\opMono=\lambda\,
\opBare{\draw (-0.55,-0.5).. controls (-0.55,0.6) and (0.55,0.6)..(0.55,-0.5);}\),
where $\lambda$ is the coefficient of the applicable one-face relation.
\item Reduce a bigon:
\(\opBigon=\sum_a\lambda_a\,\opCross\),
using the ordinary product or convolution relation. The coefficients $\lambda_a$ include all normalization factors
required by the corresponding local tangle.
\item Expand a triangle by the triangle relation:
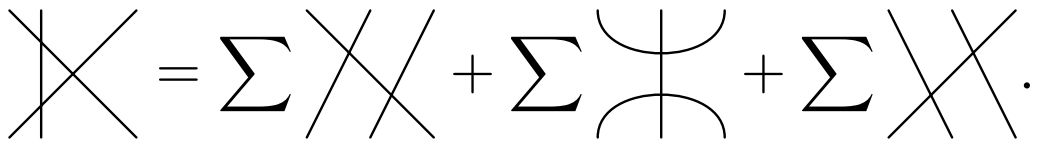
\[
\opSchematicTriangle=
\sum\opOnThree{12}+\sum\opOnThree{13}+\sum\opOnThree{14}.
\]
Each sum is weighted by the appropriate triangle coefficients.
\item Apply the resolution of identity to replace
\(
\opBare{
  \draw (-0.3,-0.55)--(-0.3,0.55);
  \draw (0.3,-0.55)--(0.3,0.55);
}
\)
by
\(
\sum\opCross
\).
\end{enumerate}
An \emph{evaluation path} is a finite succession of
the above moves, applied from left to right and
termwise to linear combinations of diagrams.
An application of Move~0, 1, 2, or~3 is called
a \emph{canonical reduction}.
A \emph{canonical evaluation path} uses only
canonical reductions.
An evaluation path is complete when every resulting
term is a coefficient times the empty diagram.
\end{definition}

\begin{definition}
A linear function
$Z:\mathcal P(L_2)_{0,\pm}\longrightarrow\mathbb K$
with $Z(\varnothing)=1$ is a \emph{pre-partition function} if it is invariant under Moves
$0$--$3$, and a \emph{partition function} if it is also invariant under
Move $4$. For a disconnected closed diagram, its value is defined to be
the product of the values of its connected components.
\end{definition}

\begin{definition}
The reduction rules are \emph{computable} if every closed diagram admits
a complete canonical evaluation. They are \emph{canonically consistent}
if all complete canonical evaluations of each closed diagram have the
same scalar value, and \emph{consistent} if this holds for all complete
computations using Moves $0$--$4$.
\end{definition}

Moves $0$--$3$ will be used to define a candidate
evaluation by averaging complete canonical evaluations.
Canonical consistency ensures that every such computation
gives the same value and that these moves preserve
the resulting evaluation.
Move~4 applies the prescribed identity-resolution relation. It is not used in defining the candidate
evaluation, and we must also verify that it preserves
the evaluation in every closed context.

Moves $1$--$3$ strictly decrease the face number
in each resulting term and preserve connectedness
of each nonzero term. Move $0$ removes
vertex-free circles. We therefore use the face number
as the induction parameter for canonical evaluations.
Identity resolution, however, can increase the face
number and join connected components.

\begin{theorem}\label{thm:computability}
The reduction rules are computable.
\end{theorem}
\begin{proof}
A circle is handled by Move $0$. If a nonempty connected closed
diagram has no monogon or bigon, Lemma~\ref{cons:eight} supplies a triangle.
There is therefore a canonical reduction whenever one is needed, and the
face number decreases. The process terminates. A diagram with several
components is evaluated componentwise.
\end{proof}

Let $\varnothing$ denote the empty closed diagram and set $Z(\varnothing)=1$.
For a nonempty connected diagram $D$, let $\Gamma(D)$ denote the finite nonempty set of complete
canonical evaluation schemes for $D$. Execution order between distinct terms is not counted again. Schemes that differ only in the order in which distinct terms
are evaluated are identified.
Finite branching and decreasing face number make
$\Gamma(D)$ finite. Let $\operatorname{ev}_\gamma(D)$ be the resulting scalar.
Define
\begin{equation}\label{eq:average-evaluation}
Z(D)=\frac{1}{|\Gamma(D)|}\sum_{\gamma\in\Gamma(D)}
\operatorname{ev}_\gamma(D),\qquad Z(\varnothing)=1,
\end{equation}
and extend linearly and multiplicatively over connected components.
The average is defined before canonical consistency is known.
Canonical consistency asserts that every complete canonical evaluation gives this same value. The canonical reductions below
are discussed on connected diagrams. Component-changing replacements
belong to Section~\ref{subsec:full-partition}.

If a first move $\alpha$ gives
$D_\alpha=\sum_i a_iD_{\alpha,i}$, the notation $Z(D_\alpha)$ means
$\sum_i a_iZ(D_{\alpha,i})$. The index $\alpha$ specifies the actual face
and the chosen reduction, not merely the degree of that face.

\begin{lemma}\label{lem:path-invariance}
For the function in \eqref{eq:average-evaluation}, the following conditions
are equivalent:
\begin{enumerate}[label=\textup{(\roman*)},leftmargin=2.4em]
\item All complete canonical evaluations of every closed diagram have the same value.
\item $Z$ is invariant under every canonical move.
\item $Z$ is a pre-partition function.
\end{enumerate} 
\end{lemma}
\begin{proof}
Suppose all complete computations have the same value. After any first
move $D\mapsto\sum_i a_iD_{\alpha,i}$, choose a complete computation
of each term. Its value is $Z(D_{\alpha,i})$, so the resulting complete
computation of $D$ has value $\sum_i a_iZ(D_{\alpha,i})=Z(D_\alpha)$.
It also equals the average $Z(D)$. Thus each move preserves $Z$.
Conversely, a complete canonical evaluation ends in a scalar $a$.
Each step preserves $Z$, and linearity together with $Z(\varnothing)=1$
gives $Z(D)=a$. Thus every complete computation has the same value.
Conditions \textup{(ii)} and \textup{(iii)} are equivalent
by the definition of a pre-partition function.
\end{proof}

\begin{lemma}\label{lem:first-step}
All complete canonical evaluations of every closed diagram have the
same value if and only if, for every connected closed diagram $D$ and
any two canonical first moves $\alpha,\beta$,
\begin{equation}\label{eq:first-step-comparison}
Z(D_\alpha)=Z(D_\beta).
\end{equation}
\end{lemma}
\begin{proof}
The forward implication follows from Lemma~\ref{lem:path-invariance}:
$Z(D_\alpha)=Z(D)=Z(D_\beta)$.

For the converse, induct on the face number. The empty diagram and
circle computations supply the initial cases. Assume path independence
for all connected closed diagrams with fewer than $N$ faces, and let
$D$ have $N$ faces. Split a complete evaluation into its first move
and all subsequent moves. If the first move gives
$D_\alpha=\sum_i a_iD_{\alpha,i}$, then each $D_{\alpha,i}$ is smaller.
The induction hypothesis says that, once the first move is fixed,
\emph{every choice from the second step onwards} has the same total value
$\sum_i a_iZ(D_{\alpha,i})=Z(D_\alpha)$.

Now change the first move to $\beta$. The assumed equality
\eqref{eq:first-step-comparison} identifies the two resulting values.
Thus the hypothesis controls the choice of the first step, while induction
controls all subsequent choices. Together they give independence of the
entire evaluation path. Its common value is the average $Z(D)$.
\end{proof}

Repeated internal-face reductions express every low-box diagram as a
linear combination of reduced diagrams once the local relations are
imposed. For box numbers at most three the terminal diagrams have no
internal face. For box number four only internal $0$-, $1$-, $2$-, and
$3$-faces are removed. In any planar algebra satisfying the triangle relations, these
reduced diagrams span the corresponding box spaces.
For $0\le r\le3$, let $O_r$ be the set of internal-face-free $r$-box
diagrams. In particular $O_0=\{\varnothing\}$. The diagrams are considered up to planar isotopy fixing
the boundary pointwise. Vertex labels are suppressed in the
pictures but restored in each test.
\begin{equation}\label{eq:On-zero-one}
O_0=\{\varnothing\},\qquad O_1=\{\OnOne\}.
\end{equation}
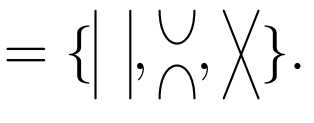
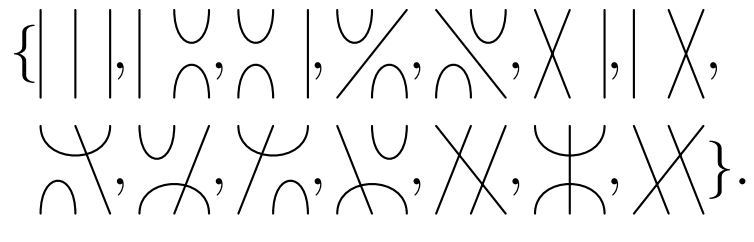
\begin{equation}\label{eq:On-two}
O_2=\{\OnTwo{1},\OnTwo{2},\OnTwo{3}\}.
\end{equation}
\begin{equation}\label{eq:On-three}
\begin{aligned}
O_3=\{&\OnThree{1},\OnThree{2},\OnThree{3},\OnThree{4},
\OnThree{5},\OnThree{6},\OnThree{7},\\
&\OnThree{8},\OnThree{9},\OnThree{10},\OnThree{11},
\OnThree{12},\OnThree{13},\OnThree{14}\}.
\end{aligned}
\end{equation}
For $O_4$ we use a different stopping condition: it consists of 4-box
diagrams with no internal $1$-, $2$-, or $3$-face, with vertex-free circles
already removed. It has $84$ underlying types, including the two
quadrilateral types. A detailed table is in Appendix~\ref{app:O4}. 
The pairing $\langle H,B\rangle$ below means planar closure, without positivity.

\begin{theorem}\label{cons:graphthm}\label{thm:finite-consistency}\label{thm:quasi-partition}
Let $Z$ be given by \eqref{eq:average-evaluation}. Suppose that for each local configuration of type
\begin{equation}\label{cons:list}
\begin{gathered}
(1,1),\ (2,1),\ (2,2)\text{-e},\ (2,2)\text{-v},\\
(3,1),\ (3,2)\text{-e},\ (3,2)\text{-v},\ (3,3)\text{-e},
\end{gathered}
\end{equation}
$Z$ is invariant under Moves $0$--$3$ on its closures with every
compatible element of $O_r$, $0\le r\le3$.
Then $Z$ is a pre-partition function. Conversely, every pre-partition
function satisfies these tests.
\end{theorem}

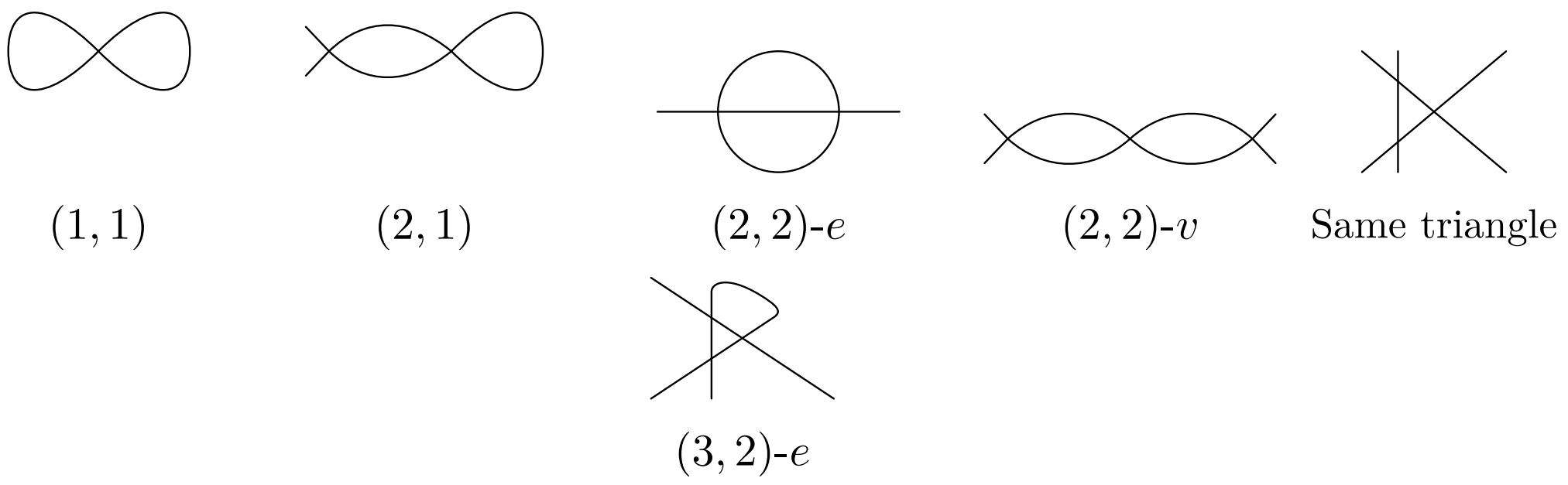
\begin{figure}[htbp]
\centering
\begin{tikzpicture}[
  x=0.011cm,
  y=-0.011cm,
  line width=0.45pt,
  line cap=round,
  line join=round,
  every node/.style={
    inner sep=0pt,
    font=\normalsize
  }
]
% (1,1)
\begin{scope}[shift={(150,88)}]
  \draw (0,0)
    .. controls (-43,-43) and (-70,-37) .. (-70,0)
    .. controls (-70,37) and (-43,43) .. (0,0)
    .. controls (43,-43) and (71,-37) .. (71,0)
    .. controls (71,37) and (43,43) .. (0,0);
\end{scope}

% (2,1)
\begin{scope}[shift={(424,88)}]
  \draw (-113,-19)--(-95,0)
    .. controls (-67,27) and (-33,27) .. (0,0)
    .. controls (43,-43) and (71,-37) .. (71,0)
    .. controls (71,37) and (43,43) .. (0,0)
    .. controls (-33,-27) and (-67,-27) .. (-95,0)
    --(-113,19);
\end{scope}

% (2,2)-e
\begin{scope}[shift={(678,135)}]
  \draw (0,0) circle[radius=47];
  \draw (-94,0)--(94,0);
\end{scope}

% (2,2)-v
\begin{scope}[shift={(951,156)}]
  \draw (-113,-19)--(-95,0)
    .. controls (-66,26) and (-29,26) .. (0,0)
    .. controls (29,-26) and (66,-26) .. (95,0)
    --(113,19);

  \draw (-113,19)--(-95,0)
    .. controls (-66,-26) and (-29,-26) .. (0,0)
    .. controls (29,26) and (66,26) .. (95,0)
    --(113,-19);
\end{scope}

% Different expansion choices on the same triangle
\begin{scope}[shift={(1187,135)}]
  \draw (-56,-47)--(56,47);
  \draw (-56,47)--(56,-47);
  \draw (-28,-47)--(-28,47);
\end{scope}

% (3,2)-e
\begin{scope}[shift={(650,311)}]
  \draw (-71,-47)--(71,47);
  \draw (-24,47)--(-24,-36)
    .. controls (-24,-43) and (-12,-49) .. (13,-34)
    .. controls (26,-26) and (32,-21) .. (24,-16)
    --(-71,47);
\end{scope}

\node at (150,225)  {$(1,1)$};
\node at (403,225)  {$(2,1)$};
\node at (678,225)  {$(2,2)\text{-}e$};
\node at (951,225)  {$(2,2)\text{-}v$};
\node at (1187,225) {\small Same triangle};
\node at (650,401)  {$(3,2)\text{-}e$};
\end{tikzpicture}
\caption{Local configurations for the consistency analysis.}
\label{fig:consistency-local-configurations}
\end{figure}

\begin{figure}[htbp]
\centering
\begin{tikzpicture}[
  x=0.013cm,
  y=-0.013cm,
  line width=0.45pt,
  line cap=round,
  line join=round,
  every node/.style={
    inner sep=0pt,
    font=\normalsize
  }
]
% (3,1)
\begin{scope}[shift={(132,136)}]
  \draw (-28,47)--(-28,-49)
    .. controls (-28,-74) and (-54,-83) .. (-64,-68)
    .. controls (-77,-53) and (-37,-34) .. (0,0)
    --(57,47);
  \draw (-56,47)--(57,-47);
\end{scope}

% (3,2)-v
\begin{scope}[shift={(325,136)}]
  \draw (-70,-75)
    .. controls (-59,-73) and (-43,-70) .. (-33,-65)
    .. controls (-23,-60) and (-23,-50) .. (-23,-39)
    --(-23,47);

  \draw (-56,-84)
    .. controls (-53,-61) and (-55,-42) .. (-47,-32)
    .. controls (-39,-24) and (-18,-12) .. (0,0)
    --(71,47);

  \draw (-70,47)--(71,-47);
\end{scope}

% (3,3)-v
\begin{scope}[shift={(533,136)}]
  \draw (-71,-47)--(70,47);
  \draw (-71,47)--(70,-47);
  \draw (-24,-47)--(-24,47);
  \draw (23,-47)--(23,47);
\end{scope}

% (3,3)-e
\begin{scope}[shift={(764,136)}]
  \draw (-95,-47)
    .. controls (-95,-5) and (-57,24) .. (0,24)
    .. controls (57,24) and (94,-5) .. (94,-47);

  \draw (-95,47)
    .. controls (-95,5) and (-57,-24) .. (0,-24)
    .. controls (57,-24) and (94,5) .. (94,47);

  \draw (0,-47)--(0,47);
\end{scope}

\node at (132,225) {$(3,1)$};
\node at (325,225) {$(3,2)\text{-}v$};
\node at (533,225) {$(3,3)\text{-}v$};
\node at (764,225) {$(3,3)\text{-}e$};
\end{tikzpicture}
\caption{Local configurations involving 3-faces.}
\label{fig:consistency-triangle-configurations}
\end{figure}
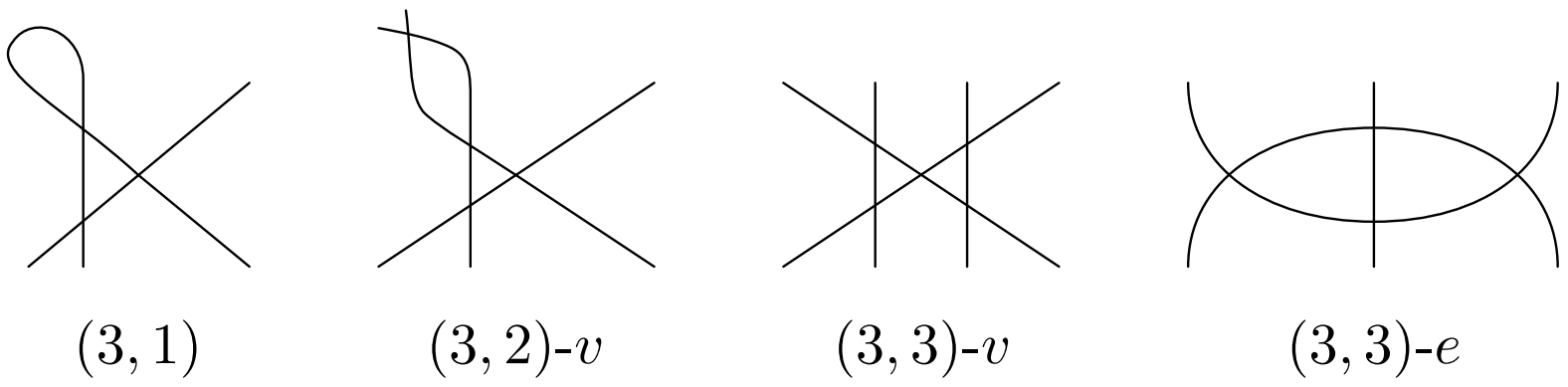

\begin{proof}[Proof of Theorem~\ref{cons:graphthm}]
Induct on the face number, starting with the empty diagram and circle
reductions. Assume that every connected closed diagram with fewer
than $N$ faces has path-independent canonical evaluation, and let $D$
have $N$ faces. Each term after a first reduction has fewer than $N$
faces. By Lemma~\ref{lem:path-invariance}, further canonical reductions
on these terms preserve $Z$. It remains to compare the choices of the
first reduction, as in the induction proving Lemma~\ref{lem:first-step}.

If the two move discs are disjoint, locality gives equality of the
resulting linear combinations $D_{m,k}=D_{k,m}$. Hence
\begin{equation}\label{eq:disjoint-first-steps}
Z(D_m)=Z(D_{m,k})=Z(D_{k,m})=Z(D_k).
\end{equation}
The end equalities are applications of the induction hypothesis through
Lemma~\ref{lem:path-invariance}. The middle equality is graphical.

For overlapping faces, cut out their local configuration $H$ and write
$D=\langle H,B\rangle$. For every type in \eqref{cons:list}, $B$ has box number at most
three. If it has no internal face, the corresponding finite test in the hypotheses of the theorem gives $Z(D_m)=Z(D)=Z(D_k)$.
Otherwise, remove a circle if one occurs, or use Lemma~\ref{cons:lowbox}
to find a reducible internal face $F_c$ in the complement. Its move
disc is disjoint from both prescribed moves. Apply
\eqref{eq:disjoint-first-steps} twice to obtain
\[
Z(D_m)=Z(D_c)=Z(D_k).
\]

For $(3,3)$-v, write $D=\langle H,R\rangle$, where $H$ is the shared-vertex configuration with 3-faces $A,B$. An internal reducible face of $R$ gives a disjoint intermediate reduction, as above. After removing vertex-free circles, it remains to consider $R\in O_4$. The comparisons in \eqref{cons:list} have already been established at the current face number. Disjoint reductions commute, and their subsequent evaluations agree by the smaller-face induction hypothesis. We will find a face $C$ for which these comparisons give
\[
 Z(D_A)=Z(D_C)=Z(D_B).
\]
The boundary labels of $H$ are shown in Figure~\ref{fig:33v-redundancy}. All intersections are four-valent vertices. Labels, distinguished intervals and shadings are carried along with the rotated or reflected representatives.

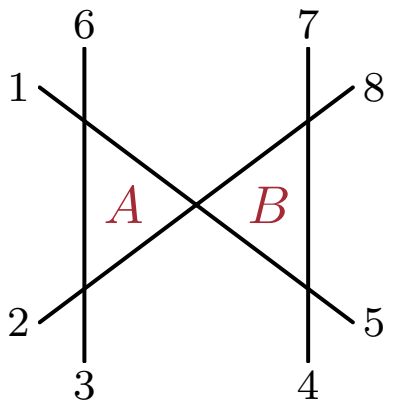
\begin{figure}[htbp]\centering
\RVPanel{.43}{\RVLocal{0}}{The shared-vertex configuration $H$.}
\begin{minipage}{.53\linewidth}\small
The cyclic order is $1,2,3,4,5,8,7,6$.
Up to the symmetries of $H$, adjacent endpoint pairs have representatives
$12$, $23$ and $34$. Their respective classes are
\[
 \{12,58\},\quad\{16,23,45,78\},\quad\{34,67\}.
\]
Here $ij$ denotes the unordered pair $\{i,j\}$.
\end{minipage}
\caption{The boundary labels used in the $(3,3)$-v comparison. Vertex labels and checkerboard shading are suppressed.}
\label{fig:33v-redundancy}
\end{figure}

We use the classification of $O_4$ in Appendix~\ref{app:O4}. Its first 82 diagrams are forests. The remaining two have one quadrilateral with two boundary legs at each vertex. A boundary leg here is an edge from an internal vertex to the boundary with no intervening vertex. We divide the diagrams according to whether they contain a vertex-free boundary-to-boundary arc.

First suppose that $R$ contains such an arc $\gamma$.
View $R$ in the complementary disc $\Delta_R$ of the
disc containing $H$.
The arc $\gamma$ divides $\Delta_R$ into two regions,
each bounded by $\gamma$ and a portion of
$\partial\Delta_R$.
The six boundary points other than the endpoints of
$\gamma$ are distributed between these two regions.

In either region, consider the part of $R$ excluding
$\gamma$ and its endpoints.
Let $v$, $b$ and $e$ denote its numbers of interior
vertices, boundary points and edges, respectively.
Each interior vertex has degree four, and each boundary
point has degree one.
Counting edge ends gives $4v+b=2e$, so $b$ is even.
The region with fewer boundary points contains at most
three of the six remaining points, and hence contains
zero or two.

Figure~\ref{fig:33v-complementary-arc} illustrates the
case where the smaller region contains two boundary points.

\begin{figure}[htbp]
\centering
\begin{tikzpicture}[
  scale=.72,
  line cap=round,
  line join=round,
  strand/.style={black,line width=.8pt},
  boundary/.style={gray,dashed,line width=.6pt},
  gammaarc/.style={green!45!black,line width=1.3pt},
  every node/.style={font=\small}
]

% Left panel: H and its complementary disc on the sphere.
\begin{scope}
  \fill[blue!5] (-3,-2.9) rectangle (3,2.7);

  % Region cut off by gamma, containing boundary points 2 and 3.
  \fill[orange!15]
    (150:1.9)
    .. controls (-3.3,.2) and (-2.1,-3.6) ..
    (300:1.9)
    arc[start angle=300,end angle=150,radius=1.9]
    -- cycle;

  % The disc containing H.
  \fill[white] (0,0) circle (1.9);
  \draw[boundary] (0,0) circle (1.9);

  \coordinate (Ltop) at (-.83,.62);
  \coordinate (Lbot) at (-.83,-.62);
  \coordinate (O) at (0,0);
  \coordinate (Rtop) at (.83,.62);
  \coordinate (Rbot) at (.83,-.62);

  % Two triangles sharing a vertex.
  \draw[strand] (Ltop)--(Lbot)--(O)--cycle;
  \draw[strand] (Rtop)--(Rbot)--(O)--cycle;

  % Eight boundary legs, with the original cyclic labelling.
  \draw[strand]
    (Ltop)--(150:1.9)
    (Ltop)--(120:1.9)
    (Lbot)--(210:1.9)
    (Lbot)--(240:1.9)
    (Rtop)--(60:1.9)
    (Rtop)--(30:1.9)
    (Rbot)--(300:1.9)
    (Rbot)--(330:1.9);

  \foreach \v in {Ltop,Lbot,O,Rtop,Rbot}{
    \fill (\v) circle (1.6pt);
  }

  \node at (-.53,0) {$A$};
  \node at (.53,0) {$B$};
  \node[font=\large] at (0,1.05) {$H$};

  \foreach \lab/\ang in
    {1/150,2/210,3/240,4/300,5/330,8/30,7/60,6/120}{
    \fill (\ang:1.9) circle (1.2pt);
    \node[fill=white,inner sep=1pt]
      at (\ang:1.68) {$\lab$};
  }

  \draw[gammaarc]
    (150:1.9)
    .. controls (-3.3,.2) and (-2.1,-3.6) ..
    (300:1.9);

  \node[text=green!45!black,font=\large]
    at (-2.65,-.65) {$\gamma$};
  \node at (-1.55,-1.7) {$b=2$};
  \node at (2.35,0) {$b=4$};
  \node[font=\large] at (2.4,1.1) {$\infty$};
  \node at (0,-3.25) {(a) The arc outside $H$};
\end{scope}

% Right panel: the same complementary disc drawn as a bounded disc.
\begin{scope}[xshift=7cm]
  \fill[blue!5] (0,0) circle (2.35);

  \fill[orange!15]
    (150:2.35)
    arc[start angle=150,end angle=300,radius=2.35]
    -- cycle;

  \draw[boundary] (0,0) circle (2.35);
  \draw[gammaarc] (150:2.35)--(300:2.35);

  \foreach \lab/\ang in
    {1/150,2/210,3/240,4/300,5/330,8/30,7/60,6/120}{
    \fill (\ang:2.35) circle (1.4pt);
    \node at (\ang:2.6) {$\lab$};
  }

  \node[font=\large] at (.45,1.35) {$\Delta_R$};
  \node at (.65,.45) {$b=4$};
  \node at (-1,-1.15) {$b=2$};
  \node[text=green!45!black,font=\large]
    at (-.15,-.25) {$\gamma$};
  \node at (0,-3.25) {(b) The complementary disc};
\end{scope}

\end{tikzpicture}
\caption{A vertex-free arc $\gamma$ joining boundary points
$1$ and $4$ of the $(3,3)$-v configuration.
The complementary disc is shown outside $H$ in (a),
including the point at infinity, and as a bounded disc
in (b).
Excluding the endpoints of $\gamma$, the two regions
contain two and four boundary points.
The remaining strands of $R$ are omitted.
The colours distinguish the regions and do not indicate
checkerboard shading.}
\label{fig:33v-complementary-arc}
\end{figure}
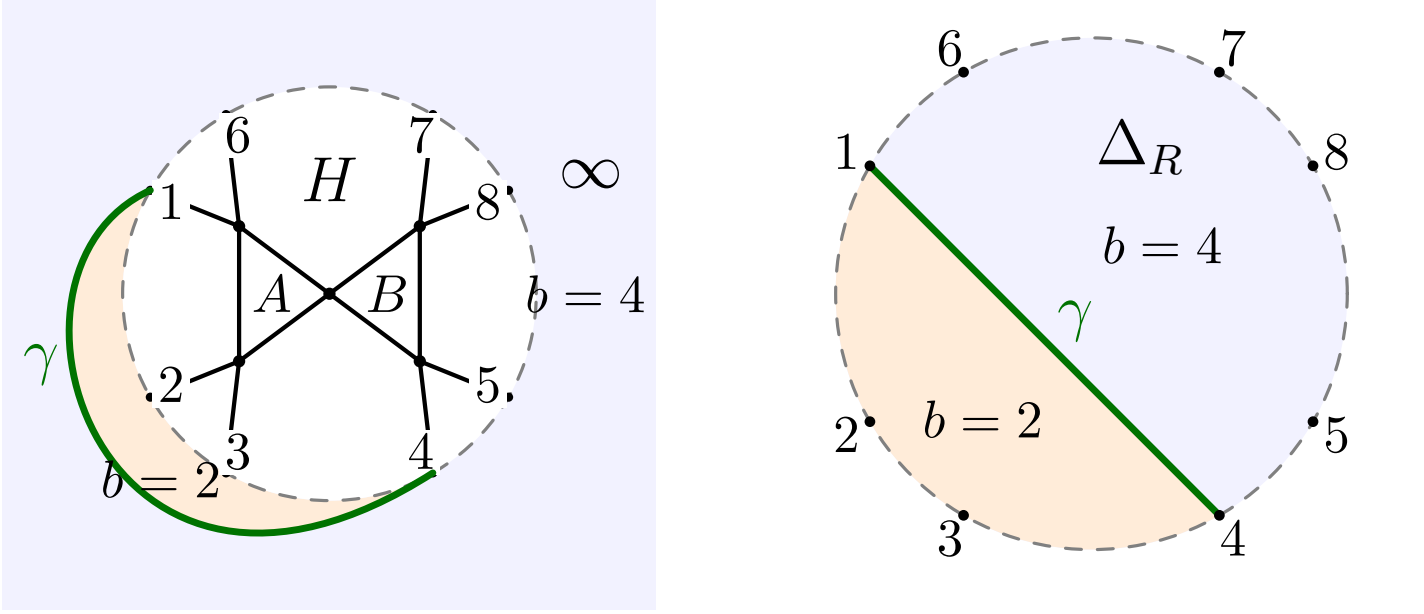

If it contains zero, the region is empty. Otherwise a closed component would contain a reducible face by the closed-diagram small-face result, contradicting reducedness. Thus $\gamma$ joins consecutive boundary points and cuts off an empty region. If it contains two, its interior is a reduced 1-box. By Lemma~\ref{cons:lowbox}, with closed components excluded as above, this is the single arc $\OnOne$ in $O_1$. Choose this inner arc instead of $\gamma$. Its endpoints are consecutive on the original boundary, and it again cuts off an empty region. In both cases it therefore suffices to consider the three attachments in Figure~\ref{fig:33v-arcs}.

\begin{figure}[htbp]\centering
\RVPanel{.32}{\RVLocal{1}}{$12$: $C$ is a 2-face.}\hfill
\RVPanel{.32}{\RVLocal{2}}{$23$: $C$ is a 1-face.}\hfill
\RVPanel{.32}{\RVLocal{3}}{$34$: $C$ is a 3-face.}
\caption{The three vertex-free arc attachments. The added arc is green.}
\label{fig:33v-arcs}
\end{figure}

For $12$, the faces $B,C$ are disjoint, and $A,C$ share an edge, forming a $(3,2)$-e configuration. The disjoint comparison and the $(3,2)$-e comparison give
$
 Z(D_B)=Z(D_C)=Z(D_A).
$ 

For $23$, the faces $B,C$ are disjoint, and $A,C$ share a vertex, forming the $(3,1)$ configuration. The disjoint comparison and the $(3,1)$ comparison give
$
 Z(D_B)=Z(D_C)=Z(D_A).
$

For $34$, the 3-face $C$ shares an edge with each of $A,B$. The two $(3,3)$-e comparisons give $Z(D_A)=Z(D_C)=Z(D_B)$. This treats every diagram containing a vertex-free arc.

Now suppose that $R$ contains no vertex-free arc.
First consider the case where $R$ is a forest.
Let $v$, $e$ and $c$ denote its numbers of internal
vertices, edges and connected components, respectively.
The internal vertices have degree four, and the eight
boundary points have degree one, so counting edge ends
gives $4v+8=2e$.
Since each tree has one fewer edge than vertices,
the forest has $e=(v+8)-c$ edges.
Substituting this into the degree sum yields
\[
4v+8=2(v+8-c),
\qquad\text{hence}\qquad
v+c=4.
\]
Every component contains an internal vertex, so $c\le v$. Consequently $(v,c)=(2,2)$ or $(3,1)$. In the first case, the two components are single four-valent vertices with four boundary legs each. In the second case, the three internal vertices form a path, whose endpoints each have three boundary legs. Figure~\ref{fig:33v-forests} displays representatives of both structures.

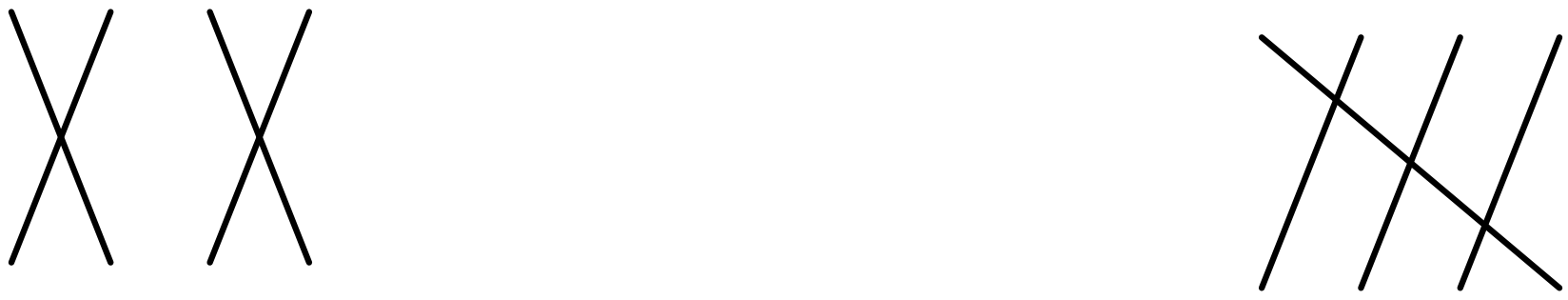
\begin{figure}[htbp]\centering
\begin{minipage}{.47\linewidth}\centering
\scalebox{3}{$\OnFour{45}$}\par\medskip\small Two connected components, each a four-valent vertex.
\end{minipage}\hfill
\begin{minipage}{.47\linewidth}\centering
\scalebox{3}{$\OnFour{71}$}\par\medskip\small Three internal vertices forming a path.
\end{minipage}
\caption{The two forest structures without a vertex-free arc. These representatives are diagrams 45 and 71 of Appendix~\ref{app:O4}.}
\label{fig:33v-forests}
\end{figure}

In the two-component case, each component meets a consecutive block of four boundary points, since the components are disjoint in the disc. In the connected case, the rest of the tree meets an endpoint vertex through its fourth edge. Its three boundary legs are therefore consecutive. In either case, choose a vertex $q$ with three consecutive boundary legs, ending at $u,w,t$ in this order. The two regions between these three edges and the boundary are empty, as illustrated in Figure~\ref{fig:33v-three-legs}.

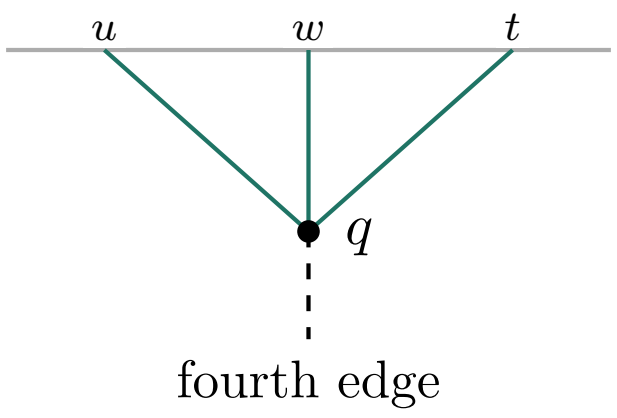
\begin{figure}[htbp]
\centering
\begin{tikzpicture}[x=1cm,y=.75cm,line width=.8pt]
  \draw[gray!65] (-2,1.5)--(2,1.5);
  \coordinate (Q) at (0,-.1);

  \draw[good,thick] (-1.35,1.5)--(Q)--(0,1.5);
  \draw[good,thick] (1.35,1.5)--(Q);
  \draw[dashed] (Q)--(0,-1.05);

  \node[vtx] at (Q) {};
  \node[right] at (.1,-.15) {$q$};

  \node[leg,above] at (-1.35,1.5) {$u$};
  \node[leg,above] at (0,1.5) {$w$};
  \node[leg,above] at (1.35,1.5) {$t$};

  \node[font=\small,below] at (0,-1.05) {fourth edge};
\end{tikzpicture}
\caption{Three consecutive boundary legs at $q$.
The regions bounded by the legs to $u,w$ and to $w,t$,
together with the intervening boundary arcs, contain
no other part of the diagram.
The endpoint pairs $\{u,w\}$ and $\{w,t\}$ share $w$.
In the four-leg case, the fourth edge also ends on
the boundary.}
\label{fig:33v-three-legs}
\end{figure}

The pairs $\{u,w\}$ and $\{w,t\}$ cannot both be $34$ or $67$, since those two pairs are disjoint. Choose a pair other than $34$ and $67$. After applying a symmetry of $H$, its attachment has type $12$ or $23$ in Figure~\ref{fig:33v-vertices}.

\begin{figure}[htbp]\centering
\RVPanel{.32}{\RVLocal{4}}{$12$: $C$ is a 3-face.}\hfill
\RVPanel{.32}{\RVLocal{5}}{$23$: $C$ is a 2-face.}\hfill
\RVPanel{.32}{\RVLocal{6}}{$34$: $C$ is a 4-face.}
\caption{Attachments through a four-valent vertex. If a chosen pair gives the orange 4-face, another pair provides one of the first two comparisons.}
\label{fig:33v-vertices}
\end{figure}

For $12$, the 3-face $C$ is disjoint from $B$ and shares an edge with $A$. The disjoint and $(3,3)$-e comparisons give $Z(D_B)=Z(D_C)=Z(D_A)$. For $23$, the 2-face $C$ is disjoint from $B$ and shares a vertex with $A$. The disjoint and $(3,2)$-v comparisons give the same equality. This treats every forest without a vertex-free arc.

Finally consider either quadrilateral diagram, numbered 83 and 84 in Appendix~\ref{app:O4}. Its four internal vertices each meet a pair of consecutive boundary points, with an empty region between the two boundary legs. These four pairs are pairwise disjoint. At most two are $34$ or $67$, so at least two give one of the comparisons just established. Figure~\ref{fig:33v-quads} shows the two boundary pairing patterns. We use a suitable boundary attachment without reducing the internal quadrilateral.

\begin{figure}[htbp]
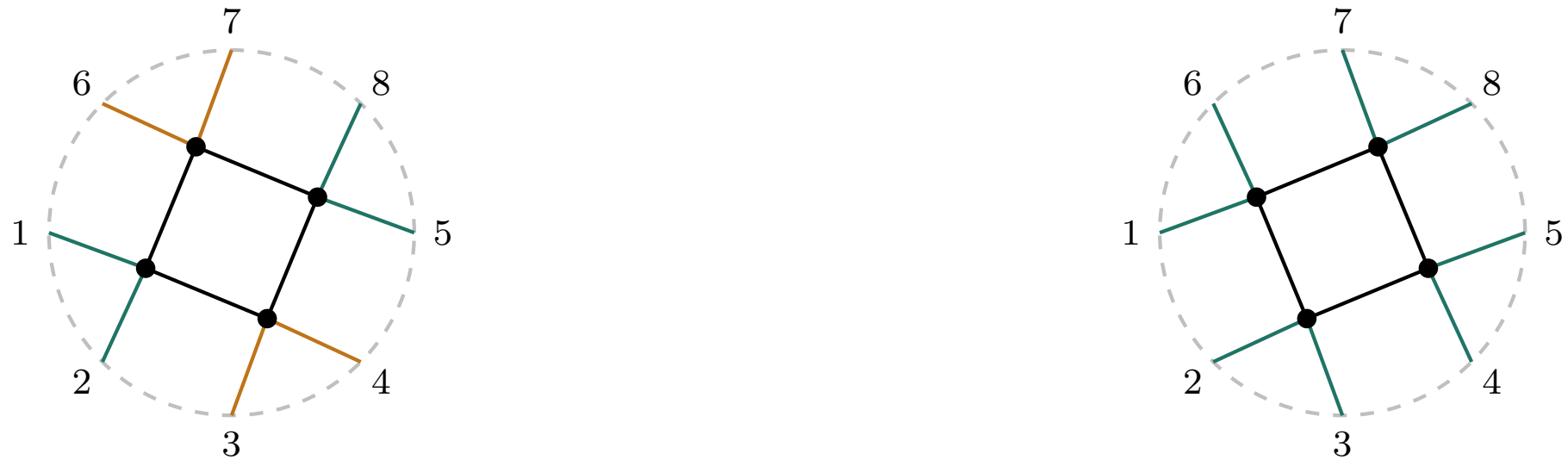
\centering
\RVPanel{.48}{\RVQuad{0}}{$12,34,58,67$}\hfill
\RVPanel{.48}{\RVQuad{1}}{$23,45,78,16$}
\caption{The two quadrilateral attachment patterns. The numbers are boundary labels of $H$. Every green pair gives a $(3,3)$-e or $(3,2)$-v comparison with a disjoint intermediate comparison.}
\label{fig:33v-quads}
\end{figure}

Thus $Z(D_A)=Z(D_B)$ for every $R\in O_4$. The $(3,3)$-v configuration introduces no additional finite test. All first reductions now have the same value. Lemma~\ref{lem:first-step} completes the induction, and Lemma~\ref{lem:path-invariance} gives the pre-partition function. Conversely, a pre-partition function satisfies every stated test.
\end{proof}
\begin{remark}
The finiteness of the consistency conditions can also be understood
from an algebraic viewpoint. With finitely many labels and reduction
rules, path independence amounts to a possibly infinite system of
polynomial equations in finitely many local coefficients.
By the Hilbert basis theorem, the ideal generated by these equations
is finitely generated, so only finitely many equations are needed.
This abstract argument does not identify which equations suffice.
Theorem~\ref{cons:graphthm} gives an explicit finite set of
consistency tests through the analysis of local configurations.
\end{remark}

The shared-vertex argument uses only the local skein rules. After its
elimination, all nontrivial coordinate comparisons occur in 2-boxes and
3-boxes. Scalar closures provide the remaining normalizations.

\label{subsec:full-partition}
\begin{theorem}\label{thm:partition-function}\label{cons:partition}
Let $Z$ be a pre-partition function on $\mathcal P(L_2)$.
If the closure of the 2-box identity $\OnTwo{1}$ with every compatible labelled
$B\in O_2$, in either shading, has the same $Z$-value before and after
Move $4$, then $Z$ is a partition function.
\label{eq:unit-finite-test}
\end{theorem}
\begin{proof}
Let $D$ be a closed diagram, and let $U$ be a disc in which
Move $4$ is to be applied. Denote the resulting linear
combination by $\widetilde D$, and let $B$ be the complementary
$2$-box outside $U$. We prove $Z(D)=Z(\widetilde D)$ by induction
on the number of faces of $D$, allowing disconnected diagrams
throughout.

Suppose that $D$ has $n$ faces and that the assertion holds
for every closed diagram with fewer than $n$ faces, at every
choice of identity-resolution disc.
If $B$ has no internal face, then $B$ is a compatibly labelled
element of $O_2$, so the assertion follows from the hypothesis
of the theorem.

Otherwise, closed components of $B$ can first be evaluated by canonical
reductions, without affecting $U$. The remaining diagram is
either a connected $2$-box or the disjoint union of two connected
$1$-boxes. If it has an internal face, Lemma~\ref{cons:lowbox}
therefore provides an internal $k$-face $F_k$ for some
$k\in\{1,2,3\}$. In either case, the chosen reduction
takes place outside $U$.

Apply this reduction to $D$, obtaining a linear combination
$D_k$, each term of which has fewer than $n$ faces, with the
disc $U$ unchanged. Applying Move $4$ at $U$ gives
$\widetilde D_k$. Since the two moves have disjoint
supports, $\widetilde D_k$ is also obtained by applying the
same reduction to $\widetilde D$.

As $Z$ is a pre-partition function, the chosen reduction
preserves its value on both sides. The induction hypothesis,
applied term by term to $D_k$, therefore gives
\begin{equation}\label{eq:unit-induction}
  Z(D)=Z(D_k)=Z(\widetilde D_k)=Z(\widetilde D).
\end{equation}
Thus $Z$ is invariant under Move $4$ in every closed
context, and hence is a partition function.
\end{proof}
The proof does not require identity resolution to preserve connectedness
or to decrease the face number. The decreasing quantity is the face
number of the common complement. In particular, the finite tests include
the cases in which the two strands belong to different closed components
and an inserted 2-box joins them.

\subsection{Algebraic consistency equations}
\label{sec:algebraic-consistency}

Let $P$ be a YBR planar algebra whose $2$-box space $P_{2,+}$
carries a fusion bialgebra structure $B$.
Using the minimal idempotents $p_i$ of $B$, we translate the
graphical consistency conditions into algebraic equations
in the quantum dimensions $d_i$, the convolution coefficients
$N_{ij}^{k}$, and the triangle coefficients $c$ and $\widetilde c$.
The latter two tables correspond to the two triangle orientations,
both expressed in the $P_{2,+}$ labels.

Let $H$ be one of the local overlap diagrams listed in
\eqref{cons:list}, with $2r$ boundary points.
Denote by $H_m$ and $H_k$ the outputs of the two canonical
moves on the overlapping faces.
Previously, these outputs were paired with a complementary
diagram to form closed diagrams. We now compare them directly
as $r$-box diagrams.

\begin{theorem}\label{thm:local-to-standard-consistency}
Let $Z$ be given by \eqref{eq:average-evaluation}. Suppose that
$H_m=H_k$ as formal linear combinations of labelled diagrams
for every overlapping pair of canonical first moves $m,k$.
Then $Z$ is a pre-partition function.
Conversely, if $Z$ is a pre-partition function, then
\[
Z\bigl(\langle H_m,B\rangle\bigr)
=Z\bigl(\langle H_k,B\rangle\bigr)
\]
for every compatible complementary diagram $B$.
If $Z$ is induced by a nondegenerate planar algebra $P$ generated
by the prescribed labels, then $H_m=H_k$ in $P_{r,\pm}$.
\end{theorem}
\begin{proof}
We prove the first assertion by induction on the number of faces
of a closed diagram $D$. The terminal case is immediate.
Suppose that $D$ has $N$ faces and that canonical evaluation is path
independent for every closed diagram with fewer than $N$ faces.
Let $m,k$ be two canonical first moves on $D$.
Every term of $D_m,D_k$ has fewer than $N$ faces, so its subsequent
canonical evaluation is path independent by the induction hypothesis.

If the supports are disjoint, each move remains available after
the other, and $(D_m)_k=(D_k)_m$. The induction hypothesis gives
$Z(D_m)=Z(D_k)$.
If the supports overlap, write $D=\langle H,B\rangle$.
Locality gives $D_m=\langle H_m,B\rangle$ and
$D_k=\langle H_k,B\rangle$. The assumed local equality and linearity imply
\[
Z(D_m)=Z\bigl(\langle H_m,B\rangle\bigr)
      =Z\bigl(\langle H_k,B\rangle\bigr)=Z(D_k).
\]
These expressions are understood linearly when a move has
more than one output term. Thus every canonical first move gives
the same value. Together with the induction hypothesis for subsequent
evaluation, this proves path independence for $D$.
Lemma~\ref{lem:path-invariance} implies that $Z$ is a pre-partition function.

Conversely, invariance under canonical moves gives
\[
Z\bigl(\langle H_m,B\rangle\bigr)
=Z\bigl(\langle H,B\rangle\bigr)
=Z\bigl(\langle H_k,B\rangle\bigr)
\]
for every compatible $B$. Hence $H_m-H_k$ pairs to zero with every
such diagram. Since the given labels generate $P$, the labelled complementary
diagrams span the complementary box space.
Nondegeneracy of the closure pairing therefore implies
$H_m-H_k=0$ in $P_{r,\pm}$.
\end{proof}

Recall the diagrams $T(i,j,k)$ and $W(i,j,k)$ and the three
face-free families $E_1(a,b)$, $E_2(a,b)$ and $E_3(a,b)$
from Section~\ref{sec:TR-definitions}.
The active triangles
$\{T(x,y,z):N_{zy}^{x}\ne0\}$
form the 3-box basis in the YBR setting under consideration. The identity resolutions read
\begin{equation}\label{cons:E1}
E_1(x,y)=\sum_{z\in I}T(x,y,z),\quad E_2(x,z)=\sum_{y\in I}T(x,y,z),\quad E_3(y,z)=\sum_{x\in I}T(x,y,z).
\end{equation}
The triangle relations for $T(i,j,k)$ and $W(i,j,k)$ take the form
\begin{align}
 T(i,j,k)&=\sum_{a,b\in I}\bigl(
 c_{1,ab}^{ijk}E_1(a,b)+c_{2,ab}^{ijk}E_2(a,b)
 +c_{3,ab}^{ijk}E_3(a,b)\bigr),\label{cons:TR}\\
 W(i,j,k)&=\sum_{a,b\in I}\bigl(
 \tc_{1,ab}^{ijk}E_1(a,b)+\tc_{2,ab}^{ijk}E_2(a,b)
 +\tc_{3,ab}^{ijk}E_3(a,b)\bigr).\label{cons:oppositeTR}
\end{align}

The labelled diagrams and their boundary conventions are fixed in \Cref{fig:triangle-convention,fig:E-diagrams}. Their expansions are shown in \eqref{pre:TR-native-picture} and \eqref{pre:TR-W-picture}.

\begin{theorem}\label{cons:algthm}\label{thm:corrected-33e}
Let $P$ be a YBR planar algebra and $B=P_{2,+}$. Write $N_{ij}^{k}$
and $d_i$ for its fusion coefficients and quantum dimensions, and
$c_{\alpha,lm}^{ijk},\widetilde c_{\alpha,lm}^{ijk}$ for the coefficients
of the proposed TR relations, where $\alpha=1,2,3$ and all labels lie in $I$.
The TR consistency tests of Theorem~\ref{cons:graphthm},
together with the identity-resolution tests of
Theorem~\ref{thm:partition-function}, are equivalent
to the following polynomial equations:
\begin{gather}
d_i=d_{\bar i},
\label{cons:dim-dual}\\[4pt]
N_{ij}^{0}=\delta_{i,\bar j}d_i,
\label{cons:vacuum}\\[4pt]
N_{ij}^{k}=N_{\bar j\,\bar i}^{\bar k},
\label{cons:duality}\\[4pt]
\sum_{s\in I}N_{ij}^{s}d_s=d_i d_j,
\label{cons:dimension}\\[4pt]
\sum_{s\in I}N_{ij}^{s}N_{sk}^{\ell}
=\sum_{s\in I}N_{is}^{\ell}N_{jk}^{s},
\label{cons:associativity}\\[4pt]
\delta^2=\sum_{i\in I}d_i,
\label{cons:loop-completion}\\[4pt]
c_{1,0\bar s}^{ijk}+c_{2,0s}^{ijk}+c_{3,\bar s s}^{ijk}
=\delta_{i0}\delta_{k\bar j}\delta_{sk},
\label{cons:31}\\[4pt]
\delta\sum_m\tc_{1,sm}^{ijk}d_{\bar m}
+\delta\sum_q\tc_{2,sq}^{ijk}d_q
+\delta\sum_{r,t}\tc_{3,rt}^{ijk}N_{tr}^{s}
=d_kN_{ji}^{s},
\label{cons:31-W}\\[4pt]
\Bigl[
\delta_{x\ell}
(c_{1,xy}^{ijk}+c_{2,xz}^{ijk}+c_{3,yz}^{ijk})
-\delta_{xi}\delta_{yj}\delta_{zk}\delta_{i\ell}
\Bigr]N_{zy}^{x}=0,
\label{cons:32v}\\[4pt]
\begin{gathered}
\Biggl[
\sum_a
(\tc_{1,xa}^{ijk}N_{\ell a}^{y}
+\tc_{2,xa}^{ijk}N_{a\bar\ell}^{z})
+\delta\sum_{r,s}\tc_{3,rs}^{ijk}
(\tc_{1,xy}^{rs\ell}+\tc_{2,xz}^{rs\ell}+\tc_{3,yz}^{rs\ell})
\\
-\sum_sN_{\ell k}^{s}
(\tc_{1,xy}^{ijs}+\tc_{2,xz}^{ijs}+\tc_{3,yz}^{ijs})
\Biggr]N_{zy}^{x}=0,
\end{gathered}
\label{cons:32v-W}\\[4pt]
\sum_{a,b}c_{1,ab}^{ijk}N_{a\bar b}^{s}
+\sum_a c_{2,as}^{ijk}d_a
+\sum_a c_{3,as}^{ijk}d_a
=\delta_{sk}N_{i\bar j}^{k},
\label{cons:32e}\\[4pt]
\delta\sum_m\tc_{1,sm}^{ijk}\delta_{sm}
+\delta\sum_q\tc_{2,sq}^{ijk}\delta_{q0}
+\delta\sum_t\tc_{3,st}^{ijk}\delta_{t0}
=\delta_{jk}N_{ki}^{s},
\label{cons:32e-W}\\[4pt]
\begin{gathered}
\Biggl[
\sum_a c_{1,ay}^{ij\bar k}N_{\ell a}^{x}
+\delta\sum_{p,q}c_{2,pq}^{ij\bar k}
(\tc_{1,xy}^{p\ell\bar q}
+\tc_{2,xz}^{p\ell\bar q}
+\tc_{3,yz}^{p\ell\bar q})
\\
+\sum_s c_{3,ys}^{ij\bar k}N_{\ell s}^{z}
-\delta
(\tc_{1,xy}^{i\ell k}+\tc_{2,xz}^{i\ell k}+\tc_{3,yz}^{i\ell k})
\delta_{jy}
\Biggr]N_{zy}^{x}=0.
\end{gathered}
\label{cons:33e}
\end{gather}
\end{theorem}

\begin{proof}
The $(1,1)$ comparison gives \eqref{cons:dim-dual},
and the $(2,1)$ comparison gives \eqref{cons:dimension}.
The two $(2,2)$ comparisons give
\eqref{cons:vacuum} and \eqref{cons:associativity},
together with the ordinary multiplication identities.
Equation~\eqref{cons:duality} expresses compatibility
of the convolution relation with the half-turn.
Frobenius reciprocity follows from these fusion identities,
as explained in Proposition~\ref{prop:basic-2box-equations}.
Closing the identity resolution gives
\eqref{cons:loop-completion}.

We compute the triangle comparisons below using $T$ and $W$,
both labelled by the minimal idempotents $p_i$ of $B$.

\par\Needspace{7\baselineskip}\medskip\noindent\textbf{$(3,1)$: the $T$ orientation.}

\localstep{First reduce the one-face.}
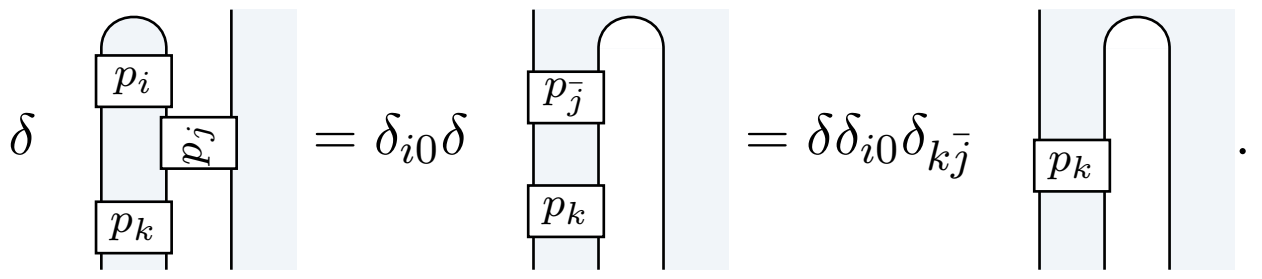
\[
 \delta\rT{1}
 =\delta_{i0}\delta\rProduct{p_{\bar j}}{p_k}
 =\delta\delta_{i0}\delta_{k\bar j}\rOut{p_k}.
\]
\localstep{First reduce the three-face.}
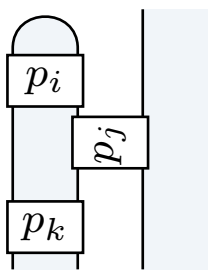
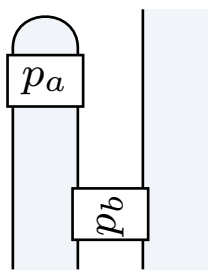
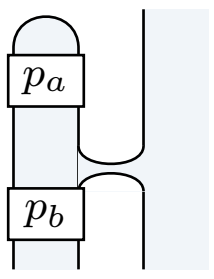
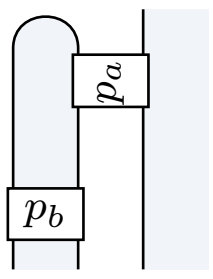
\[
 \delta\rT{1}
 =\delta\sum_{a,b}c_{1,ab}^{ijk}\gE{1}{1}{p_a}{p_b}
 +\delta\sum_{a,b}c_{2,ab}^{ijk}\gE{1}{2}{p_a}{p_b}
 +\delta\sum_{a,b}c_{3,ab}^{ijk}\gE{1}{3}{p_a}{p_b}.
\]
\[
 =\delta\sum_s\left(
 \sum_a c_{1,a\bar s}^{ijk}\delta_{a0}
 +\sum_a c_{2,as}^{ijk}\delta_{a0}
 +\sum_a c_{3,as}^{ijk}\delta_{a\bar s}\right)\rOut{p_s}.
\]
\localstep{Comparison.}
Equation~\eqref{cons:31} follows.

\par\Needspace{7\baselineskip}\medskip\noindent\textbf{$(3,1)$: the $W$ orientation.}

\localstep{First reduce the one-face.}
\[
 \delta^2\gW{4}{p_i}{p_j}{p_k}
 =d_k\delta\wGone
 =\sum_s d_kN_{ji}^{s}\bOut{p_s}.
\]
\localstep{First reduce the three-face.}
\[
 \delta^2\gW{4}{p_i}{p_j}{p_k}
 =\delta^2\sum_{a,b}\tc_{1,ab}^{ijk}\gE{4}{1}{p_a}{p_b}
 +\delta^2\sum_{a,b}\tc_{2,ab}^{ijk}\gE{4}{2}{p_a}{p_b}
 +\delta^2\sum_{a,b}\tc_{3,ab}^{ijk}\gE{4}{3}{p_a}{p_b}.
\]
\[
 =\sum_s\left(
 \delta\sum_m\tc_{1,sm}^{ijk}d_{\bar m}
 +\delta\sum_q\tc_{2,sq}^{ijk}d_q
 +\delta\sum_{r,t}\tc_{3,rt}^{ijk}N_{tr}^{s}\right)\bOut{p_s}.
\]
\localstep{Comparison.}
Equation~\eqref{cons:31-W} follows.

\par\Needspace{7\baselineskip}\medskip\noindent\textbf{$(3,2)$-v: the $T$ orientation.}

All sums below are over the label set $I$.

\localstep{First reduce the two-face.}
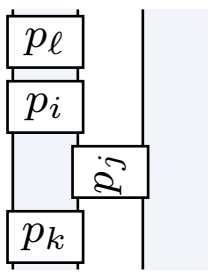
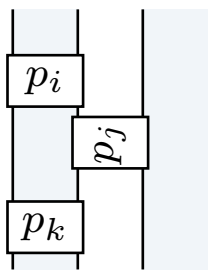
\[
 \rV=\delta_{i\ell}\eT{p_i}{p_j}{p_k}
       =\delta_{i\ell}T(i,j,k).
\]
\localstep{First reduce the three-face.}
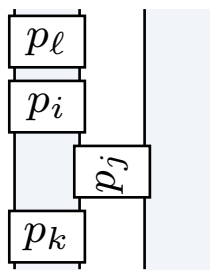
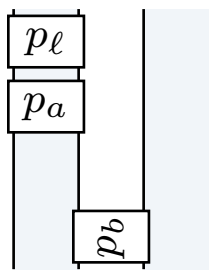
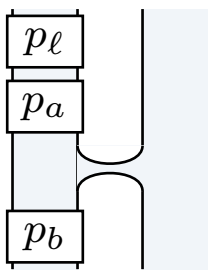
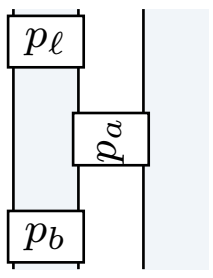
\[
\begin{aligned}
 \rV={}&\sum_{a,b}c_{1,ab}^{ijk}\rVE{1}
       +\sum_{a,b}c_{2,ab}^{ijk}\rVE{2}
       +\sum_{a,b}c_{3,ab}^{ijk}\rVE{3}\\[3pt]
 ={}&\sum_{a,b}c_{1,ab}^{ijk}\delta_{a\ell}E_1(a,b)
      +\sum_{a,b}c_{2,ab}^{ijk}\delta_{a\ell}E_2(a,b)
      +\sum_{a,b}c_{3,ab}^{ijk}T(\ell,a,b)\\
 ={}&\sum_{x,y,z}\delta_{x\ell}
       (c_{1,xy}^{ijk}+c_{2,xz}^{ijk}+c_{3,yz}^{ijk})T(x,y,z).
\end{aligned}
\]
\localstep{Comparison.}
Equation~\eqref{cons:32v} follows.

\par\Needspace{7\baselineskip}\medskip\noindent\textbf{$(3,2)$-v: the $W$ orientation.}

\localstep{First reduce the two-face.}
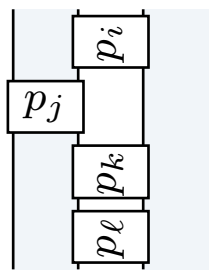
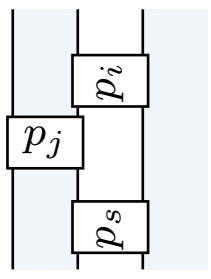
\[
 \delta\wV=\sum_s N_{\ell k}^{s}\eW{p_i}{p_j}{p_s}
     =\sum_s N_{\ell k}^{s}W(i,j,s).
\]
The resulting coefficient of $T(x,y,z)$ is
$\sum_sN_{\ell k}^s(\tc_{1,xy}^{ijs}+\tc_{2,xz}^{ijs}+\tc_{3,yz}^{ijs})$.
\localstep{First reduce the three-face.}
\[
 \delta\wV
 =\delta\sum_{a,b}\tc_{1,ab}^{ijk}\wVE{1}
 +\delta\sum_{a,b}\tc_{2,ab}^{ijk}\wVE{2}
 +\delta\sum_{a,b}\tc_{3,ab}^{ijk}\wVE{3}.
\]
\localstep{Coefficient comparison.}
Equation~\eqref{cons:32v-W} follows.

\par\Needspace{7\baselineskip}\medskip\noindent\textbf{$(3,2)$-e: the $T$ orientation.}

\localstep{First reduce the two-face.}
\[
 \delta\rT{2}
 =\sum_sN_{i\bar j}^{s}\rProduct{p_s}{p_k}
 =N_{i\bar j}^{k}\rOut{p_k}.
\]
\localstep{First reduce the three-face.}
\[
 \delta\rT{2}
 =\delta\sum_{a,b}c_{1,ab}^{ijk}\gE{2}{1}{p_a}{p_b}
 +\delta\sum_{a,b}c_{2,ab}^{ijk}\gE{2}{2}{p_a}{p_b}
 +\delta\sum_{a,b}c_{3,ab}^{ijk}\gE{2}{3}{p_a}{p_b}.
\]
\[
 =\sum_s\left(\sum_{a,b}c_{1,ab}^{ijk}N_{a\bar b}^{s}
 +\sum_a c_{2,as}^{ijk}d_a+\sum_a c_{3,as}^{ijk}d_a\right)\rOut{p_s}.
\]
\localstep{Comparison.}
Equation~\eqref{cons:32e} follows.

\par\Needspace{7\baselineskip}\medskip\noindent\textbf{$(3,2)$-e: the $W$ orientation.}

\localstep{First reduce the two-face.}
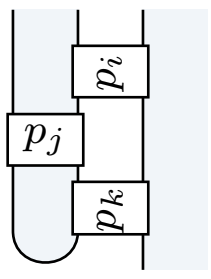
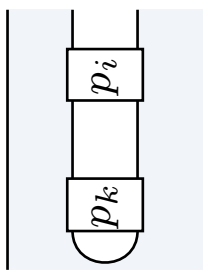
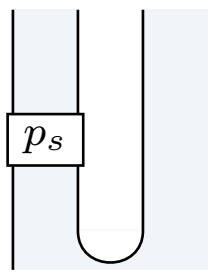
\[
 \delta\gW{3}{p_i}{p_j}{p_k}
 =\delta_{jk}\delta\wProduct
 =\sum_s\delta_{jk}N_{ki}^{s}\bOut{p_s}.
\]
\localstep{First reduce the three-face.}
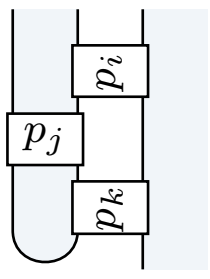
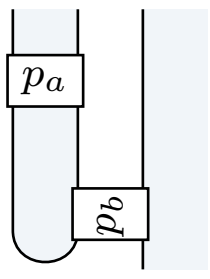
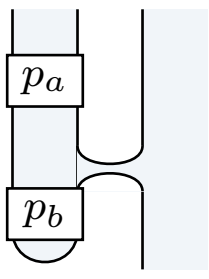
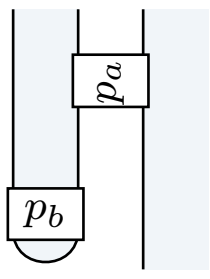
\[
 \delta\gW{3}{p_i}{p_j}{p_k}
 =\delta\sum_{a,b}\tc_{1,ab}^{ijk}\gE{3}{1}{p_a}{p_b}
 +\delta\sum_{a,b}\tc_{2,ab}^{ijk}\gE{3}{2}{p_a}{p_b}
 +\delta\sum_{a,b}\tc_{3,ab}^{ijk}\gE{3}{3}{p_a}{p_b}.
\]
\[
 =\sum_s \delta \left(\sum_m\tc_{1,sm}^{ijk}\delta_{sm}
 +\sum_q\tc_{2,sq}^{ijk}\delta_{q0}
 +\sum_t\tc_{3,st}^{ijk}\delta_{t0}\right)\bOut{p_s}.
\]
\localstep{Coefficient comparison.}
Equation~\eqref{cons:32e-W} follows.

\par\Needspace{7\baselineskip}\medskip\noindent\textbf{$(3,3)$-e: first reduce the right triangle.}

\localstep{Triangle expansion.}
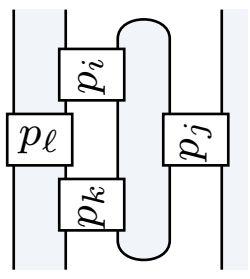
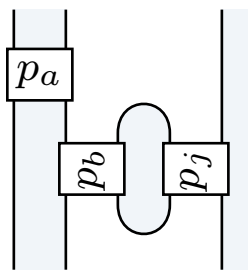
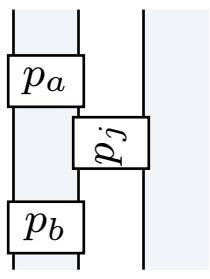
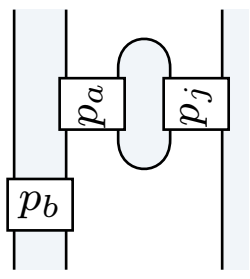
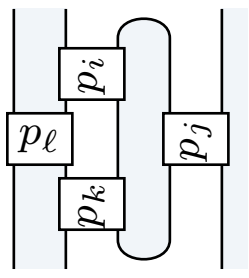
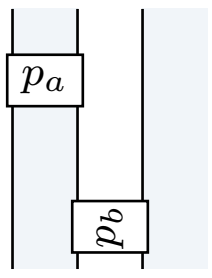
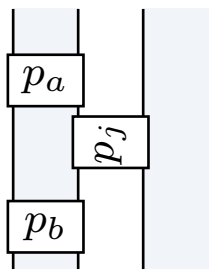
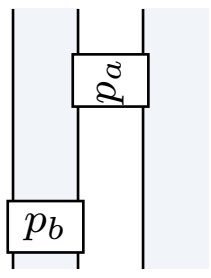
\[
 \delta\eInput
 =\delta\sum_{a,b}\tc_{1,ab}^{i\ell k}\eWraw{1}
 +\delta\sum_{a,b}\tc_{2,ab}^{i\ell k}\eT{p_a}{p_j}{p_b}
 +\delta\sum_{a,b}\tc_{3,ab}^{i\ell k}\eWraw{3}.
\]
\localstep{Primitive multiplication.}
\[
\begin{aligned}
 \delta\eInput=\delta\sum_{a,b}\tc_{1,ab}^{i\ell k}\delta_{jb}\eE{1}{p_a}{p_b}
 +\delta\sum_{a,b}\tc_{2,ab}^{i\ell k}\eT{p_a}{p_j}{p_b}+\delta\sum_{a,b}\tc_{3,ab}^{i\ell k}\delta_{ja}\eE{3}{p_a}{p_b}.
\end{aligned}
\]
\localstep{Expansion in triangle diagrams.}
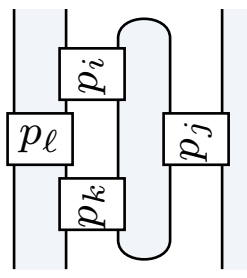
\[
 \delta\eInput
 =\delta\sum_{x,y,z}
 (\tc_{1,xy}^{i\ell k}+\tc_{2,xz}^{i\ell k}+\tc_{3,yz}^{i\ell k})
 \delta_{jy}T(x,y,z).
\]

\par\Needspace{7\baselineskip}\medskip\noindent\textbf{$(3,3)$-e: first reduce the left triangle.}

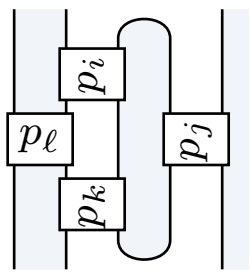
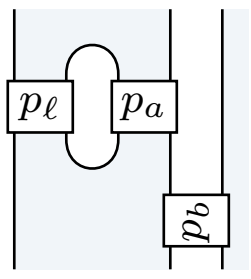
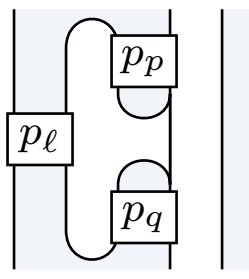
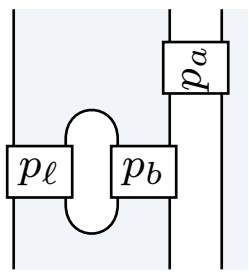
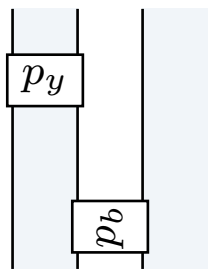
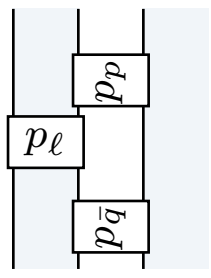
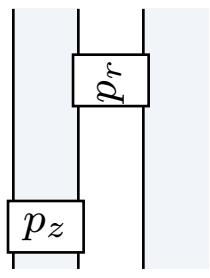
\[
 \delta\eInput
 =\delta\sum_{a,b}c_{1,ab}^{ij\bar k}\eTraw{1}
 +\delta\sum_{p,q}c_{2,pq}^{ij\bar k}\eTmiddle
 +\delta\sum_{a,b}c_{3,ab}^{ij\bar k}\eTraw{3}.
\]
\[
\begin{aligned}
 ={}&\sum_{a,b,y}c_{1,ab}^{ij\bar k}N_{\ell a}^{y}\eE{1}{p_y}{p_b}
 +\delta\sum_{p,q}c_{2,pq}^{ij\bar k}\eW{p_p}{p_\ell}{p_{\bar q}}+\sum_{r,s,z}c_{3,rs}^{ij\bar k}N_{\ell s}^{z}\eE{3}{p_r}{p_z}.
\end{aligned}
\]
The middle term is expanded as
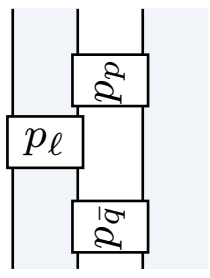
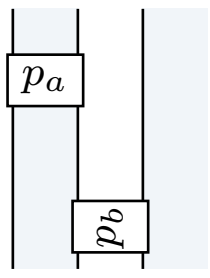
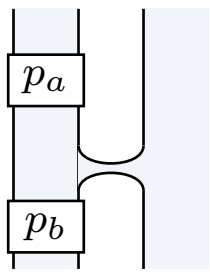
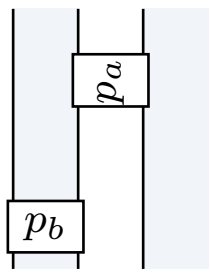
\[
 \delta\eW{p_p}{p_\ell}{p_{\bar q}}
 =\delta\sum_{a,b}\tc_{1,ab}^{p\ell\bar q}\eE{1}{p_a}{p_b}
 +\delta\sum_{a,b}\tc_{2,ab}^{p\ell\bar q}\eE{2}{p_a}{p_b}
 +\delta\sum_{a,b}\tc_{3,ab}^{p\ell\bar q}\eE{3}{p_a}{p_b}.
\]
The final common-edge comparison is recorded in Equation~\eqref{cons:33e}.

In the fixed YBR realization, the active triangles form an orthogonal
basis by Proposition~\ref{prop:triangle-basis}. Thus the coefficient
comparisons above express exactly the equality of the two reduced
outputs. Theorem~\ref{thm:local-to-standard-consistency} relates these local
equalities to the first-step comparisons used in the graphical theorem.

Conversely, in the fixed YBR setting, the weighted equations force
all active coefficients of each comparison to vanish, and inactive
triangles are zero. The corresponding local outputs therefore agree.
The 2-box comparisons follow in the same way from the primitive basis.
This proves the asserted coordinate equivalence. \end{proof}

\begin{remark}
\label{cons:complexity}
Let $n=|I|$. Treating all entries of
$c,\widetilde c,N,d$ and $\delta$ as unknowns gives
$6n^5+n^3+n+1$ scalar variables before imposing
the defining identities and normalizations.
The equations in Theorem~\ref{cons:algthm} form
a polynomial system with $O(n^7)$ scalar equations.
In particular, the $(3,3)$-e equation has seven
free labels and therefore gives $n^7$ indexed equations.

For given coefficient data, each $(3,3)$-e equation
can be evaluated using $O(n^2)$ field operations.
A direct check of all these equations therefore
takes $O(n^9)$ field operations, counting each
field move at unit cost.
The large number of triangle coefficients, coupled through
nonlinear equations, can make direct symbolic solution
computationally demanding even at relatively small ranks.
\end{remark}

\section{Fusion complexes and reconstruction}
\label{sec:fusioncomplex}
\subsection{Fusion complexes and the cohomological criterion}
\label{sec:cohomology}
% ============================================================

Let $B$ be a fusion bialgebra with label set $I$.  Introduce three disjoint
copies
$I_X=\{i_X:i\in I\},
I_Y=\{j_Y:j\in I\},
I_Z=\{k_Z:k\in I\}$.

\begin{definition}
\label{def:fusion-complex}
The \emph{fusion complex} $K(B)$ is the tripartite two-dimensional
simplicial complex with vertex set
$K(B)^{(0)}=I_X\sqcup I_Y\sqcup I_Z$
and faces
\[
[i_X,j_Y,k_Z]\in K(B)^{(2)}
\quad\Longleftrightarrow\quad
N_{kj}^{i}\neq 0.
\]
We abbreviate $[i_X,j_Y,k_Z]$ by $(i,j,k)$ whenever the three copies are
understood.
\end{definition}
\begin{proposition}
The $1$-skeleton of $K(B)$ is the complete tripartite graph on
$I_X\sqcup I_Y\sqcup I_Z$.
\end{proposition}
\begin{proof}
The dimension equation and Frobenius reciprocity give
\[
 \sum_i d_iN_{kj}^{i}=d_kd_j,\qquad
 \sum_kN_{kj}^{i}=d_j,\qquad \sum_jN_{kj}^{i}=d_k.
\]
All right-hand sides are nonzero. Every pair of vertices from different
parts is therefore contained in a face. No edge joins two vertices
in the same part.
\end{proof}

For candidate data, a fusion complex is specified by a face subset of
$I_X\times I_Y\times I_Z$, with the complete tripartite $1$-skeleton
and the prescribed label duality. A realization identifies this complex
with $K(B)$ for a fusion bialgebra $B$.

\begin{figure}[htbp]
\centering
\resizebox{.76\linewidth}{!}{%
% Requires \usepackage{tikz,amsmath,amssymb}.
% Face (i,j,k) has vertices i_X, j_Y, k_Z.
\begin{tikzpicture}[x=1cm,y=1cm,line join=round,line cap=round,
 vertex/.style={circle,inner sep=0pt,minimum size=4.5pt,draw=white,line width=.45pt},
 every node/.style={font=\small}]
\definecolor{fcface}{RGB}{75,139,182}
\definecolor{fcnew}{RGB}{224,151,60}
\definecolor{fcX}{RGB}{154,63,88}
\definecolor{fcY}{RGB}{43,126,118}
\definecolor{fcZ}{RGB}{82,84,151}
\begin{scope}[xshift=0.0cm]
\coordinate (X0) at (-1.75,-0.65);
\coordinate (X1) at (1.75,0.65);
\coordinate (Y0) at (1.4,-0.85);
\coordinate (Y1) at (-1.4,0.85);
\coordinate (Z0) at (0,2.05);
\coordinate (Z1) at (0,-2.05);
\fill[fcface,opacity=.18] (X0)--(Y0)--(Z0)--cycle;
\fill[fcface,opacity=.40] (X0)--(Y1)--(Z1)--cycle;
\fill[fcface,opacity=.40] (X1)--(Y0)--(Z1)--cycle;
\fill[fcface,opacity=.40] (X1)--(Y1)--(Z0)--cycle;
\draw[gray!65,densely dashed,line width=.55pt] (X0)--(Y0);
\draw[black!75,line width=.65pt] (X0)--(Y1);
\draw[black!75,line width=.65pt] (X1)--(Y0);
\draw[black!75,line width=.65pt] (X1)--(Y1);
\draw[gray!65,densely dashed,line width=.55pt] (X0)--(Z0);
\draw[black!75,line width=.65pt] (X0)--(Z1);
\draw[black!75,line width=.65pt] (X1)--(Z0);
\draw[black!75,line width=.65pt] (X1)--(Z1);
\draw[gray!65,densely dashed,line width=.55pt] (Y0)--(Z0);
\draw[black!75,line width=.65pt] (Y0)--(Z1);
\draw[black!75,line width=.65pt] (Y1)--(Z0);
\draw[black!75,line width=.65pt] (Y1)--(Z1);
\node[vertex,fill=fcX] at (X0) {};
\node[anchor=north east,text=fcX,inner sep=4pt] at (X0) {$ 0_X $};
\node[vertex,fill=fcX] at (X1) {};
\node[anchor=south west,text=fcX,inner sep=4pt] at (X1) {$ 1_X $};
\node[vertex,fill=fcY] at (Y0) {};
\node[anchor=north west,text=fcY,inner sep=4pt] at (Y0) {$ 0_Y $};
\node[vertex,fill=fcY] at (Y1) {};
\node[anchor=south east,text=fcY,inner sep=4pt] at (Y1) {$ 1_Y $};
\node[vertex,fill=fcZ] at (Z0) {};
\node[anchor=south,text=fcZ,inner sep=4pt] at (Z0) {$ 0_Z $};
\node[vertex,fill=fcZ] at (Z1) {};
\node[anchor=north,text=fcZ,inner sep=4pt] at (Z1) {$ 1_Z $};
\node[font=\large] at (0,2.75) {$\mathcal T_4$};
\node at (0,-2.85) {$6\text{ vertices},\quad12\text{ edges},\quad4\text{ faces}$};
\node at (0,-3.3) {$ \{000,011,101,110\} $};
\end{scope}
\begin{scope}[xshift=6.2cm]
\coordinate (X0) at (-1.75,-0.65);
\coordinate (X1) at (1.75,0.65);
\coordinate (Y0) at (1.4,-0.85);
\coordinate (Y1) at (-1.4,0.85);
\coordinate (Z0) at (0,2.05);
\coordinate (Z1) at (0,-2.05);
\fill[fcface,opacity=.18] (X0)--(Y0)--(Z0)--cycle;
\fill[fcface,opacity=.40] (X0)--(Y1)--(Z1)--cycle;
\fill[fcface,opacity=.40] (X1)--(Y0)--(Z1)--cycle;
\fill[fcface,opacity=.40] (X1)--(Y1)--(Z0)--cycle;
\fill[fcnew,opacity=.66] (X1)--(Y1)--(Z1)--cycle;
\draw[gray!65,densely dashed,line width=.55pt] (X0)--(Y0);
\draw[black!75,line width=.65pt] (X0)--(Y1);
\draw[black!75,line width=.65pt] (X1)--(Y0);
\draw[black!75,line width=.65pt] (X1)--(Y1);
\draw[gray!65,densely dashed,line width=.55pt] (X0)--(Z0);
\draw[black!75,line width=.65pt] (X0)--(Z1);
\draw[black!75,line width=.65pt] (X1)--(Z0);
\draw[black!75,line width=.65pt] (X1)--(Z1);
\draw[gray!65,densely dashed,line width=.55pt] (Y0)--(Z0);
\draw[black!75,line width=.65pt] (Y0)--(Z1);
\draw[black!75,line width=.65pt] (Y1)--(Z0);
\draw[black!75,line width=.65pt] (Y1)--(Z1);
\draw[fcnew!80!black,line width=.9pt] (X1)--(Y1)--(Z1)--cycle;
\node[vertex,fill=fcX] at (X0) {};
\node[anchor=north east,text=fcX,inner sep=4pt] at (X0) {$ 0_X $};
\node[vertex,fill=fcX] at (X1) {};
\node[anchor=south west,text=fcX,inner sep=4pt] at (X1) {$ 1_X $};
\node[vertex,fill=fcY] at (Y0) {};
\node[anchor=north west,text=fcY,inner sep=4pt] at (Y0) {$ 0_Y $};
\node[vertex,fill=fcY] at (Y1) {};
\node[anchor=south east,text=fcY,inner sep=4pt] at (Y1) {$ 1_Y $};
\node[vertex,fill=fcZ] at (Z0) {};
\node[anchor=south,text=fcZ,inner sep=4pt] at (Z0) {$ 0_Z $};
\node[vertex,fill=fcZ] at (Z1) {};
\node[anchor=north,text=fcZ,inner sep=4pt] at (Z1) {$ 1_Z $};
\node[font=\large] at (0,2.75) {$\mathcal T_5$};
\node at (0,-2.85) {$6\text{ vertices},\quad12\text{ edges},\quad5\text{ faces}$};
\node at (0,-3.3) {$ \{000,011,101,110,\boldsymbol{111}\} $};
\end{scope}
\fill[fcface,opacity=.4] (1.1,-4.0) rectangle (1.45,-3.82);
\node[anchor=west] at (1.55,-3.91) {Faces of $\mathcal T_4$};
\fill[fcnew,opacity=.66] (4.35,-4.0) rectangle (4.7,-3.82);
\node[anchor=west] at (4.8,-3.91) {Additional face $111$};
\end{tikzpicture}
}
\caption{The rank-two fusion complexes. Blue triangles are the four faces
of $\mathcal T_4$. The orange triangle is the additional face of
$\mathcal T_5$. Dashed edges lie at the back of the octahedral drawing.}
\label{fig:fusion-complex-n2}
\end{figure}
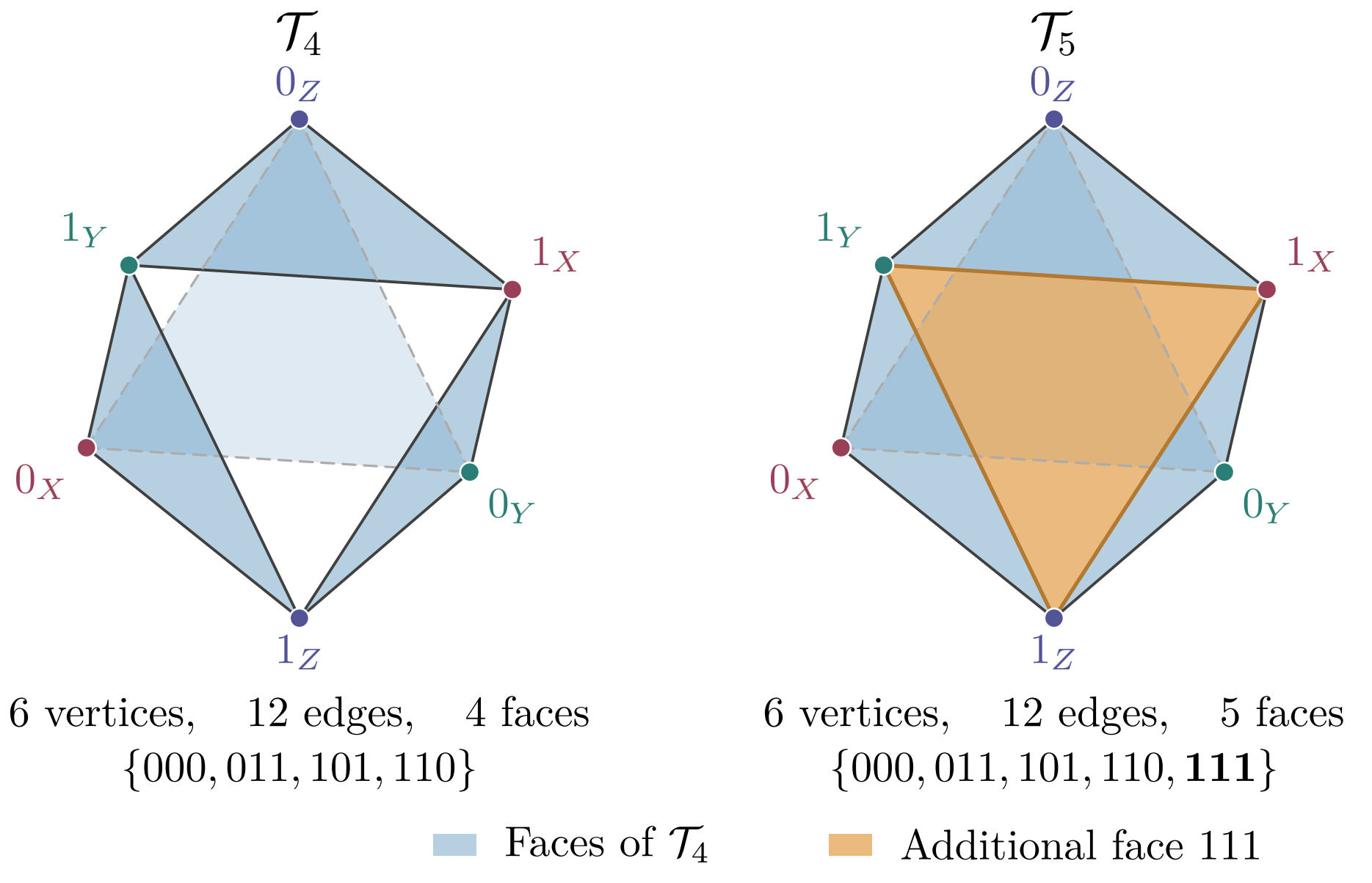

For $I=\{0,1\}$, the complete tripartite $1$-skeleton is the octahedral
graph. Figure~\ref{fig:fusion-complex-n2} shows the two rank-two examples
$\mathcal T_4=\{000,011,101,110\}$ and
$\mathcal T_5=\mathcal T_4\cup\{111\}$, where $ijk$ denotes the face
$[i_X,j_Y,k_Z]$.

\begin{remark}
The fusion complex records only the zero--nonzero support of the fusion
coefficients.  It need not determine the numerical fusion bialgebra.  A
fixed complex can have no realization, finitely many realizations, or a
parameter family of realizations.
\end{remark}

\begin{definition}
The fusion support and face number are
$\supp(B)=\{(i,j,k):N_{kj}^{i}\neq0\},
f(B)=|\supp(B)|=|K(B)^{(2)}|$.
\end{definition}

For $i_X\in I_X$, the link $\operatorname{Lk}_{K(B)}(i_X)$ is a
bipartite graph with parts $I_Y,I_Z$.  As the fixed label
varies, these links recover the family of fusion graphs used in the
exchange-relation theory.

\begin{proposition}
\label{prop:links-fusion-graphs}\label{thm:ER-detector}
For each label $i\in I$, the coordinate fusion graph associated with $p_i$
is canonically isomorphic, up to the fixed index convention, to the link
$\operatorname{Lk}_{K(B)}(i_X)$.
\end{proposition}

\begin{proof}
The edge $(j_Y,k_Z)$ belongs to the link at $i_X$ exactly when
$N_{kj}^{i}\ne0$. Weighted Frobenius reciprocity gives
$d_iN_{kj}^{i}=d_kN_{i\bar j}^{k}$. Since both quantum dimensions are nonzero,
this is the edge from input $\bar j$ to output $k$ in the fusion graph
for $p_i$. Relabel the $Y$ vertices by duality. The other two coordinate
links follow from the same reciprocity and the prescribed cyclic
permutations of the three labels.
\end{proof}

\begin{remark}
For a YBR fusion bialgebra, the forest criterion of
Lu and Liu \cite{LuThesis,LuLiuER} states that exchange relations
are equivalent to all the fusion graphs being forests. In the link
notation this is
$H^1(\operatorname{Lk}_{K(B)}(i_X);\C)=0$ for every $i\in I$.
Frobenius reciprocity and nonzero quantum dimensions identify the corresponding
links in the three colour classes, up to duality relabelling and exchange
of the two parts. Thus the choice of colour does not change the forest
condition. In the YBR class, the ER criterion concerns first cohomology of the links,
whereas the TR criterion concerns second cohomology of the whole complex.

\end{remark}

Because $K$ is tripartite, its oriented edges split into three types:
\[
K_{XY}^{(1)}
=
\{[x_X,y_Y]\in K^{(1)}\},
\qquad
K_{XZ}^{(1)}
=
\{[x_X,z_Z]\in K^{(1)}\},
\qquad
K_{YZ}^{(1)}
=
\{[y_Y,z_Z]\in K^{(1)}\}.
\]
A $1$-cochain $\alpha\in C^1(K;\C) = \operatorname{Hom}_{\C}(C_1(K;\C),\C)$ is therefore determined by its three edge components.  We write
\[
\alpha_{XY}(x,y):=\alpha([x_X,y_Y]),\quad
\alpha_{XZ}(x,z):=\alpha([x_X,z_Z]),\quad
\alpha_{YZ}(y,z):=\alpha([y_Y,z_Z]).
\]
These expressions are defined whenever the indicated edge belongs to
$K^{(1)}$.

Orient every face by
$[x_X,y_Y,z_Z]$.
Then
\[
\partial[x_X,y_Y,z_Z]
=
[y_Y,z_Z]-[x_X,z_Z]+[x_X,y_Y].
\]
Consequently, the simplicial coboundary
$d_K:C^1(K;\C)\longrightarrow C^2(K;\C)$
satisfies
\begin{equation}
\label{eq:coboundary-components}
(d_K\alpha)(x,y,z)
=
\alpha_{YZ}(y,z)-\alpha_{XZ}(x,z)+\alpha_{XY}(x,y)
\end{equation}
for every face $(x,y,z)\in K^{(2)}$.

For identities in an actual 3-box space, fix a YBR realization
and use Proposition~\ref{prop:triangle-basis}. For an abstract support,
the cochain equations below instead construct candidate coefficients.
Their solvability alone does not assert a planar-algebra realization.

Recall from \eqref{cons:E1} that each $E_\alpha$ is a sum of
the active triangle basis vectors $T(x,y,z)$, indexed by
$(x,y,z)\in K^{(2)}$.

Fix a face $\tau=(i,j,k)\in K^{(2)}$.  A proposed expansion of
the corresponding triangle has the form
\begin{equation}
\label{eq:fixed-triangle-expansion}
T(i,j,k)
=
\sum_{x,y}c_{1,xy}^{ijk}E_1(x,y)
+
\sum_{x,z}c_{2,xz}^{ijk}E_2(x,z)
+
\sum_{y,z}c_{3,yz}^{ijk}E_3(y,z).
\end{equation}
Substituting
\eqref{cons:E1} into
\eqref{eq:fixed-triangle-expansion}, the right-hand side becomes
\[
\sum_{x,y,z}
\left(
c_{1,xy}^{ijk}
+
c_{2,xz}^{ijk}
+
c_{3,yz}^{ijk}
\right)T(x,y,z).
\]
On the other hand,
\[
T(i,j,k)
=
\sum_{x,y,z}
\delta_{x,i}\delta_{y,j}\delta_{z,k}\,T(x,y,z).
\]
Since the nonzero triangle vectors form a basis, the proposed expansion is
valid if and only if the following \emph{prior equations} hold:
\begin{equation}
\label{eq:weighted-prior-equations}
\left(
c_{1,xy}^{ijk}
+
c_{2,xz}^{ijk}
+
c_{3,yz}^{ijk}
-
\delta_{x,i}\delta_{y,j}\delta_{z,k}
\right)N_{zy}^{x}
=
0
\end{equation}
for every $x,y,z\in I$.  Equivalently, on the active faces
$(x,y,z)\in K^{(2)}$, this is the system
\begin{equation}
\label{eq:prior-equations}
c_{1,xy}^{ijk}
+
c_{2,xz}^{ijk}
+
c_{3,yz}^{ijk}
=
\delta_{x,i}\delta_{y,j}\delta_{z,k}.
\end{equation}

For the fixed face $\tau=(i,j,k)$, define a $1$-cochain
$\alpha_\tau\in C^1(K;\C)$ by its three edge components:
\begin{align}
(\alpha_\tau)_{XY}(x,y)
&=
c_{1,xy}^{ijk},
\label{eq:alpha-c1}\\
(\alpha_\tau)_{XZ}(x,z)
&=
-c_{2,xz}^{ijk},
\label{eq:alpha-c2}\\
(\alpha_\tau)_{YZ}(y,z)
&=
c_{3,yz}^{ijk}.
\label{eq:alpha-c3}
\end{align}
The minus sign in the $XZ$ component compensates for the alternating sign
in the oriented boundary.  By
\eqref{eq:coboundary-components},
\[
(d_K\alpha_\tau)(x,y,z)
=
c_{1,xy}^{ijk}
+
c_{2,xz}^{ijk}
+
c_{3,yz}^{ijk}.
\]
For each $\tau\in K^{(2)}$, let $e_\tau\in C^2(K;\C)$ be the
$2$-cochain that takes the value $1$ on $\tau$ and $0$ on every
other triangle, with respect to the chosen orientations.
We then obtain the following reformulation.

\begin{proposition}
\label{prop:prior-as-coboundary}
Let $\tau=(i,j,k)\in K^{(2)}$.  The following are equivalent:
\begin{enumerate}[label=\textup{(\roman*)},leftmargin=2.4em]
  \item the triangle $T(i,j,k)$ admits the expansion
  \eqref{eq:fixed-triangle-expansion}.
  \item its coefficients satisfy the weighted prior equations
  \eqref{eq:weighted-prior-equations}.
  \item the associated $1$-cochain $\alpha_\tau$ satisfies
  $d_K\alpha_\tau=e_\tau$.
\end{enumerate}
The triangle coefficients are recovered from $\alpha_\tau$ by
\eqref{eq:alpha-c1}--\eqref{eq:alpha-c3}.
\end{proposition}

\begin{proof}
Assume first that the vector identity in \textup{(i)} holds.  Substituting
\eqref{cons:E1} into its right-hand
side gives
\[
\sum_{x,y,z}
\left(
 c_{1,xy}^{ijk}+c_{2,xz}^{ijk}+c_{3,yz}^{ijk}
\right)T(x,y,z).
\]
The left-hand side can be written as
\[
\sum_{x,y,z}
\delta_{x,i}\delta_{y,j}\delta_{z,k}T(x,y,z).
\]
The nonzero vectors $T(x,y,z)$ form an orthogonal basis, so comparison of
the active coefficients gives \eqref{eq:prior-equations}, equivalently the
weighted equations \eqref{eq:weighted-prior-equations} on all triples.
Thus \textup{(i)} implies \textup{(ii)}.

Conversely, suppose \eqref{eq:weighted-prior-equations} holds and let $R$
denote the proposed right-hand side of
\eqref{eq:fixed-triangle-expansion}.  The same identity resolutions give
\[
R-T(i,j,k)
=
\sum_{x,y,z}
\left(
 c_{1,xy}^{ijk}+c_{2,xz}^{ijk}+c_{3,yz}^{ijk}
 -\delta_{x,i}\delta_{y,j}\delta_{z,k}
\right)T(x,y,z).
\]
If $N_{zy}^{x}\neq0$, the coefficient in parentheses vanishes by the
weighted prior equation.  If $N_{zy}^{x}=0$, then $T(x,y,z)=0$.  Hence every
summand vanishes and $R=T(i,j,k)$, proving \textup{(ii)} implies
\textup{(i)}.

Finally, the weighted equations impose a condition only on the active
triples, which are exactly the $2$-faces of $K$.  By
\eqref{eq:alpha-c1}--\eqref{eq:alpha-c3} and
\eqref{eq:coboundary-components}, their left-hand side on a face
$(x,y,z)$ is
\[
(d_K\alpha_\tau)(x,y,z)
=
 c_{1,xy}^{ijk}+c_{2,xz}^{ijk}+c_{3,yz}^{ijk}.
\]
Thus the active prior system says exactly that this value is $1$ on
$\tau=(i,j,k)$ and $0$ on every other face. In other words,
$d_K\alpha_\tau=e_\tau$.  This proves \textup{(ii)} and \textup{(iii)} are
equivalent, and the recovery formulas are built into the definition of
$\alpha_\tau$.
\end{proof}

\begin{theorem}
\label{thm:cohomological-criterion}\label{lem:rotation-transfer}
Let $B$ be a YBR fusion bialgebra and let $K(B)$ be its fusion complex.  Then
\[
B\text{ is a TR fusion bialgebra}
\quad\Longleftrightarrow\quad
H^2(K(B);\C)=0.
\]
\end{theorem}

\begin{proof}
Assume first that $B$ is a TR fusion bialgebra, and choose a TR planar
algebra $P$ realizing $B\cong P_{2,+}$.  For each active face
$\tau\in K(B)^{(2)}$, the plus-shaded triangle $T(\tau)$ has a normalized
TR expansion.  By \Cref{prop:prior-as-coboundary}, its coefficient table
determines a cochain $\alpha_\tau$ satisfying
$d_K\alpha_\tau=e_\tau$.
Thus every standard basis vector of $C^2(K;\C)$ lies in the image of
$d_K$, so
$d_K:C^1(K;\C)\longrightarrow C^2(K;\C)$
is surjective.  Since $K$ has no cells above dimension two,
$H^2(K;\C)=C^2(K;\C)/\im d_K=0$.

Conversely, assume $H^2(K(B);\C)=0$, and choose a YBR planar algebra $P$
realizing $B\cong P_{2,+}$.  Again because $K$ is two-dimensional, the
vanishing of $H^2$ is equivalent to the surjectivity of $d_K$.  Hence, for
each face $\tau=(i,j,k)$, choose
$\alpha_\tau\in C^1(K;\C) \;\text{with}\; d_K\alpha_\tau=e_\tau$.
Define the coefficients of the plus-shaded triangle expansion by
\[
c_{1,xy}^{ijk}=(\alpha_\tau)_{XY}(x,y),\qquad
c_{2,xz}^{ijk}=-(\alpha_\tau)_{XZ}(x,z),\qquad
c_{3,yz}^{ijk}=(\alpha_\tau)_{YZ}(y,z).
\]
Then \Cref{prop:prior-as-coboundary} shows that these coefficients satisfy
the prior equations and therefore give a valid expansion of
$T(i,j,k)$ in the three $E$-families.  Thus every active $T$-triangle
lies in $\Espan$. Since the YBR relation expresses every $W$ as a
linear combination of $T$-triangles, all $W$-triangles also lie in $\Espan$.  The chosen YBR realization is therefore a TR planar
algebra, and $B$ is a TR fusion bialgebra.
\end{proof}

\begin{corollary}
\label{cor:base-field}
For a finite fusion complex $K$,
\[
H^2(K;\C)=0
\quad\Longleftrightarrow\quad
H^2(K;\Q)=0.
\]
\end{corollary}

\begin{proof}
Let $\Bbbk$ be either $\Q$ or $\C$.  The universal coefficient theorem for
cohomology gives a natural short exact sequence
\[
0
\longrightarrow
\operatorname{Ext}_{\Z}^{1}
\bigl(H_1(K;\Z),\Bbbk\bigr)
\longrightarrow
H^2(K;\Bbbk)
\longrightarrow
\operatorname{Hom}_{\Z}
\bigl(H_2(K;\Z),\Bbbk\bigr)
\longrightarrow
0.
\]
Both $\Q$ and $\C$ are divisible, hence injective, $\Z$-modules, so the
$\operatorname{Ext}$ term vanishes.  Moreover, since $K$ is a finite
two-dimensional simplicial complex,
$H_2(K;\Z)\subseteq C_2(K;\Z)$
is a finitely generated free abelian group, say of rank $r$.  Therefore
$H^2(K;\Q)\cong\Q^r, \; H^2(K;\C)\cong\C^r$.
The two cohomology groups consequently vanish simultaneously.
\end{proof}

\begin{remark}
Since the coboundary matrix has integer entries, its rank over
$\Q$ equals its rank over $\C$.
Thus the condition $H^2(K;\C)=0$ can be checked by exact
row reduction over $\Q$.
\end{remark}

\begin{definition}
Assume $H^2(K;\Q)=0$.  The \emph{Hodge representative} of the
triangle coefficient system for a face $\tau$ is the unique solution
$\alpha^{\mathrm H}_\tau$ of
$d_K\alpha^{\mathrm H}_\tau=e_\tau$
that is orthogonal to $\ker d_K$ with respect to the standard inner product
on the edge-cochain space.
\end{definition}

\begin{theorem}
\label{thm:hodge-gauge}
For each face $\tau\in K^{(2)}$, the Hodge representative exists, is unique,
and is given by
\[
\boxed{
\alpha^{\mathrm H}_\tau
=
d_K^*(d_Kd_K^*)^{-1}e_\tau.
}
\]
Moreover, $\alpha^{\mathrm H}_\tau$ has rational coordinates and is the
unique minimum-norm exact solution:
\[
\alpha^{\mathrm H}_\tau
=
\underset{d_K\alpha=e_\tau}{\operatorname{argmin}}\,\|\alpha\|_2.
\]
\end{theorem}

\begin{proof}
The surjectivity of $d_K$ implies that $d_Kd_K^*$ is invertible on
$C^2(K;\Q)$.  The displayed vector is a solution and belongs to
$\im d_K^*=(\ker d_K)^\perp$.  If two solutions are orthogonal to
$\ker d_K$, their difference belongs to both $\ker d_K$ and its orthogonal
complement, and is therefore zero.  Rationality follows because $d_K$ is an integer matrix.  Finally, every
solution is uniquely of the form
$\alpha^{\mathrm H}_\tau+\beta$ with $\beta\in\ker d_K$, and orthogonality
gives
$
\|\alpha^{\mathrm H}_\tau+\beta\|_2^2
=
\|\alpha^{\mathrm H}_\tau\|_2^2+\|\beta\|_2^2.
$
Hence the Hodge representative is the unique solution of smallest norm.
\end{proof}

\begin{remark}
The triangle expansion coefficients are generally not unique.
For each active triangle $\tau$, every solution of
$d_K\alpha_\tau=e_\tau$ is uniquely of the form
\[
\alpha_\tau=\alpha_\tau^{\mathrm H}+\beta,
\qquad \beta\in Z^1(K)=\ker d_K.
\]
Thus adding a $1$-cocycle changes the coefficients without
changing the represented $3$-box vector.
The full coefficient ambiguity is parametrized by $Z^1(K)$.

The Hodge gauge selects the unique solution orthogonal to
$\ker d_K$ with respect to the standard cellular inner product.
This solution is rational and compatible with relabellings
of the fusion complex, giving a unique exact coefficient table
for each labelled fusion complex in the classification algorithm. This choice need not minimize the number of nonzero coefficients.
\end{remark}

Choose orderings of the active faces and oriented edges of $K$.
Define the matrix $D$ with rows indexed by
$\sigma\in K^{(2)}$ and columns indexed by $e\in K^{(1)}$ by
\[
D_{(x,y,z),e}=
\begin{cases}
 1,  & e=[x_X,y_Y],\\
 -1, & e=[x_X,z_Z],\\
 1,  & e=[y_Y,z_Z],\\
 0,  & \text{otherwise}.
\end{cases}
\]
Thus $D$ has $|K^{(2)}|$ rows and $3n^2$ columns.
Representing a $1$-cochain by the column vector of its edge
values, we have
$(D\alpha)_{(x,y,z)}=(d_K\alpha)(x,y,z)$.
For each fixed target face $\tau$, the system
$D\alpha_\tau=e_\tau$ has one equation for every active face:
the right-hand side is $1$ in the row indexed by $\tau$
and $0$ in all other rows.
The same matrix $D$ is used for every target face.

\begin{theorem}\label{thm:canonical-readout}
Let $K$ be a fusion complex with $H^2(K;\mathbb Q)=0$, and
let $D:C^1(K;\mathbb Q)\to C^2(K;\mathbb Q)$ be its oriented coboundary
matrix. For each input $(i,j,k)$, let $b_{ijk}=e_{ijk}$ if $(i,j,k)$
is a face and $b_{ijk}=0$ otherwise. Define
\[
 \alpha_{ijk}=D^{\mathsf T}(DD^{\mathsf T})^{-1}b_{ijk}.
\]
\[
 c_{1,xy}^{ijk}=(\alpha_{ijk})_{XY}(x,y),\quad
 c_{2,xz}^{ijk}=-(\alpha_{ijk})_{XZ}(x,z),\quad
 c_{3,yz}^{ijk}=(\alpha_{ijk})_{YZ}(y,z).
\]
This gives a rational coefficient table satisfying prior and having zero
columns for every inactive input. It is the Hodge representative on
each active input. If $K=K(B)$ for a TR realization, this table gives
its $T$-triangle expansions.
\end{theorem}
\begin{proof}
Since $D$ is surjective, $DD^{\mathsf T}$ is invertible over $\mathbb Q$.
Thus $D\alpha_{ijk}=b_{ijk}$, which is prior with the stated edge signs.
For an inactive input the right-hand side is zero and the displayed
formula gives $\alpha_{ijk}=0$. In a TR realization, expand the
$E_\alpha$ using the identity resolution. Prior gives precisely the
input active triangle, and inactive triangles vanish.
\end{proof}

\subsection{Reconstruction}
\label{sec:twosided}
\label{sec:reconstruction}
% ============================================================
We now study the reconstruction of bi-commutative TR planar
algebras from fusion complex pairs. The preceding cohomological
construction provides triangle coefficients on both sides.
We determine when these coefficients and compatible fusion data
can be realized by a TR planar algebra.

\begin{definition}
A fusion bialgebra is \emph{bi-commutative} if both its ordinary product and
its convolution product are commutative.
\end{definition}

If every label is self-dual, then
$N_{ij}^{k}=N_{\bar j\bar i}^{\bar k}=N_{ji}^{k}$.
Thus every self-dual fusion bialgebra is bi-commutative.
Bi-commutativity allows us to choose minimal-projection bases on both
shadings and to describe both sides by fusion complexes.

Let $B_+$ and $B_-$ have distinguished bases of minimal idempotents
$\{p_i\}_{i=0}^{n-1}$ and $\{p'_a\}_{a=0}^{n-1}$,
with convolution units $p_0$ and $p'_0$, respectively.
The labels on the two sides are initially independent.

\begin{definition}
\label{def:fusion-bialgebra-pair}
An ordered pair $(B_+,B_-)$ of bi-commutative fusion bialgebras is called a
\emph{fusion bialgebra pair} if
\begin{enumerate}[label=\textup{(\roman*)},leftmargin=2.4em]
  \item $\dim_{\C}B_+=\dim_{\C}B_-=n$.
  \item after relabelling, their duality involutions have the same cycle
  type and are represented by the same convention $i\mapsto\dual{i}$.
\end{enumerate}
\end{definition}

\begin{definition}
\label{def:bare-fusion-complex-pair}
A \emph{fusion complex pair} is an ordered pair $(K_+,K_-)$
of fusion complexes of the same rank and 
duality type after relabelling.
A fusion-bialgebra realization is a pair $(B_+,B_-)$ together with
support identifications $K(B_\pm)\cong K_\pm$. 
\end{definition}
\begin{corollary}
\label{cor:opposite-vanishing}
Let $P$ be a bi-commutative YBR planar algebra with associated fusion
complexes $K_+$ and $K_-$.  Then
\[
H^2(K_+;\C)=0
\quad\Longleftrightarrow\quad
H^2(K_-;\C)=0.
\]
\end{corollary}

\begin{proof}
Either condition is equivalent, by \Cref{thm:cohomological-criterion}, to
$P$ being of TR type.
\end{proof}

\begin{lemma}
\label{lem:equal-faces}
If $B_+$ and $B_-$ arise as $P_{2,+}$ and $P_{2,-}$ of one YBR planar
algebra, then
$|K_+^{(2)}|=|K_-^{(2)}|$.
\end{lemma}

\begin{proof}
The nonzero triangles form orthogonal bases of $P_{3,+}$ and $P_{3,-}$,
respectively.  Hence
$|K_+^{(2)}|=\dim P_{3,+}, \; |K_-^{(2)}|=\dim P_{3,-}$.
The quarter-rotation $\F$ is a linear isomorphism between the two box
spaces, so the dimensions are equal.
\end{proof}

The TR relations were previously formulated as expansions of
the two triangle families $T$ and $W$, both labelled by
minimal idempotents of $B_+$.
For a bi-commutative fusion bialgebra pair $(B_+,B_-)$,
we instead use one triangle relation on each shading,
with labels from the corresponding fusion bialgebra.

On the plus side, we retain $T(i,j,k)$, labelled by
$p_i,p_j,p_k\in B_+$.
On the minus side, let $T'(a,b,c)$ denote the corresponding
triangle labelled by $p'_a,p'_b,p'_c\in B_-$.
Here $a,b,c$ index the minimal idempotents of $B_-$.
Write $E'_\alpha(d,e)$ for the corresponding diagrams
without internal faces, labelled by $p'_d,p'_e$.

These diagrams are obtained from the corresponding
$E_\alpha$ diagrams by reversing the shading and
replacing the labels, with the box orientations unchanged.

\begin{figure}[htbp]
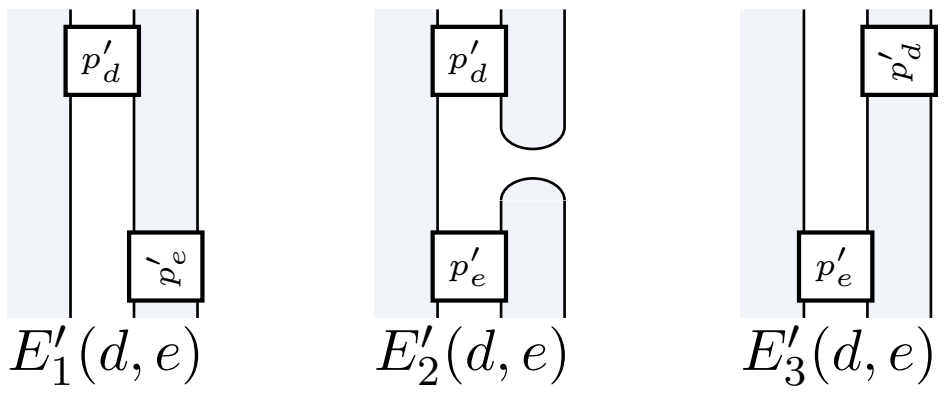

\centering
\[
\begin{array}{ccc}
\opEprime{1}
&
\qquad\opEprime{2}\qquad
&
\opEprime{3}
\\[8pt]
E'_1(d,e)
&
E'_2(d,e)
&
E'_3(d,e)
\end{array}
\]
\caption{The minus-shaded diagrams $E'_\alpha(d,e)$.}
\label{fig:Eprime-diagrams}
\end{figure}
The minus-side triangle relation is
\begin{equation}\label{eq:minus-tr-expansion}
 T'(a,b,c)
 =\sum_{d,e\in I_-}\Bigl(
 c'{}_{1,de}^{abc}E'_1(d,e)
 +c'{}_{2,de}^{abc}E'_2(d,e)
 +c'{}_{3,de}^{abc}E'_3(d,e)\Bigr).
\end{equation}
Graphically, this reads
\begin{equation}\label{eq:minus-tr-expansion-picture}
 \opTprime
 =\sum_{d,e\in I_-}\left(
 c'{}_{1,de}^{abc}\,\opEprime{1}
 +c'{}_{2,de}^{abc}\,\opEprime{2}
 +c'{}_{3,de}^{abc}\,\opEprime{3}
 \right).
\end{equation}
The fusion complexes $K_+$ and $K_-$ provide the
coefficient tables $c$ and $c'$ for these two relations.

The earlier formulation expresses the second triangle
relation in $P_{3,+}$ using $W$ and plus-side labels.
The present formulation expresses it directly in
$P_{3,-}$ using $T'$ and minus-side labels.
Under Fourier matching, the rotation and label changes
specified by our diagram conventions convert the $T'$
family into the $W$ family and the $E'_\alpha$ families
into the $E_\alpha$ families.
These changes are invertible, so $T'$ admits an expansion
in the $E'_\alpha$ if and only if $W$ admits an expansion
in the $E_\alpha$.
Lemma~\ref{prop:cprime-to-ctilde} gives the explicit
formulas for the table $\widetilde c$ obtained by
transporting $c'$.

With the triangle coefficient tables $c,c'$ fixed, the
trace closures of the two triangle relations express the
fusion coefficients in terms of the quantum dimensions.
Write $I_+$ and $I_-$ for the two label sets, and put
\[
Q_{ij}^{k}:=d_kN_{ij}^{k},
\qquad
Q'{}_{ab}^{c}:=d'_cN'{}_{ab}^{c}.
\]
We first carry out the calculation on the plus side.

Put $Q_{ij}^{k}=d_kN_{ij}^{k}$. 
We evaluate the same closure of $T(i,j,k)$ in two ways.
First, applying the normalized convolution relation and then taking
the unnormalized trace gives
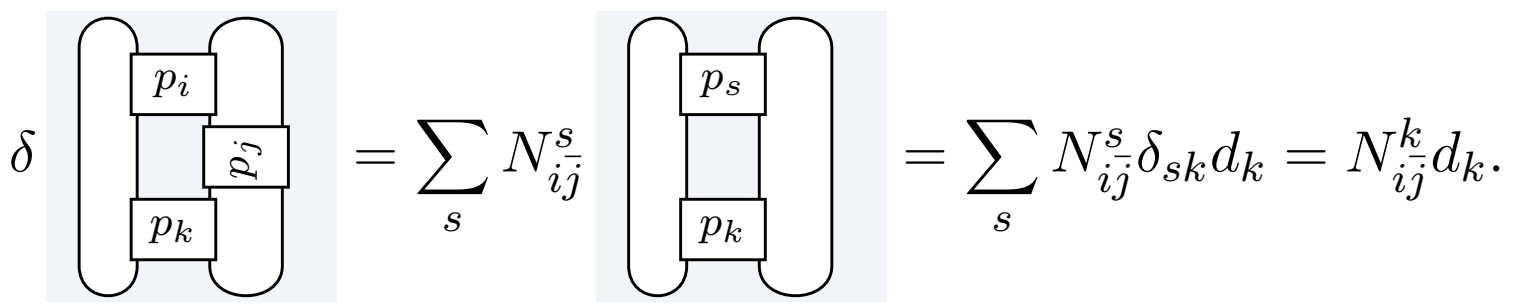
\begin{equation}\label{eq:Q-closure-direct}
\begin{aligned}
  \delta\,\qClosedT
  =
  \sum_s N_{i\bar j}^{s}\,\qProductTrace
  =
  \sum_s N_{i\bar j}^{s}\delta_{sk}d_k
  =
  N_{i\bar j}^{k}d_k.
\end{aligned}
\end{equation}

Alternatively, we first apply the triangle relation and then close
each summand using the same boundary pairing:
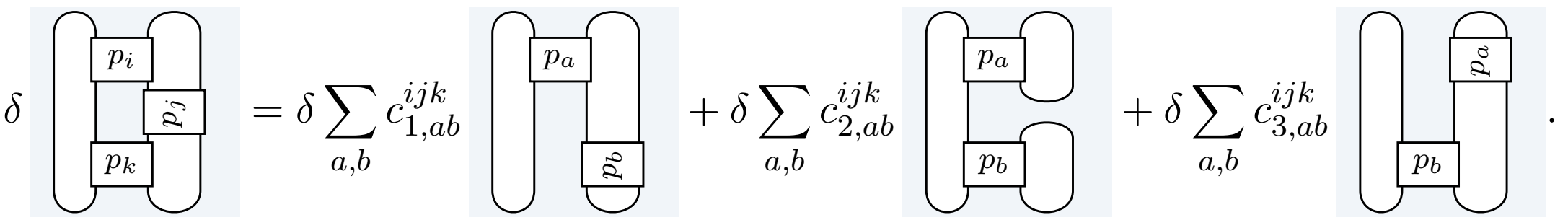
\begin{equation}\label{eq:Q-closure-TR}
\begin{aligned}
  \delta\,\qClosedT
  ={}
  \delta\sum_{a,b}c_{1,ab}^{ijk}\,\qClosedE{1}
  +
  \delta\sum_{a,b}c_{2,ab}^{ijk}\,\qClosedE{2}
  +
  \delta\sum_{a,b}c_{3,ab}^{ijk}\,\qClosedE{3}.
\end{aligned}
\end{equation}
The elementary closure relations give
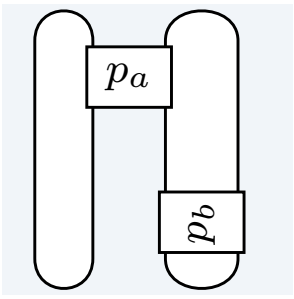
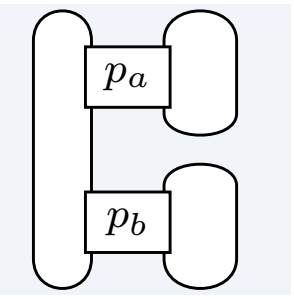
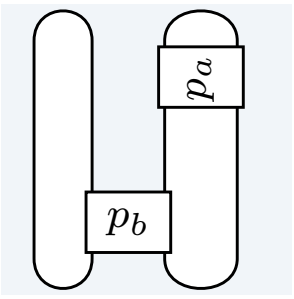
\begin{equation}\label{eq:Q-closure-E}
  \delta\,\qClosedE{1}
  =
  \delta\,\qClosedE{2}
  =
  \delta\,\qClosedE{3}
  =
  d_a d_b.
\end{equation}
Comparing the two evaluations, we obtain
\begin{equation}\label{eq:Q-closure-comparison}
  N_{i\bar j}^{k}d_k
  =
  \sum_{a,b}
  \left(
    c_{1,ab}^{ijk}
    +c_{2,ab}^{ijk}
    +c_{3,ab}^{ijk}
  \right)d_a d_b.
\end{equation}
Replacing $j$ by $\bar j$ gives
\begin{equation}\label{eq:Q-closure-quadratic}
 Q_{ij}^{k}(d)
 =
 \sum_{a,b\in I_+}
 \left(
 c_{1,ab}^{i\bar j k}
 +c_{2,ab}^{i\bar j k}
 +c_{3,ab}^{i\bar j k}
 \right)d_a d_b.
\end{equation}
The same closure calculation for $T'$ on the minus side gives
\begin{equation}\label{eq:Qprime-closure-quadratic}
 Q'{}_{ab}^{c}(d')
 =
 \sum_{u,v\in I_-}
 \left(
 c'{}_{1,uv}^{a\bar b c}
 +c'{}_{2,uv}^{a\bar b c}
 +c'{}_{3,uv}^{a\bar b c}
 \right)d'_u d'_v.
\end{equation}
Thus, for nonzero quantum dimensions,
\[
N_{ij}^{k}=\frac{Q_{ij}^{k}(d)}{d_k},
\qquad
N'{}_{ab}^{c}=\frac{Q'{}_{ab}^{c}(d')}{d'_c}.
\]

For fixed $K_\pm$ and their canonical tables $c,c'$, the
fusion coefficients $N,N'$ are determined by $d,d'$.
Substitution into the fusion structure equations gives
polynomial equations over $\Q$ after clearing denominators.
We impose $d_0=d'_0=1$, nonzero dimensions and total
dimensions, and require the supports of $Q,Q'$ to be
exactly those specified by $K_+,K_-$.

We next match the two candidate fusion bialgebras
$B_+$ and $B_-$.
In a planar algebra, the Fourier maps exchange ordinary
multiplication and convolution, and their compositions
implement duality.
The following definition formulates these requirements
for the abstract pair.

\begin{definition}
\label{def:ssprime-matched}
Let $(B_+,B_-)$ be a fusion bialgebra pair.  An
\emph{$(S,S')$-Fourier matching} consists of two linear isomorphisms
\[
\F_+:B_+\longrightarrow B_-,
\qquad
\F_-:B_-\longrightarrow B_+,
\]
whose matrices in the distinguished bases are
\begin{equation}
\label{eq:S-definitions}
\F_+(p_i)=\sum_{a\in I_-}S_{ai}p'_a
\quad (i\in I_+),
\qquad
\F_-(p'_a)=\sum_{i\in I_+}S'_{ia}p_i
\quad (a\in I_-),
\end{equation}
and which satisfy the following conditions.
\begin{enumerate}[label=\textup{(\roman*)},leftmargin=2.4em]
  \item A common scalar $\delta\in\C^\times$ satisfies
\[
\delta^2=\sum_{i\in I_+}d_i
        =\sum_{a\in I_-}d'_a.
\]

  \item The Fourier maps exchange the ordinary and convolution products
  with the planar-algebra normalization:
  \begin{align*}
  \F_+(x\cdot_+y)
  &=\frac1\delta\F_+(x)*_-\F_+(y),
  &
  \F_+(x*_+y)
  &=\delta\F_+(x)\cdot_-\F_+(y)
  &&(x,y\in B_+),\\
  \F_-(u\cdot_-v)
  &=\frac1\delta\F_-(u)*_+\F_-(v),
  &
  \F_-(u*_-v)
  &=\delta\F_-(u)\cdot_+\F_-(v)
  &&(u,v\in B_-).
  \end{align*}

  \item Two counterclockwise quarter-turns implement duality:
\[
\F_-\circ\F_+(p_i)=p_{\dual{i}},
\qquad
\F_+\circ\F_-(p'_a)=p'_{\dual{a}}.
\]
Let $C_{\mathrm{dual},+}$ and $C_{\mathrm{dual},-}$
be the permutation matrices defined by
\[
(C_{\mathrm{dual},+})_{ji}=\delta_{j,\dual{i}},
\qquad
(C_{\mathrm{dual},-})_{ba}=\delta_{b,\dual{a}}.
\]
In the column-matrix convention, the duality identities
are equivalent to
\begin{equation}
\label{eq:Fourier-square}
S'S=C_{\mathrm{dual},+},
\qquad
SS'=C_{\mathrm{dual},-}.
\end{equation}

  \item The Fourier maps satisfy the unit normalizations:
  \begin{equation}
  \label{eq:Jones-Fourier-normalization}
  \F_+(p_0)=\frac1\delta\,1_-,
  \qquad
  \F_-(p'_0)=\frac1\delta\,1_+.
  \end{equation}
\end{enumerate}
The pair equipped with these matrices is called an
\emph{$(S,S')$-matched fusion bialgebra pair}.
A fusion bialgebra pair is \emph{Fourier-matchable}
if such a matching exists.
\end{definition}

\begin{definition}
\label{def:ssprime-matched-complex-pair}
A fusion complex pair is
\emph{$(S,S')$-matched through a chosen realization}
if it is equipped with a fusion-bialgebra realization
$(B_+,B_-;S,S')$ satisfying
Definition~\ref{def:ssprime-matched}.
It is \emph{Fourier-matchable} if such a realization exists.
\end{definition}

\begin{remark}
\label{rmk:duality-type-necessity}
The compatibility of the duality involutions in
Definition~\ref{def:fusion-bialgebra-pair} is forced by Fourier matching rather than
being merely a convenient labelling convention.  Indeed, for any
$(S,S')$-Fourier matching, since $S$ and $S'$ are invertible, one has
$C_{\mathrm{dual},-} =S C_{\mathrm{dual},+}S^{-1}$.
Thus the two duality permutation matrices are similar.  In particular, the
duality involutions on the two sides have the same cycle type.  After a
relabeling preserving the distinguished label $0$, they may therefore be
written using one common convention $i\mapsto\dual{i}$.
This gives the basic motivation for requiring a fusion bialgebra pair, and
hence its associated fusion complex pair, to carry compatible duality data
before any Fourier matching is attempted.
\end{remark}

\begin{remark}
The common labelling is only a bookkeeping convention.  At this stage no
individual minimal projection $p_i\in B_+$ is identified with an
individual minimal projection $p'_a\in B_-$.  The actual relation between
the two distinguished bases is determined only after an $(S,S')$-Fourier
matching has been chosen.
The label $0$ denotes the convolution unit on each side.
It does not mean that $\F_+(p_0)=p'_0$.
Under Fourier matching, the convolution unit on each side
is sent to $\delta^{-1}$ times the opposite ordinary unit.
\end{remark}

\begin{proposition}\label{prop:S-character-tables}
For an $(S,S')$-matched pair, there are fusion characters
$\chi_a^+:(B_+,*)\to\C$ and $\chi_i^-:(B_-,*)\to\C$ such that
\begin{equation}\label{eq:S-from-fusion-characters}
 S_{ai}=\delta^{-1}\chi_a^+(p_i),\qquad
 S'_{ia}=\delta^{-1}\chi_i^-(p'_a).
\end{equation}
Thus $\delta S$ and $\delta S'$ are the fusion-character tables, with
rows matched to the opposite minimal idempotents. The row indexed by $0$ in each fusion-character table
is the dimension character. The following normalizations hold:
\label{eq:S-first-columns}\label{eq:S-row-sums}
\begin{equation}\label{eq:unit-Fourier-normalization}
 \F_+(1_+)=\delta p'_0,\quad \F_-(1_-)=\delta p_0,\quad
 S_{a0}=S'_{i0}=\delta^{-1},\quad
 \sum_i S_{ai}=\delta\delta_{a0},\quad
 \sum_a S'_{ia}=\delta\delta_{i0}.
\end{equation}
\end{proposition}
\begin{proof}
The product-exchange identities make $\delta\F_+$ and $\delta\F_-$
algebra isomorphisms from convolution to ordinary multiplication.
Composing with the ordinary coordinate characters gives the stated
fusion characters. The Jones normalization and Fourier-square identities
give the displayed unit, column and row formulas.
\end{proof}

For computation, let $\mathrm{Char}_\pm$ be the two character tables.
Enumerate row permutations $P_+,P_-$ matching their characters to the
opposite minimal idempotents, and form
\begin{equation}\label{eq:raw-character-matching}
 S=\delta^{-1}P_+\mathrm{Char}_+,\qquad
 S'=\delta^{-1}P_-\mathrm{Char}_-.
\end{equation}
There is no further row rescaling. The resulting candidates
must satisfy all the conditions of
Definition~\ref{def:ssprime-matched}.

\begin{remark}
If $B_+$ and $B_-$ are isomorphic as fusion bialgebras
via $p_i\mapsto p'_{\sigma(i)}$, then
$d_i=d'_{\sigma(i)}$ for every $i\in I_+$.
Indeed, these dimensions are the coefficients of the
convolution unit in $p_i*p_{\dual i}$ and
$p'_{\sigma(i)}*p'_{\dual{\sigma(i)}}$, respectively.
An isomorphism of the fusion complexes alone does not
imply this equality.
\end{remark}

\begin{lemma}\label{lem:inverse-Fourier-bases}\label{eq:inverse-minus-basis}\label{eq:inverse-plus-basis}
For $a\in I_-$ and $i\in I_+$,
\[
 p'_a=\sum_iS'_{i,\bar a}\F_+(p_i),\qquad
 p_i=\sum_aS_{a,\bar i}\F_-(p'_a).
\]
\end{lemma}
\begin{proof}
These identities follow from $SS'=C_{\mathrm{dual},-}$ and
$S'S=C_{\mathrm{dual},+}$.
\end{proof}

We now give the explicit coefficient formulas for transporting
the $T'$ expansion to the $W$ expansion under Fourier matching.
We use the half-turn $\rho_3$ of
Figure~\ref{fig:annular-half-turn} and the label conventions of
Figure~\ref{fig:triangle-convention}.

\begin{lemma}
\label{prop:cprime-to-ctilde}
For a fixed $(S,S')$-Fourier matching, the $T'$ expansion
with coefficients $c'$ induces a $W$ expansion whose
coefficients are given by
\begin{align}
\widetilde c_{1,xy}^{ijk}
&=
\sum_{a,b,c,d,e\in I_-}
S_{a,\dual{k}}S_{b,j}S_{c,\dual{i}}
S'_{\dual{y},\dual{d}}S'_{\dual{x},e}
 c'{}_{1,de}^{abc},
\label{eq:ctilde-1}\\
\widetilde c_{2,xz}^{ijk}
&=
\sum_{a,b,c,d,e\in I_-}
S_{a,\dual{k}}S_{b,j}S_{c,\dual{i}}
S'_{z,\dual{d}}S'_{\dual{x},\dual{e}}
 c'{}_{2,de}^{abc},
\label{eq:ctilde-2}\\
\widetilde c_{3,yz}^{ijk}
&=
\sum_{a,b,c,d,e\in I_-}
S_{a,\dual{k}}S_{b,j}S_{c,\dual{i}}
S'_{\dual{z},d}S'_{\dual{y},\dual{e}}
 c'{}_{3,de}^{abc}.
\label{eq:ctilde-3}
\end{align}
The first identity holds for all $i,j,k,x,y\in I_+$, the second for all
$i,j,k,x,z\in I_+$, and the third for all $i,j,k,y,z\in I_+$. 
\end{lemma}

\begin{proof}
Apply the Fourier changes of labels to the $T'$ relation
in \eqref{eq:minus-tr-expansion-picture}, using
Lemma~\ref{lem:inverse-Fourier-bases}, and then apply
the half-turn $\rho_3$.
Before the half-turn, the three input labels contribute
$S_{ai}S_{b,\bar j}S_{ck}$.
The half-turn replaces $(i,j,k)$ by
$(\bar k,\bar j,\bar i)$, so the common input factor is
$
S_{a,\bar k}S_{b,j}S_{c,\bar i}.
$
For the output diagrams, the same changes of labels
and half-turn give
\[
\begin{aligned}
E'_1(d,e)
&\longmapsto
\sum_{x,y\in I_+}
S'_{\bar y,\bar d}S'_{\bar x,e}\,E_1(x,y),\\
E'_2(d,e)
&\longmapsto
\sum_{x,z\in I_+}
S'_{z,\bar d}S'_{\bar x,\bar e}\,E_2(x,z),\\
E'_3(d,e)
&\longmapsto
\sum_{y,z\in I_+}
S'_{\bar z,d}S'_{\bar y,\bar e}\,E_3(y,z).
\end{aligned}
\]
Here the arrows denote the combined change of labels
and rotation. The duals in each factor follow from
the orientations of the two labelled boxes.

For each $\alpha$, multiply the corresponding output
factor by $S_{a,\bar k}S_{b,j}S_{c,\bar i}$ and
$c'{}_{\alpha,de}^{abc}$, and sum over
$a,b,c,d,e\in I_-$.
The resulting coefficients of the $W$ expansion are
exactly \eqref{eq:ctilde-1}--\eqref{eq:ctilde-3}.
\end{proof}

% Local replacement fragment; this is not a standalone document.
% Replace the entire uploaded passage, from the theorem carrying
% thm:final-consistency-system through the final elimination remark.
% Existing definitions, graphical conventions, and preamble are unchanged.
% The finite-dimensionality citation placeholder at the proof's end is
% retained explicitly and still requires the author's verification.

The graphical consistency comparisons can now be expressed
using the two native coefficient tables $c,c'$ and the
Fourier matrices $S,S'$.

\begin{theorem}\label{thm:final-consistency-system}
\label{thm:unique-33e}
Let $P$ be a bi-commutative YBR planar algebra and let $B_\pm=P_{2,\pm}$.
Write $N,d$ and $N',d'$ for their fusion coefficients and quantum dimensions,
$S,S'$ for the Fourier matrices, and $c,c'$ for the coefficients of the
proposed native TR relations.
The TR consistency tests of Theorem~\ref{cons:graphthm},
together with the identity-resolution tests of
Theorem~\ref{thm:partition-function}, are equivalent
to the following polynomial equations:
\begingroup
\small
\allowdisplaybreaks[4]
\setlength{\jot}{5pt}
\setlength{\abovedisplayskip}{8pt}
\setlength{\belowdisplayskip}{8pt}
\setlength{\abovedisplayshortskip}{4pt}
\setlength{\belowdisplayshortskip}{6pt}

\begin{gather}
\makebox[.45\linewidth]{$\displaystyle
d_i=d_{\bar i}$}
\quad
\makebox[.45\linewidth]{$\displaystyle
d'_a=d'_{\bar a}$}
\label{eq:pair-dim-dual}\\
\makebox[.45\linewidth]{$\displaystyle
N_{ij}^{0}=\delta_{i,\bar j}d_i$}
\quad
\makebox[.45\linewidth]{$\displaystyle
N'{}_{ab}^{0}=\delta_{a,\bar b}d'_a$}
\label{eq:pair-vacuum}\\
\makebox[.45\linewidth]{$\displaystyle
N_{ij}^{k}=N_{\bar j\,\bar i}^{\bar k}$}
\quad
\makebox[.45\linewidth]{$\displaystyle
N'{}_{ab}^{c}=N'{}_{\bar b\,\bar a}^{\bar c}$}
\label{eq:pair-duality}\\
\makebox[.45\linewidth]{$\displaystyle
N_{ij}^{k}=N_{ji}^{k}$}
\quad
\makebox[.45\linewidth]{$\displaystyle
N'{}_{ab}^{c}=N'{}_{ba}^{c}$}
\label{eq:pair-commutativity}\\
\makebox[.45\linewidth]{$\displaystyle
\sum_s N_{ij}^{s}d_s=d_id_j$}
\quad
\makebox[.45\linewidth]{$\displaystyle
\sum_t N'{}_{ab}^{t}d'_t=d'_ad'_b$}
\label{eq:pair-dimension}\\
\makebox[.45\linewidth]{$\displaystyle
\sum_s N_{ij}^{s}N_{sk}^{\ell}
=\sum_s N_{is}^{\ell}N_{jk}^{s}$}
\quad
\makebox[.45\linewidth]{$\displaystyle
\sum_t N'{}_{ab}^{t}N'{}_{tc}^{e}
=\sum_t N'{}_{at}^{e}N'{}_{bc}^{t}$}
\label{eq:pair-associativity}
\end{gather}
\begin{equation}
\label{eq:pair-31}
\begin{gathered}
\makebox[.45\linewidth]{$\displaystyle c_{1,0\bar s}^{ijk}+c_{2,0s}^{ijk}+c_{3,\bar s s}^{ijk}=\delta_{i0}\delta_{k\bar j}\delta_{sk}$}%
\quad
\makebox[.45\linewidth]{$\displaystyle c'{}_{1,0\bar t}^{abc}+c'{}_{2,0t}^{abc}+c'{}_{3,\bar t t}^{abc}=\delta_{a0}\delta_{c\bar b}\delta_{tc}$}
\end{gathered}
\end{equation}
\begin{equation}
\label{eq:pair-32v}
\begin{gathered}
\makebox[.45\linewidth]{$\displaystyle \begin{gathered}\bigl[\delta_{x\ell}(c_{1,xy}^{ijk}+c_{2,xz}^{ijk}+c_{3,yz}^{ijk})\\-\delta_{xi}\delta_{yj}\delta_{zk}\delta_{i\ell}\bigr]N_{zy}^{x}=0\end{gathered}$}%
\quad
\makebox[.45\linewidth]{$\displaystyle \begin{gathered}\bigl[\delta_{ue}(c'{}_{1,uv}^{abc}+c'{}_{2,uw}^{abc}+c'{}_{3,vw}^{abc})\\-\delta_{ua}\delta_{vb}\delta_{wc}\delta_{ae}\bigr]N'{}_{wv}^{u}=0\end{gathered}$}
\end{gathered}
\end{equation}
\begin{equation}
\label{eq:pair-32e}
\begin{gathered}
\makebox[.45\linewidth]{$\displaystyle \begin{gathered}\sum_{p,q}c_{1,pq}^{ijk}N_{p\bar q}^{s}+\sum_p(c_{2,ps}^{ijk}+c_{3,ps}^{ijk})d_p\\=\delta_{sk}N_{i\bar j}^{k}\end{gathered}$}%
\quad
\makebox[.45\linewidth]{$\displaystyle \begin{gathered}\sum_{u,v}c'{}_{1,uv}^{abc}N'{}_{u\bar v}^{t}+\sum_u(c'{}_{2,ut}^{abc}+c'{}_{3,ut}^{abc})d'_u\\=\delta_{tc}N'{}_{a\bar b}^{c}\end{gathered}$}
\end{gathered}
\end{equation}
\begin{equation}
\label{eq:final-fourier-plus-product}
\begin{gathered}
\makebox[.45\linewidth]{$\displaystyle \sum_{a,b}S_{ai}S_{bj}N'{}_{ab}^{c}=\delta\delta_{ij}S_{ci}$}%
\quad
\makebox[.45\linewidth]{$\displaystyle \sum_{i,j}S'_{ia}S'_{jb}N_{ij}^{k}=\delta\delta_{ab}S'_{ka}$}
\end{gathered}
\end{equation}
\begin{equation}
\label{eq:final-fourier-plus-convolution}
\begin{gathered}
\makebox[.45\linewidth]{$\displaystyle \sum_kN_{ij}^{k}S_{ak}=\delta S_{ai}S_{aj}$}%
\quad
\makebox[.45\linewidth]{$\displaystyle \sum_cN'{}_{ab}^{c}S'_{ic}=\delta S'_{ia}S'_{ib}$}
\end{gathered}
\end{equation}
\begin{equation}
\label{eq:pair-fourier-square}
\begin{gathered}
\makebox[.45\linewidth]{$\displaystyle S'S=C_{\mathrm{dual},+}$}%
\quad
\makebox[.45\linewidth]{$\displaystyle SS'=C_{\mathrm{dual},-}$}
\end{gathered}
\end{equation}
\begin{equation}
\label{eq:pair-fourier-unit}
\begin{gathered}
\makebox[.45\linewidth]{$\displaystyle \delta S_{a0}=1$}%
\quad
\makebox[.45\linewidth]{$\displaystyle \delta S'_{i0}=1$}
\end{gathered}
\end{equation}\begin{equation}\label{eq:pair-modulus}\delta^2=\sum_{i\in I_+}d_i=\sum_{a\in I_-}d'_a.\end{equation}
\begin{equation}\label{eq:pair-33e}
\begin{gathered}
\Bigg[
 \sum_{r\in I_+}\left(c_{1,ry}^{ij\bar k}N_{\ell r}^{x}
                         +c_{3,yr}^{ij\bar k}N_{\ell r}^{z}\right)\\
+\delta\sum_{\substack{p,q\in I_+\\a,b,c,d,e\in I_-}}
 c_{2,pq}^{ij\bar k}S_{a,q}S_{b,\ell}S_{c,\bar p}\,\Bigl(
 S'_{\bar y,\bar d}S'_{\bar x,e}c'{}_{1,de}^{abc}
 +S'_{z,\bar d}S'_{\bar x,\bar e}c'{}_{2,de}^{abc}
 +S'_{\bar z,d}S'_{\bar y,\bar e}c'{}_{3,de}^{abc}\Bigr)\\
-\delta\delta_{jy}\sum_{a,b,c,d,e\in I_-}
 S_{a,\bar k}S_{b,\ell}S_{c,\bar i}\,\Bigl(
 S'_{\bar y,\bar d}S'_{\bar x,e}c'{}_{1,de}^{abc}
 +S'_{z,\bar d}S'_{\bar x,\bar e}c'{}_{2,de}^{abc}
 +S'_{\bar z,d}S'_{\bar y,\bar e}c'{}_{3,de}^{abc}\Bigr)
 \Bigg]N_{zy}^{x}=0.
\end{gathered}
\end{equation}
\endgroup
Here $\alpha\in\{1,2,3\}$. The indices $i,j,k,\ell,p,q,r,s,x,y,z$
range over $I_+$, and $a,b,c,d,e,t,u,v,w$ range over $I_-$.
All equations hold for every choice of the free indices, and every
sum is over the corresponding label set.
\end{theorem}

\begin{proof}
On each side, the fusion structure equations give
\eqref{eq:pair-dim-dual}--\eqref{eq:pair-duality},
\eqref{eq:pair-dimension} and \eqref{eq:pair-associativity}.
Bi-commutativity gives \eqref{eq:pair-commutativity}.

The computations in Theorem~\ref{cons:algthm} identify
the native $(3,1)$, $(3,2)$-v and $(3,2)$-e comparisons
with the plus-side equations in
\eqref{eq:pair-31}--\eqref{eq:pair-32e}.
Applying the same computations on the minus side,
with $(N,d,c)$ replaced by $(N',d',c')$, gives
the corresponding primed equations.

Substituting the minimal idempotents into the
product-exchange identities of
Definition~\ref{def:ssprime-matched} and comparing
coefficients gives
\eqref{eq:final-fourier-plus-product} and
\eqref{eq:final-fourier-plus-convolution}.
The Fourier-square and unit-normalization identities give
\eqref{eq:pair-fourier-square} and
\eqref{eq:pair-fourier-unit}.
Closing the identity resolutions on the two sides gives
\eqref{eq:pair-modulus}.

By Lemma~\ref{prop:cprime-to-ctilde}, the $T'$ relation
is transported to the $W$ relation with coefficients
$\widetilde c$ given by
\eqref{eq:ctilde-1}--\eqref{eq:ctilde-3}.
Substituting these expressions into \eqref{cons:33e}
gives \eqref{eq:pair-33e}.
The corresponding comparison in the other shading is
obtained by rotation and the invertible Fourier changes
of labels, so it imposes the same condition.

Conversely, the displayed equations give the corresponding
local comparisons by the same coefficient calculations.
The minimal idempotents form bases of the two $2$-box
spaces, and the active triangles form bases of the
two $3$-box spaces by
Proposition~\ref{prop:triangle-basis}.
Thus the coefficient equalities are equivalent to the
local diagram equalities. The invertibility of Fourier
transport then identifies these with the comparisons
of Theorem~\ref{cons:algthm}, proving the equivalence.
\end{proof}

The native $(3,1)$ and $(3,2)$ comparisons follow from
the prior equations, as the following proposition shows.

\begin{proposition}
\label{prop:prior-redundancy}
Let $B_+$ and $B_-$ be fusion bialgebras with nonzero
quantum dimensions.
Suppose that the coefficient tables $c,c'$ satisfy
the prior equations on their respective fusion supports:
\[
\begin{gathered}
\left(
c_{1,xy}^{ijk}+c_{2,xz}^{ijk}+c_{3,yz}^{ijk}
-\delta_{xi}\delta_{yj}\delta_{zk}
\right)N_{zy}^{x}=0,\\
\left(
c'{}_{1,uv}^{abc}+c'{}_{2,uw}^{abc}+c'{}_{3,vw}^{abc}
-\delta_{ua}\delta_{vb}\delta_{wc}
\right)N'{}_{wv}^{u}=0.
\end{gathered}
\]
Then the $(3,1)$, $(3,2)$-v and $(3,2)$-e equations
\eqref{eq:pair-31}--\eqref{eq:pair-32e}
hold on both sides.
In particular, these equations hold for the canonical
tables read from $K_\pm$ when $H^2(K_\pm;\C)=0$
and $B_\pm$ realize the corresponding fusion supports.
\end{proposition}

\begin{proof}
We prove the unprimed equations. The same argument
applies on the minus side.

Multiplying the prior equation by $\delta_{x\ell}$ and using
$\delta_{x\ell}\delta_{xi}=\delta_{i\ell}\delta_{xi}$
gives the unprimed equation in \eqref{eq:pair-32v}.

Next, set $x=0$, $y=\bar s$ and $z=s$ in prior.
Since $N_{s\bar s}^{0}=d_s\ne0$, this gives
\[
c_{1,0\bar s}^{ijk}
+c_{2,0s}^{ijk}
+c_{3,\bar s s}^{ijk}
=\delta_{i0}\delta_{\bar s,j}\delta_{sk}
=\delta_{i0}\delta_{k\bar j}\delta_{sk},
\]
which is the unprimed equation in \eqref{eq:pair-31}.

For \eqref{eq:pair-32e}, set $z=s$ in prior,
multiply by $d_x$, and sum over $x,y$:
\[
\sum_{x,y}d_xN_{sy}^{x}
\left(
c_{1,xy}^{ijk}+c_{2,xs}^{ijk}+c_{3,ys}^{ijk}
\right)
=\delta_{sk}d_iN_{sj}^{i}.
\]
Frobenius reciprocity and the dimension identities give
\[
\begin{gathered}
d_xN_{sy}^{x}=d_sN_{x\bar y}^{s},\\
\sum_y d_xN_{sy}^{x}
=\sum_y d_yN_{\bar s x}^{y}
=d_sd_x,
\qquad
\sum_x d_xN_{sy}^{x}=d_sd_y,\\
\delta_{sk}d_iN_{sj}^{i}
=\delta_{sk}d_sN_{i\bar j}^{k}.
\end{gathered}
\]
Substituting these identities and dividing by $d_s$ yields
\[
\sum_{x,y}c_{1,xy}^{ijk}N_{x\bar y}^{s}
+\sum_x c_{2,xs}^{ijk}d_x
+\sum_y c_{3,ys}^{ijk}d_y
=\delta_{sk}N_{i\bar j}^{k}.
\]
This is the unprimed equation in \eqref{eq:pair-32e}.
Multiplying by $d_s$ and summing over $s$ also gives
\begin{equation}\label{cons:Qclosure}
\sum_{x,y}
\left(
c_{1,xy}^{ijk}+c_{2,xy}^{ijk}+c_{3,xy}^{ijk}
\right)d_xd_y
=d_kN_{i\bar j}^{k}.
\end{equation}

Finally, the canonical tables obtained from
$H^2(K_\pm;\C)=0$ satisfy prior by
Proposition~\ref{prop:prior-as-coboundary}
and Theorem~\ref{thm:hodge-gauge}.
The preceding calculations therefore apply to these tables.
\end{proof}

Starting with $(K_+,K_-)$ satisfying $H^2(K_\pm;\C)=0$,
we obtain the rational coefficient tables $c,c'$ from
Theorem~\ref{thm:canonical-readout}.
By Proposition~\ref{prop:prior-redundancy}, the
$(3,1)$ and $(3,2)$ comparisons follow from prior.
The following theorem reconstructs a bi-commutative
TR planar algebra from compatible fusion data and
Fourier matching satisfying the common-edge equation,
with the prescribed fusion supports.

\begin{theorem}\label{thm:construction-after-consistency}
Let $(K_+,K_-)$ be a fusion complex pair with
$H^2(K_\pm;\mathbb Q)=0$, and let $c,c'$ be the canonical
coefficient tables of Theorem~\ref{thm:canonical-readout}.
Suppose that the data $(N,N',d,d',S,S',\delta)$ satisfy
$d_0=d'_0=1$ and the following equations of
Theorem~\ref{thm:final-consistency-system}:
\[
\eqref{eq:pair-dim-dual}--
\eqref{eq:pair-associativity},
\qquad
\eqref{eq:final-fourier-plus-product}--
\eqref{eq:pair-fourier-unit},
\qquad
\eqref{eq:pair-modulus},
\qquad
\eqref{eq:pair-33e}.
\]
Assume that the fusion coefficients have exactly the
supports specified by $K_+,K_-$, that all quantum
dimensions and $\delta$ are nonzero, and that $S,S'$
are invertible.
Let $B_+,B_-$ be the fusion bialgebras defined by
$(N,d)$ and $(N',d')$, respectively.
Then there exists a bi-commutative TR planar algebra
$P$ with $P_{2,\pm}\cong B_\pm$.
Its fusion coefficients are $N,N'$, its quantum
dimensions are $d,d'$, its native TR relations have
coefficient tables $c,c'$, and its Fourier matrices
are $S,S'$.
\end{theorem}

\begin{proof}
First form the abstract plus-side vector space with basis
$\{p_i:i\in I_+\}$ and prescribe
\[
 p_ip_j=\delta_{ij}p_i,\qquad
 1_+=\sum_i p_i,\qquad
 p_i*_+p_j=\sum_kN_{ij}^{k}p_k,\qquad
 \operatorname{tr}_+(p_i)=d_i.
\]
The plus-side fusion structure equations, together with $d_0=1$ and the nonvanishing assumptions, define the fusion bialgebra $B_+$. Its total dimension is nonzero by \eqref{eq:pair-modulus}.
For each distinguished basis element $p_i$, introduce a
formal plus-shaded $2$-box with the boundary markings
fixed by our graphical conventions.
Use the multiplication, convolution and trace of $B_+$
to prescribe the ordinary-product, normalized-convolution,
cap and partial-trace relations, with the conventions of
Section~\ref{subsec:presentations}.
A contractible circle evaluates to $\delta$,
$p_0$ is identified with the normalized Jones diagram,
and a half-turn sends $p_i$ to $p_{\bar i}$.
We will prove below that the images of these labelled
$2$-boxes in the resulting planar algebra are nonzero
and linearly independent.

We define the minus labels by actual rotation of plus-side diagrams.
Let $\mathcal F_+$ and $\mathcal F_-$ denote the counterclockwise
quarter-turn tangles from plus to minus and from minus to plus,
respectively. Thus $\mathcal F_-^{-1}$ is the clockwise quarter-turn
from plus to minus. Put
\begin{equation}\label{eq:reconstruction-minus-labels}
 x_a:=\sum_{i\in I_+}S'_{ia}p_i,\qquad
 p'_a:=\mathcal F_-^{-1}(x_a)
      =\sum_{i\in I_+}S'_{ia}\mathcal F_-^{-1}(p_i).
\end{equation}
The label $p'_a$ denotes this linear combination of
rotated plus-side boxes, with boundary markings
inherited from the rotation. In particular,
\[
 \mathcal F_-(p'_a)=\sum_iS'_{ia}p_i.
\]
Using $S'S=C_{\mathrm{dual},+}$ from
\eqref{eq:pair-fourier-square} and the prescribed
half-turn gives
\[
 \sum_aS_{ai}p'_a
 =\mathcal F_-^{-1}\!\left(\sum_{j,a}S'_{ja}S_{ai}p_j\right)
 =\mathcal F_-^{-1}(p_{\bar i})
 =\mathcal F_+(p_i).
\]
The other identity in \eqref{eq:pair-fourier-square},
$SS'=C_{\mathrm{dual},-}$, similarly gives
$\mathcal F_+\mathcal F_-(p'_a)=p'_{\bar a}$.

We next verify the low-box rules for the minus labels.
The right-hand equations in
\eqref{eq:final-fourier-plus-product} and
\eqref{eq:final-fourier-plus-convolution} give, in $B_+$,
\begin{equation}\label{eq:reconstruction-x-products}
 x_a*_+x_b=\delta\,\delta_{ab}x_a,\qquad
 \delta x_ax_b=\sum_cN'{}_{ab}^{c}x_c.
\end{equation}
The geometric product-exchange identities, with our normalization,
are
\[
 \mathcal F_-(uv)=\delta^{-1}\mathcal F_-(u)*_+\mathcal F_-(v),
 \qquad
 \mathcal F_-(u*_-v)=\delta\mathcal F_-(u)\mathcal F_-(v).
\]
Applying these tangle identities to
\eqref{eq:reconstruction-minus-labels} and
\eqref{eq:reconstruction-x-products} yields the diagrammatic rules
\[
 p'_ap'_b=\delta_{ab}p'_a,\qquad
 p'_a*_-p'_b=\sum_cN'{}_{ab}^{c}p'_c.
\]

Since $S'$ is invertible, the $x_a$ form a basis of the abstract
space $B_+$. Their convolution relations imply
$(\sum_a x_a)*_+x_b=\delta x_b$ for every $b$. As $p_0$ is the
convolution unit, this proves
\[
 \sum_a x_a=\delta p_0.
\]
Rotating and using the Jones diagram gives
\[
 \sum_a p'_a=\mathcal F_-^{-1}(\delta p_0)=1_-.
\]
The unit normalization
\eqref{eq:pair-fourier-unit} gives
\[
 \mathcal F_+(p_0)
 =\sum_a S_{a0}p'_a
 =\delta^{-1}1_-,
 \qquad
 \mathcal F_-(p'_0)
 =\sum_i S'_{i0}p_i
 =\delta^{-1}1_+.
\]
Hence
$p'_0=\delta^{-1}\mathcal F_-^{-1}(1_+)$
is the normalized minus Jones diagram.

We next verify that the closure of $p'_a$ has the
prescribed value $d'_a$.
Frobenius reciprocity in
Proposition~\ref{prop:basic-2box-equations}
and the dimension equations
\eqref{eq:pair-dimension} give
\[
 d'_c\sum_b N'{}_{ab}^{c}
 =\sum_b d'_b N'{}_{\bar a c}^{b}
 =d'_a d'_c.
\]
Since $d'_c\ne0$, we obtain
$\sum_b N'{}_{ab}^{c}=d'_a$.
Using $\sum_b x_b=\delta p_0$, sum the second
identity of \eqref{eq:reconstruction-x-products}
over $b$:
\[
 \delta^2 x_a p_0
 =\sum_c\left(\sum_b N'{}_{ab}^{c}\right)x_c
 =\delta d'_a p_0.
\]
Comparing the $p_0$ coefficients in $B_+$ gives
$\delta S'_{0a}=d'_a$.
The plus-side cap rules give
\[
 \operatorname{tr}_-\bigl(\mathcal F_-^{-1}(p_i)\bigr)
 =\delta\delta_{i0}.
\]
Thus \eqref{eq:reconstruction-minus-labels} gives
\[
 \operatorname{tr}_-(p'_a)
 =\delta S'_{0a}
 =d'_a.
\]
The corresponding
cap and partial-trace rules are obtained by the same rotated
low-box computations. These verifications establish the minus-side
rules in the presentation. 

We next form the $3$-box diagrams from these labelled
$2$-boxes. Inserting the labels $p_i$ into the tangles
of Section~\ref{sec:TR-definitions} gives
$T(i,j,k)$, $W(i,j,k)$ and $E_\alpha(\ell,m)$.
Using the minus-side labels $p'_a$ in the corresponding
minus-shaded tangles gives $T'(a,b,c)$ and
$E'_\alpha(d,e)$, with the conventions of
\eqref{eq:minus-tr-expansion-picture}.

The canonical tables $c,c'$ specify the local
triangle reduction rules
\[
 T(i,j,k)\longmapsto
 \sum_{\alpha=1}^{3}\sum_{\ell,m}
 c_{\alpha,\ell m}^{ijk}E_\alpha(\ell,m),
 \qquad
 T'(a,b,c)\longmapsto
 \sum_{\alpha=1}^{3}\sum_{d,e}
 c'{}_{\alpha,de}^{abc}E'_\alpha(d,e).
\]
Applying the Fourier changes of labels and the
rotation in Lemma~\ref{prop:cprime-to-ctilde}
to the second rule gives the reduction rule for
$W(i,j,k)$, with coefficients $\widetilde c$.
Substituting these coefficients into
\eqref{cons:33e} gives \eqref{eq:pair-33e}.

We now verify the algebraic equations associated with
the local consistency comparisons.
Equation~\eqref{eq:pair-dim-dual} gives the $(1,1)$
comparison, and \eqref{eq:pair-dimension} gives the
$(2,1)$ comparison.
The vacuum and associativity equations
\eqref{eq:pair-vacuum} and \eqref{eq:pair-associativity},
together with ordinary multiplication and the
convolution unit relations, give the two $(2,2)$
comparisons.
Equation~\eqref{eq:pair-duality} makes convolution
compatible with the half-turn.
These statements hold on both sides.

The canonical tables satisfy prior by
Theorem~\ref{thm:canonical-readout}.
Proposition~\ref{prop:prior-redundancy} then gives
the $(3,1)$, $(3,2)$-v and $(3,2)$-e equations
\eqref{eq:pair-31}, \eqref{eq:pair-32v} and
\eqref{eq:pair-32e}, respectively, on both sides.
The Fourier transport formulas in
Lemma~\ref{prop:cprime-to-ctilde} convert the
minus-side equations to the common plus shading.
Under the same coefficient substitution, the
$(3,3)$-e equation takes the form
\eqref{eq:pair-33e}, which holds by assumption.

These algebraic facts will be used in the induction
below to prove equality of the corresponding closed
evaluations. Any identity resolutions used in the
coefficient calculations must also be justified
within that induction.

Let $Z$ be the average canonical evaluation of
\eqref{eq:average-evaluation}, extended linearly.
Write $v(D)$ for the number of labelled $2$-box vertices.
The Fourier changes of labels replace each labelled
$2$-box by a linear combination of single labelled
$2$-boxes, preserving the vertex count in every term.

We prove canonical consistency and compatibility with identity
resolution simultaneously. The local comparisons for canonical
consistency use identity resolutions on diagrams with fewer
vertices, while the identity-resolution argument uses canonical
reductions on smaller diagrams.
More precisely, we prove the following assertions for every
$q\ge0$:
\begin{enumerate}[label=\textup{(\roman*)},leftmargin=2.6em]
\item[$(A_q)$] Every closed diagram with at most $q$ vertices has
path-independent canonical evaluation.
\item[$(U_q)$] For every compatible complementary $2$-box diagram
$C$ with $v(C)+1\le q$,
\[
 Z\bigl(\langle U,C\rangle\bigr)
 =\sum_{i\in I_+}Z\bigl(\langle p_i,C\rangle\bigr),
\]
where $U$ is the ordinary identity diagram. The analogous assertion
holds in the other shading.
\end{enumerate}
Both assertions include disconnected diagrams.
In $(U_q)$, the diagram $\langle U,C\rangle$ has
$v(C)$ labelled vertices, whereas each
$\langle p_i,C\rangle$ has $v(C)+1$.
Thus the bound $v(C)+1\le q$ applies to every
term after the identity resolution.

The assertion $(A_0)$ follows from the circle rule, and $(U_0)$
is vacuous. Vertex-free circles contribute factors of $\delta$ and
their removal commutes with all other reductions. We handle them
separately. Every nonzero term of a move of type $1$, $2$,
or $3$ has fewer vertices than its input.

We first check the identity resolution when the
complementary diagram belongs to $O_2$.
For a labelled vertex $p_j$, the required equality
follows from
$
 \left(\sum_i p_i\right)p_j=p_j.
$
The cap and cup configurations are evaluated using
the prescribed cap and Jones relations.
Closing the identity diagram into two circles gives
$\delta^2$, while closing its expansion gives
$\sum_i d_i$; these agree by
\eqref{eq:pair-modulus}.
The same computations apply on the minus side,
using the rules verified above.
These checks establish the required equality for
the $O_2$ diagrams. The induction below extends it
to arbitrary complementary diagrams.

Suppose that $(A_{q-1})$ and $(U_{q-1})$ hold. To prove $(A_q)$,
let $D$ have $q$ vertices and choose two canonical first moves
$m,k$, with outputs $D_m,D_k$. After the separate handling of
circles, every output term has at most $q-1$ vertices.

If the two supports are disjoint, the moves commute as formal
linear combinations of diagrams. Hence $(A_{q-1})$ gives
\[
 Z(D_m)=Z(D_{m,k})=Z(D_{k,m})=Z(D_k).
\]
If they overlap, write $D=\langle H,B\rangle$, with $H$ the local
input and $B$ its unchanged complement. Thus
$D_m=\langle H_m,B\rangle$ and $D_k=\langle H_k,B\rangle$.
We separate the local comparisons according to whether
their outputs are $2$-box or $3$-box diagrams, following
the computations in the proof of
Theorem~\ref{cons:algthm}.

For $(3,1)$ and $(3,2)$-e, the overlap diagram $H$
has three vertices. After either first reduction,
each local term has at most two vertices.
The local computations reduce the two outputs to
\[
 H_m\longmapsto\sum_s a_s p_s,
 \qquad
 H_k\longmapsto\sum_s b_s p_s.
\]
The corresponding equations
\eqref{eq:pair-31} and \eqref{eq:pair-32e}
state precisely that $a_s=b_s$ for every $s$.
Thus the two resulting linear combinations of
labelled diagrams are identical.
Since $v(B)=q-3$, all reductions after the first
step take place in closed diagrams with at most
$q-1$ vertices. By $(A_{q-1})$,
\[
 Z(D_m)
 =\sum_s a_s Z(\langle p_s,B\rangle)
 =\sum_s b_s Z(\langle p_s,B\rangle)
 =Z(D_k).
\]

For $(3,2)$-v and $(3,3)$-e, the overlap diagram $H$
has four vertices. 
We now justify the diagrammatic steps used to express
the difference of the two outputs in triangle form.
For $(3,2)$-v and $(3,3)$-e, put $b=v(B)$.
The overlap diagram $H$ has four vertices, so $q=b+4$.
After either first reduction, each local term has
at most three vertices. Thus all subsequent canonical
reductions take place in closed diagrams with at most
$q-1$ vertices and preserve $Z$ by $(A_{q-1})$.

For $(3,2)$-v and $(3,3)$-e, we follow the two
calculation chains displayed in the proof of
Theorem~\ref{cons:algthm}. Starting from $H_m$ and
$H_k$, these computations express the outputs as
linear combinations of triangle diagrams, using
canonical reductions and the identity resolutions
\eqref{cons:E1}.

Put $b=v(B)$, so $q=b+4$.
Inspection of the displayed diagrams shows that
every local term after the first reduction has
at most three vertices. After closure against $B$,
each therefore has at most $b+3=q-1$ vertices,
so $(A_{q-1})$ guarantees that the canonical
reductions preserve $Z$.
For the identity resolutions, each $E_\alpha$
diagram has two local vertices, whereas each
triangle in its expansion has three.
Thus the closed diagrams before and after the
resolution have $b+2=q-2$ and $b+3=q-1$ vertices,
respectively.
The complement $C$ of the unit-resolution disc
has $v(C)=q-2$, so $v(C)+1=q-1$.
Hence $(U_{q-1})$ guarantees that these steps
also preserve $Z$.

Following the two calculation chains therefore gives
\[
 Z(\langle H_m,B\rangle)
 -Z(\langle H_k,B\rangle)
 =
 \sum_{x,y,z}R_{xyz}
 Z(\langle T(x,y,z),B\rangle).
\]
The corresponding equations
\eqref{eq:pair-32v} and \eqref{eq:pair-33e}
give $R_{xyz}N_{zy}^{x}=0$.
For an active triple, $N_{zy}^{x}\ne0$, so
$R_{xyz}=0$.
For an inactive triple, the prescribed rule is
$T(x,y,z)\longmapsto0$, so the entire closed diagram
$\langle T(x,y,z),B\rangle$ reduces to zero.
By $(A_{q-1})$, its evaluation is therefore zero.

The $(3,3)$-v type is handled by the graphical redundancy argument
in the proof of Theorem~\ref{cons:graphthm}. At the present vertex
count, that argument connects the two first moves by disjoint
comparisons and the other overlap types just established. Its subsequent
canonical reductions are on smaller diagrams and are justified by
$(A_{q-1})$. 

All canonical first moves on $D$ now give the same output value.
After the first move, $(A_{q-1})$ makes every subsequent complete
canonical evaluation path independent. This proves $(A_q)$.
Since $Z$ is the average of complete canonical evaluations, every
canonical move on a closed diagram with at most $q$ vertices
now preserves $Z$.

To prove $(U_q)$, let $v(C)+1\le q$ and remove vertex-free circles
first. A closed component with vertices has a reducible
face by the closed-diagram small-face argument. After excluding
such components, the boundary components form either a connected
$2$-box or two connected $1$-boxes. If an internal face remains,
Lemma~\ref{cons:lowbox} provides a reducible internal face.
If none remains, the comparison is one of the low-box $O_2$ tests
already verified.

In the reducible case, choose one move entirely inside $C$,
and write
\[
 C\longrightarrow\sum_j a_jC_j.
\]
Each nonzero term satisfies $v(C_j)+1\le q-1$. Apply the same
move after closing with $U$ and after closing with each $p_i$.
All these input closures have at most $q$ vertices, so $(A_q)$
and then $(U_{q-1})$ give
\[
\begin{aligned}
 Z(\langle U,C\rangle)
 &=\sum_j a_jZ(\langle U,C_j\rangle)\\
 &=\sum_j a_j\sum_iZ(\langle p_i,C_j\rangle)\\
 &=\sum_iZ(\langle p_i,C\rangle).
\end{aligned}
\]
The first and third equalities use $(A_q)$, and the second uses
$(U_{q-1})$. In particular, the third equality reads the already
proved canonical-reduction equality in reverse on the closures
with $p_i$. The opposite shading is identical. This proves $(U_q)$
and completes the simultaneous induction.

It follows that $Z$ is invariant under Moves $0$--$3$ by the
assertions $(A_q)$ and under Move~4 by the assertions $(U_q)$.
Thus $Z$ is a partition function with $Z(\varnothing)=1$.

Let $V_{n,\varepsilon}$ be the space of labelled
$n$-box diagrams of shading $\varepsilon$. Define
\[
\mathcal N_{n,\varepsilon}
=
\left\{
x\in V_{n,\varepsilon}:
Z(\langle x,B\rangle)=0
\text{ for every compatible complementary diagram }B
\right\}.
\]
These subspaces form a planar ideal: closing a tangle
containing $x$ produces a complementary closure of $x$.
Following \cite[Proposition~1.24]{JonesPA}, set
\[
P_{n,\varepsilon}
=
V_{n,\varepsilon}/\mathcal N_{n,\varepsilon}.
\]
The sphere evaluation makes $P$ spherical, and its
trace pairings are nondegenerate by
\cite[Proposition~1.29]{JonesPA}.
The normalization and low-box reductions give
$P_{0,\pm}\cong\C$ and $\dim P_{1,\pm}=1$.
The identity strand survives because its closure
is $\delta\ne0$.

Low-box reduction and unit resolution show that the images of $p_i$
span $P_{2,+}$. Moreover,
$
 \operatorname{tr}_+(p_ip_j)=\delta_{ij}d_i.
$
Since every $d_i$ is nonzero, these images are linearly independent.
Thus $P_{2,+}\cong B_+$ with its prescribed ordinary product,
convolution, trace, and duality. Likewise the graphical labels $p'_a$
span $P_{2,-}$, and
$
 \operatorname{tr}_-(p'_ap'_b)=\delta_{ab}d'_a
$
shows that they are nonzero and linearly independent. Together with
their orthogonal idempotent relations, this proves that they form
the full set of primitive idempotents and that $P_{2,-}\cong B_-$.
In particular, both ordinary $2$-box algebras are commutative.

The actual geometric quarter-turns have matrices $S,S'$ by
\eqref{eq:reconstruction-minus-labels} and the Fourier-square
calculation above. The prescribed native triangle relations have
tables $c,c'$, and their transport gives the opposite-orientation
relation. Hence $P$ is TR. The exact-support hypotheses identify its
fusion complexes with $K_+,K_-$. Finally, unit resolution expresses
each $E_\alpha$ as a sum of $T$-triangles, so the transported
$W$ relation also makes $P$ YBR.

% AUTHOR CHECK (retained from the supplied manuscript): identify the
% precise finite-dimensionality theorem and verify its hypotheses.
% The radical quotient by itself does not establish local finite dimensionality.
By construction, $P$ is generated by its $2$-box spaces,
which are finite-dimensional since $P_{2,\pm}\cong B_\pm$.
The triangle relations, together with the unit resolutions,
give the Yang--Baxter relations in both shadings.
Therefore, the finite-spanning argument in the proof of
\cite[Theorem~3.4]{LiuYBR} applies, and every box space
$P_{n,\pm}$ is finite-dimensional.
\end{proof}

\begin{proposition}
\label{prop:matched-uniqueness}
Let $P$ and $Q$ be two $2$-box-generated,
trace-nondegenerate realizations of the same labelled
local data, TR relations and Fourier matching.
Then the prescribed correspondence of labelled generators
extends to a unique planar-algebra isomorphism $P\cong Q$.
\end{proposition}

\begin{proof}
Canonical evaluations with the same local coefficients
give the same value for every closed labelled diagram.
Thus the two realizations have the same partition function
$Z$ and the same complementary closure radical $\mathcal N$.

The maps from the labelled diagram space to $P$ and $Q$
are surjective by generation. Nondegeneracy implies that
both kernels are exactly $\mathcal N$.
Hence both realizations are the quotient determined by $Z$,
as in \cite[Proposition~1.24]{JonesPA}.
The induced isomorphism is unique because the labelled
$2$-boxes generate both planar algebras.
\end{proof}

\begin{remark}
For a fixed fusion complex pair with
$H^2(K_\pm;\Q)=0$, the canonical readout fixes
$c,c'$, and trace closure determines $N$ from $d$.
The character-table construction then recovers
the remaining numerical data up to finitely many
algebraic choices.
Thus the continuous parameters are controlled by
the $n-1$ dimension coordinates $d_1,\ldots,d_{n-1}$.

Using the root coordinates described in the
algorithm below, each chart requires at most
$n-2$ additional root variables and one variable
for the nonvanishing conditions.
With $\delta$ recovered afterwards, this gives
at most $2n-2=O(n)$ polynomial variables.
The reconstruction system has $O(n^7)$ indexed
scalar equations, the largest family being the
common-edge equations with seven free labels.
For $n=4,5$, this leaves at most six and eight
polynomial variables per chart, respectively.
The cohomological readout removes the $6n^5$
triangle coefficient entries from the nonlinear
search, making exact low-rank computations practical.
In our implementation, the rank-two and rank-three
computations took less than one minute each,
and the complete rank-four computation took
approximately 34 minutes.
In rank five, the self-dual branch took several hours,
while the branches with duality $(12)$ and $(12)(34)$
took approximately 73 and 11 minutes, respectively.
These timings include verification and describe
the observed performance of our implementation,
rather than theoretical complexity bounds.
\end{remark}

\subsection{Classification algorithm}
\label{sec:primitive}

\begin{definition}\label{def:primitive}
A TR planar algebra is \emph{primitive} if it is
neither an exchange-relation planar algebra nor isomorphic to a nontrivial
tensor product or ordered free product of lower-rank TR planar algebras
in that class.
\end{definition}

\label{sec:algorithm}

The preceding theory turns classification at fixed rank into finite
support enumeration, exact linear algebra, and polynomial equations on
explicit open sets. The four mathematical layers are
\[
 \text{support}\longrightarrow\text{fusion}
 \longrightarrow\text{Fourier matching}\longrightarrow\text{consistency}.
\]
Their purpose is to reduce the unknown data before imposing the final
consistency equations, rather than solve all structure constants
simultaneously.

Fix the rank $n$ and the label involution on each side.
The procedure retains all surviving realizations,
including ER objects. Primitive objects are identified
afterwards as in Section~\ref{sec:primitive}.
Each realization includes the fusion data on both
sides and their Fourier matching.
The output is a finite list of algebraic parameter
strata, together with the required nonvanishing
conditions and the equivalences described below.

\begin{proposition}
\label{prop:finite-classification-scheme}
For a fixed rank and prescribed label involutions,
Algorithms~1 and~2 give an exhaustive procedure for
classifying bi-commutative TR planar algebras with
these data, provided that all enumerations and
algebraic solution steps are completed.
Every point of a retained solution locus determines
a planar algebra by
Theorem~\ref{thm:construction-after-consistency},
unique in the sense of
Proposition~\ref{prop:matched-uniqueness}.
\end{proposition}

All polynomial computations are performed over $\Q$,
with solutions taken over $\C$.
Nonvanishing conditions are retained throughout:
a system $F=0$, $G\ne0$ may be represented by
$F=0$, $uG-1=0$, where $u$ is an auxiliary variable.
Whenever an elimination requires division by a
parameter-dependent expression, its zero locus is
treated separately unless that expression is already
required to be nonzero.

\begin{center}
\fbox{\begin{minipage}{0.94\linewidth}
\textbf{Algorithm 1. Classification at fixed rank}
\par\smallskip
\textbf{Input:} rank $n$ and the label involution
on each side, fixing $0$.
\par
\textbf{Output:} a finite collection of algebraic
solution loci with their nonvanishing conditions,
covering all realizations with the prescribed data.

\begin{enumerate}[leftmargin=1.8em,itemsep=4pt]
\item
Partition the label triples into orbits under
commutativity, duality and Frobenius reciprocity.
Use the unit and vacuum conditions to determine
which orbits are required or forbidden, and
enumerate all choices of the remaining orbits.

\item
Take representatives under label permutations
fixing $0$ and commuting with the prescribed
involution. For each resulting complex $K$,
form the matrix $D$ of $d_K:C^1(K;\Q)\to C^2(K;\Q)$.
Retain precisely those with
$\operatorname{rank}D=|K^{(2)}|$.

\item
Compute the Hodge right inverse
\[
C=D^{\mathsf T}(DD^{\mathsf T})^{-1}
\]
and verify $DC=I$.
Read off the canonical triangle coefficients,
assigning zero coefficients to inactive inputs.

\item
Use the quadratic closure formulas to express
$Q$ and $N$ in terms of the dimensions.
Solve the fusion equations with $d_0=1$,
$d_i\ne0$, and the exact-support conditions:
$Q_{ij}^{k}$ vanishes precisely at the inactive
positions. Retain every nonempty solution locus.

\item
Form all ordered pairs compatible with the rank,
face-number and duality constraints.
Keep the parameters on the two sides independent,
including when the two support complexes coincide.
Apply Algorithm~2 to each pair.

\item
Collect all nonempty solution loci returned by
Algorithm~2, together with their support pairs,
numerical data and nonvanishing conditions.
Organize the results by $3$-box dimension and duality.

\item
Organize the resulting realizations by
$3$-box dimension and duality.
Determine their ER and product decompositions,
where applicable.
\end{enumerate}
\end{minipage}}
\end{center}

\begin{center}
\fbox{\begin{minipage}{0.94\linewidth}
\textbf{Algorithm 2. Fourier matching and consistency}
\par\smallskip
\textbf{Input:} two fusion-data solution spaces
with exact supports $K_+,K_-$ and canonical
triangle tables $c,c'$.
\par
\textbf{Output:} all nonempty Fourier-matching
and consistency loci, with their numerical data
and nonvanishing conditions.

\begin{enumerate}[leftmargin=1.8em,itemsep=4pt]
\item
Choose one side for the character calculation,
denoted by the plus side, and put
$\tau=\sum_i d_i$.
Restrict to $\tau\ne0$ and use the scaled
Fourier matrices $X=\delta S$ and $Y=\delta S'$,
where $\delta^2=\tau$.

\item
Construct exact parameter charts for the
characters of $(B_+,*)$.
Separate the dimension character, and compute
the remaining characters from
$B_+/\mathbb C E$, where $E=\sum_i p_i$.
Use primitive-element charts and all ordered
distinct roots of their defining polynomials
to cover every invertible character table.

\item
Assign the character rows to the minus-side
labels, with $X_{0i}=d_i$ and $X_{a0}=1$.
Impose duality compatibility and the necessary
equations coming from the inactive minus-side
fusion coefficients.
Discard only loci proved empty.

\item
On each remaining locus, form
\[
\kappa_a=\sum_i\frac{X_{ai}X_{\bar a i}}{d_i},
\qquad
d'_a=\frac{\tau}{\kappa_a},
\qquad
Y_{ia}=\frac{d'_a}{d_i}X_{ai},
\]
requiring $\kappa_a\ne0$.
Require the recovered minus-side data to belong
to the given fusion-data solution space,
with its exact support and nonvanishing conditions.

\item
Impose all Fourier-matching conditions of
Definition~\ref{def:ssprime-matched}, including
$\sum_a d'_a=\tau$.
Verify the original product-exchange,
Fourier-square and unit-normalization equations
for $S=X/\delta$ and $S'=Y/\delta$.

\item
Transport $c'$ using
Lemma~\ref{prop:cprime-to-ctilde} and impose
the common-edge equation \eqref{eq:pair-33e}.
If a residual is not identically zero on a
matching locus, solve its vanishing equations
on that locus and retain every nonempty solution.

\item
Return the surviving data and their
nonvanishing conditions, retaining both choices
of $\delta$ satisfying $\delta^2=\tau$.
\end{enumerate}
\end{minipage}}
\end{center}

\paragraph{Root coordinates for character tables.}
We explain how the character tables in Step~2 of
Algorithm~2 are computed algebraically.
After separating the dimension character
$p_i\mapsto d_i$, consider the remaining
$(n-1)$-dimensional semisimple convolution algebra.
On each chart, choose an element $h$ that separates
its characters.
Then $1,h,\ldots,h^{n-2}$ form a basis, and the
images of the distinguished elements can be written
as $\phi_i(h)$, with coefficients rational in $d$
on that chart.

Let $q(t)$ be the characteristic polynomial of
multiplication by $h$.
An ordering of its roots $r_1,\ldots,r_{n-1}$
determines the remaining character rows by
$
\chi_a(p_i)=\phi_i(r_a).
$
We encode these roots by comparing coefficients in
$
q(t)=\prod_{a=1}^{n-1}(t-r_a).
$
The root-sum identity eliminates the last root,
leaving at most $n-2$ root variables.
Together with the dimension row, these character
rows give $X=\delta S$ after an admissible row matching.

\begin{proof}
We prove that the procedure retains the data of every
bi-commutative TR planar algebra of the prescribed rank
and duality types, and that every retained solution
admits the asserted realization.

Let $P$ be such a planar algebra.
After admissible relabellings, its two fusion supports
are among the enumerated representatives.
Theorem~\ref{thm:cohomological-criterion} gives
$H^2(K_\pm;\mathbb Q)=0$, so neither support is removed
by the cohomology test.
By Theorem~\ref{thm:canonical-readout}, the canonical
tables $c,c'$ give valid triangle expansions in $P$.
Closing these expansions gives the trace-closure
formulas for its fusion coefficients.
Consequently, substituting those formulas into the
fusion equations preserves the dimension vectors of $P$.
Its coefficients also satisfy the prescribed
nonvanishing and exact-support conditions.
Thus the two numerical fusion realizations of $P$
belong to the solution spaces retained by the fusion stage.

These two realizations form an admissible ordered pair.
For its Fourier matrices $S,S'$, the product-exchange
identities imply that the rows of $\delta S$ and
$\delta S'$ are convolution characters.
Invertibility ensures that they enumerate the
characters on their respective sides.
Hence the actual Fourier matching of $P$ is included
among the character-row matchings considered by the
algorithm.
It satisfies all Fourier-matching conditions and,
by Theorem~\ref{thm:final-consistency-system},
the common-edge equation \eqref{eq:pair-33e}.
The data of $P$ therefore remain in the final solution
locus. This proves that no target planar algebra is omitted.

Conversely, take any retained solution.
The fusion stage supplies the fusion-structure
identities and the required dimension and support
conditions. The matching stage supplies the common
modulus and all Fourier-matching conditions, and the
final stage imposes \eqref{eq:pair-33e}.
Moreover, the canonical tables satisfy prior by
Theorem~\ref{thm:canonical-readout}.
Proposition~\ref{prop:prior-redundancy} then gives the
$(3,1)$, $(3,2)$-e and $(3,2)$-v equations on both sides.
Thus all hypotheses of
Theorem~\ref{thm:construction-after-consistency}
are satisfied, and the retained solution determines
a bi-commutative TR planar algebra realizing its data.

Finally, Proposition~\ref{prop:matched-uniqueness}
gives uniqueness for fixed matched data and justifies
the identifications induced by admissible relabellings.

By Proposition~\ref{prop:links-fusion-graphs} and
the forest criterion of Lu and Liu \cite{LuLiuER},
a retained realization
is ER precisely when all its fusion links are forests.
To identify tensor and ordered free products, we first
compare the retained supports with those of products
of lower-rank models.
A support match provides a candidate identification.
We then seek parameter substitutions and compatible
label permutations on both sides, and verify equality
of the dimension vectors, convolution coefficients,
duality, modulus, and both Fourier matrices.
For parameter families, these comparisons are made
on the stated loci, including their nonvanishing
conditions, and the claimed coverage is checked.
When the hypotheses of
Proposition~\ref{prop:matched-uniqueness}
are satisfied, the matched-data identification
extends to a planar-algebra isomorphism.
\end{proof}

\paragraph{Code and computational data.}
The implementation and computational results are available at
\url{https://github.com/jiaqi-li-mathphy/TR-planar-algebras}.
Released version contains the materials accompanying this
version of the paper.

\section{Low-rank classification and examples}\label{sec:computations}
\label{sec:scope}
We classify connected, irreducible, spherical,
trace-nondegenerate, bi-commutative TR planar algebras
over $\C$, generated by their $2$-box spaces and with
every box space finite-dimensional, at ranks $2,3,4,5$.
Throughout this section, $n=\dim P_{2,+}=\dim P_{2,-}$ and
$m=\dim P_{3,+}=\dim P_{3,-}$. The triangle-basis theorem identifies
$m$ with the support face count. The classification is algebraic.
Positivity is not imposed.

All candidates pass through the support, fusion, Fourier,
and consistency layers before the surviving models are identified.
Tensor products and ordered free products are recognized and verified
after this complete filtering.
The resulting exchange-relation models agree with the corresponding
known classifications \cite{LuThesis,LuLiuER,LiuER}.

Self-dual means that every label is fixed by duality. For ranks three and four,
we call the other branch pair-dual. In rank five the branches are
self-dual, $(12)$, and $(12)(34)$.
We use $X=\delta S$ and $Y=\delta S'$. Permitted label changes transport
both matrices, and loop signs are retained. Tables are given up to
shading reversal.

\subsection{Ranks two and three}
\label{sec:small}
\begin{proposition}
Up to relabelling and shading reversal, there are
exactly two surviving support pairs in rank two
and seven in rank three.
The rank-three pairs comprise five in the self-dual
branch and two in the $(12)$ branch.
The corresponding planar algebras are precisely
the models in Table~\ref{tab:models}, with parameters
in the stated domains.
\end{proposition}
\begin{proof}
In rank two, both enumerated supports pass the
cohomology and fusion tests.
In rank three, the cohomology test retains nine
self-dual and three pair-dual supports. The fusion
equations retain six and two, respectively.

Fourier matching with the exact-support conditions
leaves two support pairs in rank two and,
up to shading reversal, seven in rank three:
five self-dual and two pair-dual.
Every surviving Fourier locus satisfies
the common-edge consistency equation.
The resulting data are identified with the models
in Table~\ref{tab:models} on the stated parameter domains.

Proposition~\ref{prop:finite-classification-scheme},
Theorem~\ref{thm:construction-after-consistency},
and Proposition~\ref{prop:matched-uniqueness}
give the asserted classifications.
\end{proof}

For the model names, write $T(a)$ for the rank-two Temperley--Lieb data with nontrivial quantum dimension $a$ and loop square $1+a$. The loop sign is retained separately. In particular, $T(1)$ denotes the index-two specialization with loop square $2$. The symbol $*$ between models denotes ordered free product. In Table~\ref{tab:models}, every unfixed TL parameter lies in $\mathbb C\setminus\{0,-1,1\}$. Subscripts on the calligraphic complex symbols denote face counts. The explicit supports appear in Appendix~\ref{app:supports}.

\begin{table}[htbp]
\centering\small
\begin{tabular}{crcllr}
\toprule
$n$ & $m$ & Duality & Support pair & Model & Dim.\\\midrule
2&4&self-dual&$(\mathcal T_4,\mathcal T_4)$&$T(1)$&0\\
2&5&self-dual&$(\mathcal T_5,\mathcal T_5)$&$T(a)$&1\\
3&9&pair-dual&$(\mathcal Z_9,\mathcal Z_9)$&$\mathbb Z_3$&0\\
3&10&self-dual&$(\mathcal F_{10},\mathcal F_{10})$&$T(1)*T(1)$&0\\
3&11&self-dual&$(\mathcal F_{11}^L,\mathcal F_{11}^R)$, reverse&$T(1)*T(b)$, reverse&1\\
3&12&self-dual&$(\mathcal F_{12},\mathcal F_{12})$&$T(a)*T(b)$&2\\
3&13&self-dual&$(\mathcal D_{13},\mathcal D_{13})$&$\mathbb Z_2\subset\mathbb Z_5\rtimes\mathbb Z_2$&0\\
3&13&pair-dual&$(\mathcal E_{13},\mathcal E_{13})$&Binary-tree $\mathcal P^B$&0\\
3&14&self-dual&$(\mathcal W_{14},\mathcal W_{14})$&BMW&1\\\bottomrule
\end{tabular}
\caption{All surviving support-pair strata. The last column is the dimension of the parameter locus.}\label{tab:models}
\end{table}

\subsection{Rank four}
\label{sec:four}
Let $\mathcal Z_q$ denote the cyclic group planar algebra, and put
\[
\mathcal D_q=\mathcal S^{\mathbb Z_2\subset\mathbb Z_q\rtimes\mathbb Z_2}
\quad(q=5,7),\qquad \BMW(T)=\text{the rank-three BMW model in Appendix~\ref{app:rank3models}}.
\]
Here $\PB$ is Liu-Lu binary-tree exchange relation planar algebra. All unfixed TL parameters in the next table belong to
$\mathcal U=\mathbb C\setminus\{0,-1,1\}$.
The parameter $T$ of $\BMW(T)$ belongs to the exact fourteen-face domain given in Appendix~\ref{app:rank3models}. In particular, its dihedral specializations are listed in their own rows.

\begin{theorem}\label{thm:four}\label{prop:four}
Every bi-commutative TR planar algebra with $\dim P_2=4$
is, up to planar-algebra isomorphism and shading reversal,
one of the models listed in Table~\ref{tab:four},
with parameters in the stated domains.
Conversely, every listed model occurs throughout its
stated parameter domain.
There are no primitive rank-four objects in the sense of
Definition~\ref{def:primitive}.
\end{theorem}
\begin{table}[H]
\centering\small
\renewcommand{\arraystretch}{1.16}
\begin{tabular}{c p{.54\textwidth} p{.28\textwidth}}
\toprule
$m$ & Self-dual & Pair-dual\\
\midrule
16 & $T(1)\otimes T(1)=\mathcal Z_2\otimes\mathcal Z_2$
   & $\mathcal Z_4$\\
18 & --- & $\mathcal Z_3*T(1)$\\
19 & $T(1)*T(1)*T(1)$ & $\mathcal Z_3*T(t)$\\
20 & $T(1)*T(b)*T(1)$\newline
     $T(a)*T(1)*T(1)$\newline
     $T(1)\otimes T(b)$
   & ---\\
21 & $T(1)*T(b)*T(c)$\newline
     $T(a)*T(1)*T(c)$
   & ---\\
22 & $\mathcal D_5*T(1)$\newline
     $T(a)*T(b)*T(c)$
   & $\PB*T(1)$\\
23 & $\BMW(T)*T(1)$\newline
     $\mathcal D_5*T(t)$
   & $\PB*T(t)$\\
24 & $\BMW(T)*T(t)$ & ---\\
25 & $\mathcal D_7$\newline
     $T(a)\otimes T(b)$
   & ---\\
\bottomrule
\end{tabular}
\caption{Complete rank-four classification by $3$-box dimension.
Each product is ordered when appropriate.
Its shading reversal is included.
Every unfixed TL parameter lies in $\mathcal U$.
A dash indicates that no objects occur.}
\label{tab:four}
\end{table}
\begin{proof}
At rank four, support enumeration and the cohomology
test leave 186 single-sided candidates, of which
47 supports admit solutions of the fusion equations.
Fourier matching and common-edge consistency,
with the exact-support conditions retained, leave
19 support pairs up to shading reversal:
14 in the self-dual branch and 5 in the pair-dual branch.
Their numerical data are identified with the models
in Table~\ref{tab:four} on the stated parameter domains.
The enumeration and model identifications are recorded
in the accompanying computational certificates.

Proposition~\ref{prop:finite-classification-scheme}
therefore gives the classification.
Every listed model is an exchange-relation planar
algebra or a tensor or ordered free product of
lower-rank examples, so none is primitive.
\end{proof}
\begin{remark}
The parameter specializations are listed on their actual support strata.
Consistency excludes the additional Fourier-compatible cases with
$m=27,29$. The $m=25$ dihedral and tensor-product models survive.
\end{remark}

\subsection{Rank five}
\label{sec:five}
\label{prop:five}
\begin{theorem}\label{thm:five}
Every bi-commutative TR planar algebra with $\dim P_2=5$
is, up to planar-algebra isomorphism and shading reversal,
one of the models listed in Table~\ref{tab:five},
with parameters in the stated domains.
Conversely, every listed model occurs throughout its
stated parameter domain.
\end{theorem}

\begingroup\small
\setlength{\tabcolsep}{3pt}
\setlength{\LTcapwidth}{\textwidth}
\renewcommand{\arraystretch}{1.16}

\begin{longtable}{@{}r p{6.1cm} p{5.1cm} p{3.1cm}@{}}

% 第一页表头
\toprule
$m$ & Self-dual & $(12)$ & $(12)(34)$\\
\midrule
\endfirsthead

% 后续各页表头
\multicolumn{4}{@{}l}{Table \thetable\ (continued)}\\[4pt]
\toprule
$m$ & Self-dual & $(12)$ & $(12)(34)$\\
\midrule
\endhead

% 非最后一页的底线
\bottomrule
\endfoot

% 最后一页的底线与说明
\bottomrule
\addlinespace[8pt]
\caption[Complete rank-five classification.]
{\normalsize
Complete rank-five classification by $3$-box dimension.
Every unfixed TL parameter lies in $\mathcal U$.
Each BMW parameter lies in the domain specified in
Appendix~\ref{app:rank3models}.
The Hadamard-type families have $a\notin\{0,1,2\}$
in the $m=35$ rows and $a=2$ in the $m=34$ rows.
Free products are ordered, and shading reversals are included.
A dash indicates that no objects occur.}
\label{tab:five}
\label{tab:five-self}
\label{tab:five12}
\label{tab:five1234}\\
\endlastfoot
25 & --- & --- & \begin{tabular}[t]{@{}l@{}}$\mathcal Z_5$\end{tabular}\\[4pt]
28 & \begin{tabular}[t]{@{}l@{}}$(T(1)\otimes T(1))*T(1)$\end{tabular} & \begin{tabular}[t]{@{}l@{}}$\mathcal Z_4*T(1)$\end{tabular} & ---\\[4pt]
29 & \begin{tabular}[t]{@{}l@{}}$(T(1)\otimes T(1))*T(t)$\end{tabular} & \begin{tabular}[t]{@{}l@{}}$\mathcal Z_4*T(t)$\end{tabular} & \begin{tabular}[t]{@{}l@{}}$\mathcal Z_3*\mathcal Z_3$\end{tabular}\\[4pt]
30 & --- & \begin{tabular}[t]{@{}l@{}}$\mathcal Z_3*T(1)*T(1)$\\$T(1)*\mathcal Z_3*T(1)$\end{tabular} & ---\\[4pt]
31 & \begin{tabular}[t]{@{}l@{}}$T(1)*T(1)*T(1)*T(1)$\end{tabular} & \begin{tabular}[t]{@{}l@{}}$\mathcal Z_3*T(t)*T(1)$\\$\mathcal Z_3*T(1)*T(t)$\\$T(t)*\mathcal Z_3*T(1)$\end{tabular} & ---\\[4pt]
32 & \begin{tabular}[t]{@{}l@{}}$T(1)*T(1)*T(1)*T(t)$\\$(T(1)\otimes T(t))*T(1)$\\$T(1)*T(1)*T(t)*T(1)$\end{tabular} & \begin{tabular}[t]{@{}l@{}}$\mathcal Z_3*T(a)*T(b)$\\$T(a)*\mathcal Z_3*T(b)$\end{tabular} & ---\\[4pt]
33 & \begin{tabular}[t]{@{}l@{}}$T(1)*T(a)*T(1)*T(b)$\\$T(1)*T(1)*T(a)*T(b)$\\$T(1)*T(a)*T(b)*T(1)$\\$(T(1)\otimes T(a))*T(b)$\\$T(a)*T(1)*T(1)*T(b)$\end{tabular} & \begin{tabular}[t]{@{}l@{}}$\mathcal Z_3*\mathcal D_5$\end{tabular} & \begin{tabular}[t]{@{}l@{}}$\mathcal Z_3*\PB$\end{tabular}\\[4pt]
34 & \begin{tabular}[t]{@{}l@{}}$\mathcal D_5*T(1)*T(1)$\\$T(1)*\mathcal D_5*T(1)$\\$T(a)*T(1)*T(b)*T(c)$\\$T(1)*T(a)*T(b)*T(c)$\\$\mathcal H_{\mathrm{sym}}(2)$\end{tabular} & \begin{tabular}[t]{@{}l@{}}$\PB*T(1)*T(1)$\\$T(1)*\PB*T(1)$\\$\mathcal Z_3*\BMW(T)$\\$\mathcal H_{\mathrm{pair}}(2)$\end{tabular} & ---\\[4pt]
35 & \begin{tabular}[t]{@{}l@{}}$\mathcal D_5*T(1)*T(t)$\\$T(1)*\mathcal D_5*T(t)$\\$\mathcal D_5*T(t)*T(1)$\\$\BMW(T)*T(1)*T(1)$\\$T(1)*\BMW(T)*T(1)$\\$T(a)*T(b)*T(c)*T(d)$\\$\mathcal H_{\mathrm{sym}}(a)$\end{tabular} & \begin{tabular}[t]{@{}l@{}}$\PB*T(t)*T(1)$\\$\PB*T(1)*T(t)$\\$T(t)*\PB*T(1)$\\$\mathcal H_{\mathrm{pair}}(a)$\end{tabular} & ---\\[4pt]
36 & \begin{tabular}[t]{@{}l@{}}$T(1)*\BMW(T)*T(t)$\\$\mathcal D_5*T(a)*T(b)$\\$\BMW(T)*T(1)*T(t)$\\$\BMW(T)*T(t)*T(1)$\\$T(a)*\mathcal D_5*T(b)$\end{tabular} & \begin{tabular}[t]{@{}l@{}}$\PB*T(a)*T(b)$\\$T(a)*\PB*T(b)$\end{tabular} & ---\\[4pt]
37 & \begin{tabular}[t]{@{}l@{}}$(T(a)\otimes T(b))*T(1)$\\$\mathcal D_7*T(1)$\\$\mathcal D_5*\mathcal D_5$\\$\BMW(T)*T(a)*T(b)$\\$T(a)*\BMW(T)*T(b)$\end{tabular} & \begin{tabular}[t]{@{}l@{}}$\PB*\mathcal D_5$\end{tabular} & \begin{tabular}[t]{@{}l@{}}$\PB*\PB$\end{tabular}\\[4pt]
38 & \begin{tabular}[t]{@{}l@{}}$\mathcal D_5*\BMW(T)$\\$\mathcal D_7*T(t)$\\$(T(a)\otimes T(b))*T(c)$\end{tabular} & \begin{tabular}[t]{@{}l@{}}$\PB*\BMW(T)$\end{tabular} & ---\\[4pt]
39 & \begin{tabular}[t]{@{}l@{}}$\BMW(T_1)*\BMW(T_2)$\end{tabular} & --- & ---\\[4pt]
41 & $\mathbb Z_9$ & --- & ---\\
\end{longtable}
\endgroup
\begin{proof}
At rank five, solving the fusion, Fourier-matching,
and common-edge consistency equations with the
exact-support conditions gives the following results.
Up to shading reversal, 64 support pairs survive:
38 in the self-dual branch,
22 in the $(12)$ branch, and
4 in the $(12)(34)$ branch.

The surviving data are identified with the models
in Table~\ref{tab:five} on the stated parameter domains.
The enumeration, exclusions, and model identifications
are recorded in the accompanying computational
certificates.
Proposition~\ref{prop:finite-classification-scheme}
therefore gives the classification.
\end{proof}

\begin{remark}\label{sec:isolated.}
The isolated self-dual entry with $m=41$ is the known example arising
from $\mathbb Z_9$. It has quantum dimensions $(1,2,2,2,2)$ and
Jones index $9$. It is primitive but is not an essentially new example.
The essentially new families occur in the self-dual and $(12)$ branches.
\end{remark}

\subsection{Two new families}
\label{sec:examples}
\label{sec:families}

The rank-five classification yields two new
one-parameter families,
$\mathcal H_{\mathrm{sym}}(a)$ and
$\mathcal H_{\mathrm{pair}}(a)$.
We describe their fusion rules, Fourier matrices,
duality types, and specializations.

Both families have
\[
a\in\mathbb C\setminus\{0,1\},\qquad
d=d'=(1,a,a,2a-2,1),\qquad \delta^2=4a.
\]
In the distinguished bases $\{p_i\}_{i=0}^4$ and
$\{p'_i\}_{i=0}^4$, ordinary multiplication is
$p_ip_j=\delta_{ij}p_i$ and $p'_ip'_j=\delta_{ij}p'_i$.
Convolution is commutative, with units $p_0$ and $p'_0$,
and has the same coefficients on both sides.

\paragraph{The self-dual family.}

In $\mathcal H_{\mathrm{sym}}(a)$, every label is self-dual.
The products involving $p_1$ and $p_2$ are
\begin{equation}\label{eq:symfive}
p_1\star p_1=p_2\star p_2
=ap_0+\frac a2p_3,
\qquad
p_1\star p_2=\frac a2p_3+ap_4.
\end{equation}
The remaining non-unit products are
\begin{equation}\label{eq:commonfive}
\begin{gathered}
p_1\star p_3=p_2\star p_3=(a-1)(p_1+p_2),\\[3pt]
p_1\star p_4=p_2,\qquad
p_2\star p_4=p_1,\\[3pt]
p_3\star p_3
=2(a-1)(p_0+p_4)+2(a-2)p_3,\\[3pt]
p_3\star p_4=p_3,\qquad
p_4\star p_4=p_0.
\end{gathered}
\end{equation}
Together with commutativity and the unit relation,
these formulas specify the entire convolution table.
The minus-side formulas are obtained by replacing
$p_i$ by $p'_i$.

The Fourier matrices are
\begin{equation}\label{eq:fivefourier}
S=S'=\delta^{-1}X,\qquad
X=
\begin{pmatrix}
1&a&a&2a-2&1\\
1&r&-r&0&-1\\
1&-r&r&0&-1\\
1&0&0&-2&1\\
1&-a&-a&2a-2&1
\end{pmatrix},
\qquad r^2=a.
\end{equation}
Here $X^2=4aI_5$, and hence $S'S=SS'=I_5$,
as required by the self-duality of every label.

\paragraph{The family with duality $(12)$.}

In $\mathcal H_{\mathrm{pair}}(a)$, duality exchanges
labels $1$ and $2$ and fixes $0,3,4$.
The convolution table differs from the self-dual table
in the following products:
\begin{equation}\label{eq:pairfive}
p_1\star p_1=p_2\star p_2
=\frac a2p_3+ap_4,
\qquad
p_1\star p_2=ap_0+\frac a2p_3.
\end{equation}
The remaining non-unit products are given by
\eqref{eq:commonfive}, with the same formulas on
the minus side in the primed basis.

The Fourier matrices are again $S=S'=\delta^{-1}X$,
where $X$ is the matrix in \eqref{eq:fivefourier},
now with $r^2=-a$.
Thus
\[
X^2=4aD_{12},
\qquad
S'S=SS'=D_{12},
\]
where $D_{12}$ is the permutation matrix exchanging
labels $1$ and $2$.
The change in duality is therefore reflected both
in the vacuum coefficients of the convolution table
and in the square of the Fourier matrix.

\paragraph{Properties and specializations.}

The parameter restriction $a\ne0,1$ ensures that
the modulus and all quantum dimensions are nonzero.
The square identities above show that the Fourier
matrices are invertible throughout this domain.
The two choices of $r$ are related by an admissible
relabelling on one side, together with the corresponding
transport of the Fourier matching.

In both families,
$
N_{33}^{3}=2(a-2).
$
Within the allowed parameter domain, $a=2$ is the
only value at which the fusion support changes.
At this value the face $(3,3,3)$ disappears, giving
\[
\dim P_{3,\pm}=
\begin{cases}
34,&a=2,\\
35,&a\ne2.
\end{cases}
\]
The active faces of the two fusion complexes,
including their specializations at $a=2$, are listed
in Appendix~\ref{app:rank-five-supports}.

At $a=2$, both families have
$
d=d'=(1,2,2,2,1),
\delta^2=8.
$
They remain distinct at this specialization, since
their duality permutations have different cycle types.

The two families are essentially new relative to the
exchange-relation and lower-rank product constructions
in the classification, in the following precise sense.

\begin{proposition}
Neither family is an exchange-relation planar algebra,
a nontrivial tensor product, or a nontrivial ordered
free product of lower-rank planar algebras.
\end{proposition}

\begin{proof}
In both families, the fusion link at label $3$
contains a square, so the exchange-relation forest
criterion excludes the exchange relation.
A nontrivial tensor-product decomposition is excluded
by the prime rank $5$.

Suppose, for a contradiction, that $P\cong A*B$ is a
nontrivial ordered free product of lower-rank planar algebras.
Let $\{a_i\}$ and $\{b_j\}$ be the minimal idempotents
of $A_{2,+}$ and $B_{2,+}$, with $a_0=e_A$ and $b_0=e_B$
the normalized Jones projections.
The minimal idempotents of $(A*B)_{2,+}$ are
\[
P_i=a_i\otimes e_B,\qquad
Q_j=1_A\otimes b_j\quad(j\ne0).
\]
The first family is closed under normalized convolution.
Moreover, since
$a_i\star1_A=d_A(a_i)1_A$ and $e_B\star b_j=b_j$,
\[
P_i\star Q_j
=(a_i\star1_A)\otimes(e_B\star b_j)
=d_A(a_i)Q_j.
\]

A planar-algebra isomorphism permutes the minimal
idempotents and preserves the Jones projection,
convolution, and trace.
Consequently, the images of the $P_i$ determine a
fusion-closed label subset $J\subset I$ satisfying
\[
0\in J,\qquad 1<|J|<|I|,\qquad
p_s\star p_t=d_s p_t
\quad(s\in J,\ t\notin J).
\]
Here the strict inequalities follow because both
factors have $2$-box dimension greater than one.

We now determine the possible subsets $J$ from the
displayed fusion rules.
In the self-dual family, a fusion-closed subset
containing $1$ must contain $3$ by $p_1\star p_1$,
then $2$ by $p_1\star p_3$, and finally $4$ by
$p_1\star p_2$.
In the family with duality $(12)$, the product
$p_1\star p_1$ forces both $3$ and $4$, and
$p_1\star p_4=p_2$ then forces $2$.
The same conclusions hold with $1$ and $2$ exchanged.
All coefficients used here are nonzero because
$a\ne0,1$.
Thus a proper fusion-closed subset containing $0$
contains neither $1$ nor $2$.
Furthermore, if it contains $3$, then
\[
p_3\star p_3
=2(a-1)(p_0+p_4)+2(a-2)p_3
\]
forces it to contain $4$.
The only possibilities with $1<|J|<|I|$ are therefore
\[
J=\{0,4\}\qquad\text{and}\qquad J=\{0,3,4\}.
\]

For either choice, $4\in J$ and $1\notin J$.
The necessary free-product identity would give
$p_4\star p_1=d_4p_1=p_1$.
However, both families satisfy $p_4\star p_1=p_2$,
contradicting the distinctness of the minimal
idempotents $p_1$ and $p_2$.
This excludes every nontrivial ordered free-product
decomposition of the stated kind.
The argument also applies at $a=2$.
\end{proof}

\begin{remark}
At $a=2$, the fusion and Fourier data of
$\mathcal H_{\mathrm{pair}}(2)$ agree, after compatible
choices of labels and Fourier convention, with those
of the group--subgroup subfactor
\[
R^G\subset R^H,\qquad
G=\langle u,v\mid u^8=v^2=1,\;
vuv^{-1}=u^3\rangle,\qquad
H=\langle v\rangle,
\]
where $G$ acts outerly on the hyperfinite
$\mathrm{II}_1$ factor $R$.
Since $|G|=16$ and $|H|=2$, this inclusion has index $8$,
in agreement with $\delta^2=8$.

Our classification of the two families is algebraic.
We do not determine the parameter values for which
they admit a positive planar trace or a subfactor
realization. Determining these values is a natural
direction for further study.
\end{remark}

\paragraph{Acknowledgements.}
The authors thank Fan Lu for many helpful discussions.
This work was supported in part by the Beijing Natural
Science Foundation Key Program (Grant No.~Z220002)
and the Beijing Natural Science Foundation
(Grant No.~Z221100002722017).

\paragraph{Use of AI tools.}
The authors developed all mathematical ideas and proofs,
including the consistency theorem, the concept of a fusion
complex, the cohomological characterization
$H^2(K;\C)=0$, and the mathematical design of the
classification algorithm. The authors also derived all
equations, carried out the diagrammatic calculations,
analyzed the results, and determined the structure and
presentation of the paper.

Large language models assisted with implementing and
optimizing the algorithms, identifying possible errors
in code and formulas, preparing figures, and drafting
and polishing the exposition. They also assisted with
literature searches concerning possible identifications
of the resulting examples with previously known
constructions.

All AI-assisted work was carried out under the authors'
direction and review. The authors checked all proofs
and derivations, reviewed and tested all AI-assisted code,
and verified the computational results and literature
identifications reported in the paper. The authors take
full responsibility for the mathematical statements,
algorithms, computations, and conclusions.

\appendix
\section{Model conventions and diagram catalogues}

\subsection{Conventions for the known models}
\label{app:rank3models}
\label{sec:fourmatrices}

The known models in the classification tables are described
in \cite{BischJonesPA,BischJonesPAII,LiuER,BJL2014,LuThesis}.
Here we record the parameter conventions used in this paper.

We write $T(a)$ for the rank-two Temperley--Lieb model
with quantum dimensions $(1,a)$ and $\delta^2=1+a$.
The specialization $a=1$ is listed separately.
All unfixed Temperley--Lieb parameters belong to
$\mathcal U=\mathbb C\setminus\{0,-1,1\}$.

The symbol $\mathcal Z_q$ denotes the cyclic group planar
algebra, and
$\mathcal D_q=\mathcal S^{\mathbb Z_2\subset
\mathbb Z_q\rtimes\mathbb Z_2}$ for $q=5,7$.
We write $\PB$ for the binary-tree exchange-relation
planar algebra of \cite{LuLiuER}.
The symbols $\otimes$ and $*$ denote tensor product
and ordered free product, respectively.

For $\BMW(T)$, put $D=T^2+T-1$. Our parameterization is
\[
d=d'=\bigl(1,(T-1)D,(T^2-2)D\bigr),
\qquad
\delta=\pm(T^2+T-2).
\]
The fusion formulas agree with those in
\cite[Lemma~3.3]{BJL2014}, after accounting for
the convolution normalization.
In our parameterization, nonzero quantum dimensions
and modulus, together with the prescribed fourteen-face
fusion support, give the domain
\[
T(T-1)(T+1)(T+2)(T^2+T-1)(T^2-2)(T^2-T-1)\ne0.
\]
This is an algebraic support condition, without
a positivity assumption.
When two BMW factors occur in a product, their parameters
$T_1,T_2$ vary independently in this domain.
At $T^2-T-1=0$, the model specializes to
$\mathcal D_5$, which is listed separately.

\subsection{Explicit support complexes}\label{app:supports}
A triple $(x,y,z)$ denotes the face $N_{zy}^{x}\ne0$. These lists define the complex labels used above. They use the canonical labelings of the supplied calculation. The normalized model bases in the main text can differ by the recorded primitive-label permutations. The rank-two supports were given in Section~\ref{sec:cohomology}.
\subsubsection{\texorpdfstring{$\mathcal Z_9$}{Z9}}
\[\begin{gathered}
(0,0,0),\ (0,1,2),\ (0,2,1),\ (1,0,1),\ (1,1,0)\\(1,2,2),\ (2,0,2),\ (2,1,1),\ (2,2,0)
\end{gathered}\]
\subsubsection{\texorpdfstring{$\mathcal F_{10}$}{F10}}
\[\begin{gathered}
(0,0,0),\ (0,1,1),\ (0,2,2),\ (1,0,1),\ (1,1,0)\\(1,1,2),\ (1,2,1),\ (2,0,2),\ (2,1,1),\ (2,2,0)
\end{gathered}\]
\subsubsection{\texorpdfstring{$\mathcal F_{11}^{L}$}{F11L}}
\[\begin{gathered}
(0,0,0),\ (0,1,1),\ (0,2,2),\ (1,0,1),\ (1,1,0)\\(1,1,1),\ (1,1,2),\ (1,2,1),\ (2,0,2),\ (2,1,1)\\(2,2,0)
\end{gathered}\]
\subsubsection{\texorpdfstring{$\mathcal F_{11}^{R}$}{F11R}}
\[\begin{gathered}
(0,0,0),\ (0,1,1),\ (0,2,2),\ (1,0,1),\ (1,1,0)\\(1,1,1),\ (1,2,2),\ (2,0,2),\ (2,1,2),\ (2,2,0)\\(2,2,1)
\end{gathered}\]
\subsubsection{\texorpdfstring{$\mathcal F_{12}$}{F12}}
\[\begin{gathered}
(0,0,0),\ (0,1,1),\ (0,2,2),\ (1,0,1),\ (1,1,0)\\(1,1,1),\ (1,1,2),\ (1,2,1),\ (2,0,2),\ (2,1,1)\\(2,2,0),\ (2,2,2)
\end{gathered}\]
\subsubsection{\texorpdfstring{$\mathcal D_{13}$}{D13}}
\[\begin{gathered}
(0,0,0),\ (0,1,1),\ (0,2,2),\ (1,0,1),\ (1,1,0)\\(1,1,2),\ (1,2,1),\ (1,2,2),\ (2,0,2),\ (2,1,1)\\(2,1,2),\ (2,2,0),\ (2,2,1)
\end{gathered}\]
\subsubsection{\texorpdfstring{$\mathcal E_{13}$}{E13}}
\[\begin{gathered}
(0,0,0),\ (0,1,2),\ (0,2,1),\ (1,0,1),\ (1,1,0)\\(1,1,1),\ (1,1,2),\ (1,2,1),\ (2,0,2),\ (2,1,2)\\(2,2,0),\ (2,2,1),\ (2,2,2)
\end{gathered}\]
\subsubsection{\texorpdfstring{$\mathcal W_{14}$}{W14}}
\[\begin{gathered}
(0,0,0),\ (0,1,1),\ (0,2,2),\ (1,0,1),\ (1,1,0)\\(1,1,1),\ (1,1,2),\ (1,2,1),\ (1,2,2),\ (2,0,2)\\(2,1,1),\ (2,1,2),\ (2,2,0),\ (2,2,1)
\end{gathered}\]

\subsubsection{Rank-five supports of the two new families}
\label{app:new-five-supports}
\label{app:rank-five-supports}

In the model labelings of Section~\ref{sec:families}, the active
faces for $a\notin\{0,1,2\}$ are as follows.

\paragraph{\texorpdfstring{$\mathcal H^{\mathrm{sym}}_{35}$}{Hsym35}}
\[\begin{gathered}
(0,0,0),\ (0,1,1),\ (0,2,2),\ (0,3,3),\ (0,4,4),\\
(1,0,1),\ (1,1,0),\ (1,1,3),\ (1,2,3),\ (1,2,4),\\
(1,3,1),\ (1,3,2),\ (1,4,2),\ (2,0,2),\ (2,1,3),\\
(2,1,4),\ (2,2,0),\ (2,2,3),\ (2,3,1),\ (2,3,2),\\
(2,4,1),\ (3,0,3),\ (3,1,1),\ (3,1,2),\ (3,2,1),\\
(3,2,2),\ (3,3,0),\ (3,3,3),\ (3,3,4),\ (3,4,3),\\
(4,0,4),\ (4,1,2),\ (4,2,1),\ (4,3,3),\ (4,4,0)
\end{gathered}\]

\paragraph{\texorpdfstring{$\mathcal H^{\mathrm{pair}}_{35}$}{Hpair35}}
\[\begin{gathered}
(0,0,0),\ (0,1,2),\ (0,2,1),\ (0,3,3),\ (0,4,4),\\
(1,0,1),\ (1,1,0),\ (1,1,3),\ (1,2,3),\ (1,2,4),\\
(1,3,1),\ (1,3,2),\ (1,4,2),\ (2,0,2),\ (2,1,3),\\
(2,1,4),\ (2,2,0),\ (2,2,3),\ (2,3,1),\ (2,3,2),\\
(2,4,1),\ (3,0,3),\ (3,1,1),\ (3,1,2),\ (3,2,1),\\
(3,2,2),\ (3,3,0),\ (3,3,3),\ (3,3,4),\ (3,4,3),\\
(4,0,4),\ (4,1,1),\ (4,2,2),\ (4,3,3),\ (4,4,0)
\end{gathered}\]

Both sides have the listed support.
At $a=2$, only $N_{33}^{3}=2(a-2)$ vanishes, so
$\mathcal H^{\mathrm{sym}}_{34}$ and
$\mathcal H^{\mathrm{pair}}_{34}$ are obtained by deleting
$(3,3,3)$ from the respective lists.

\Needspace{9\baselineskip}
\subsection{Notation summary}
% ============================================================

{\small
\begin{longtable}{>{\raggedright\arraybackslash}p{0.20\textwidth}
                  >{\raggedright\arraybackslash}p{0.70\textwidth}}
\toprule
Symbol & Meaning \\
\midrule
\endfirsthead
\toprule
Symbol & Meaning \\
\midrule
\endhead
$\F$ & Counterclockwise $90^\circ$ Fourier rotation. Here $\F^2$ is the half-turn. \\
$p_i,p'_a$ & Plus- and minus-side minimal projections. The elements $p_0,p'_0$ are Jones projections. \\
$N_{kj}^{i}$ & Fusion coefficients in the convention $(i,j,k)\in K\iff N_{kj}^{i}\neq0$. \\
$K(B)$ & Tripartite fusion complex associated with the support of $B$. \\
$\chi_a^+,\chi_i^-$ & Fusion characters of $(P_{2,+},*_+)$ and $(P_{2,-},*_-)$ after opposite-side label matching. \\
$c$ & Native plus-shaded triangle coefficient table. \\
$c'$ & Native minus-shaded coefficient table in the $p'_a$ basis. \\
$\widetilde c$ & Common-shading coefficient table of $W(i,j,k)$ in the plus-side basis. \\
$S,S'$ & Geometric Fourier matrices in target-row/source-column convention. Here $S_{ai}=\delta^{-1}\chi_a^+(p_i)$ and $S'_{ia}=\delta^{-1}\chi_i^-(p'_a)$. \\
$\mathrm{Char}_+,\mathrm{Char}_-$ & Raw fusion-character matrices, with characters in rows and distinguished basis elements in columns. \\
$C_{\mathrm{dual},+},C_{\mathrm{dual},-}$ & Duality permutation matrices. \\
$Q_{ij}^{k}$ & Weighted fusion coefficient $N_{ij}^{k}d_k$. \\
\bottomrule
\end{longtable}
}

% ============================================================

\subsection{The reduced 4-box diagrams}\label{app:O4}
Here $O_4$ denotes 4-box diagrams with no internal $1$-, $2$-, or
$3$-face, after removing vertex-free circles. Boundary points are fixed. Rotations are not identified. Vertex labels
and checkerboard shading are suppressed in the drawings.

Types $1$--$82$ are forests: their counts by interior vertex number are
$14,28,28,12$. Types $83,84$ have one quadrilateral face. To see that no
other reduced type with an internal face occurs, a component with a cycle
must use all eight boundary points by Lemma~\ref{cons:lowbox}.
In a cycle of length $r$ enclosing only faces of degree at least four,
disc Euler counting gives $2r\le s\le3r-4$, where $s$ is the sum of
boundary vertex degrees in the enclosed subgraph. Hence $r\ge4$.
For the whole connected diagram, write $v,E,f,L$ for its interior vertex,
edge, internal face, and exterior boundary-walk counts. Then
$E=2v+4$, $f=v-3$, and $4f\le2E-L$, so $L\le20$.
The eight boundary legs contribute $16$ and the cyclic part at least $4$.
Equality forces a quadrilateral with two external legs at each vertex.
Its two boundary placements give $83,84$.

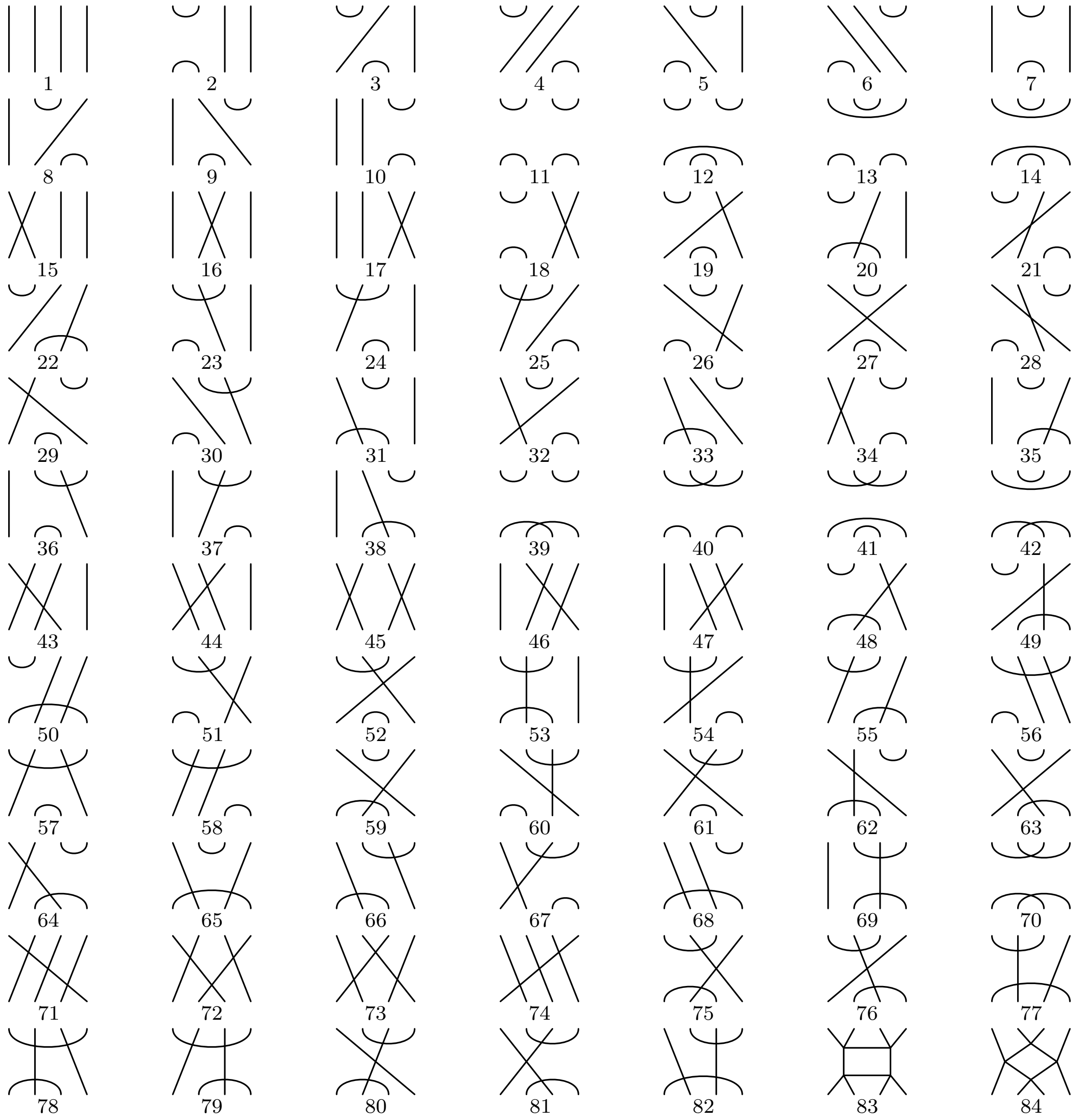
\begin{center}
\begin{minipage}[c][1.24cm][c]{.137\linewidth}\centering$\scalebox{1.55}{ $\OnFour{1}$ }$\par\scriptsize 1\end{minipage}%
\begin{minipage}[c][1.24cm][c]{.137\linewidth}\centering$\scalebox{1.55}{ $\OnFour{2}$ }$\par\scriptsize 2\end{minipage}%
\begin{minipage}[c][1.24cm][c]{.137\linewidth}\centering$\scalebox{1.55}{ $\OnFour{3}$ }$\par\scriptsize 3\end{minipage}%
\begin{minipage}[c][1.24cm][c]{.137\linewidth}\centering$\scalebox{1.55}{ $\OnFour{4}$ }$\par\scriptsize 4\end{minipage}%
\begin{minipage}[c][1.24cm][c]{.137\linewidth}\centering$\scalebox{1.55}{ $\OnFour{5}$ }$\par\scriptsize 5\end{minipage}%
\begin{minipage}[c][1.24cm][c]{.137\linewidth}\centering$\scalebox{1.55}{ $\OnFour{6}$ }$\par\scriptsize 6\end{minipage}%
\begin{minipage}[c][1.24cm][c]{.137\linewidth}\centering$\scalebox{1.55}{ $\OnFour{7}$ }$\par\scriptsize 7\end{minipage}%
\par
\begin{minipage}[c][1.24cm][c]{.137\linewidth}\centering$\scalebox{1.55}{ $\OnFour{8}$ }$\par\scriptsize 8\end{minipage}%
\begin{minipage}[c][1.24cm][c]{.137\linewidth}\centering$\scalebox{1.55}{ $\OnFour{9}$ }$\par\scriptsize 9\end{minipage}%
\begin{minipage}[c][1.24cm][c]{.137\linewidth}\centering$\scalebox{1.55}{ $\OnFour{10}$ }$\par\scriptsize 10\end{minipage}%
\begin{minipage}[c][1.24cm][c]{.137\linewidth}\centering$\scalebox{1.55}{ $\OnFour{11}$ }$\par\scriptsize 11\end{minipage}%
\begin{minipage}[c][1.24cm][c]{.137\linewidth}\centering$\scalebox{1.55}{ $\OnFour{12}$ }$\par\scriptsize 12\end{minipage}%
\begin{minipage}[c][1.24cm][c]{.137\linewidth}\centering$\scalebox{1.55}{ $\OnFour{13}$ }$\par\scriptsize 13\end{minipage}%
\begin{minipage}[c][1.24cm][c]{.137\linewidth}\centering$\scalebox{1.55}{ $\OnFour{14}$ }$\par\scriptsize 14\end{minipage}%
\par
\begin{minipage}[c][1.24cm][c]{.137\linewidth}\centering$\scalebox{1.55}{ $\OnFour{15}$ }$\par\scriptsize 15\end{minipage}%
\begin{minipage}[c][1.24cm][c]{.137\linewidth}\centering$\scalebox{1.55}{ $\OnFour{16}$ }$\par\scriptsize 16\end{minipage}%
\begin{minipage}[c][1.24cm][c]{.137\linewidth}\centering$\scalebox{1.55}{ $\OnFour{17}$ }$\par\scriptsize 17\end{minipage}%
\begin{minipage}[c][1.24cm][c]{.137\linewidth}\centering$\scalebox{1.55}{ $\OnFour{18}$ }$\par\scriptsize 18\end{minipage}%
\begin{minipage}[c][1.24cm][c]{.137\linewidth}\centering$\scalebox{1.55}{ $\OnFour{19}$ }$\par\scriptsize 19\end{minipage}%
\begin{minipage}[c][1.24cm][c]{.137\linewidth}\centering$\scalebox{1.55}{ $\OnFour{20}$ }$\par\scriptsize 20\end{minipage}%
\begin{minipage}[c][1.24cm][c]{.137\linewidth}\centering$\scalebox{1.55}{ $\OnFour{21}$ }$\par\scriptsize 21\end{minipage}%
\par
\begin{minipage}[c][1.24cm][c]{.137\linewidth}\centering$\scalebox{1.55}{ $\OnFour{22}$ }$\par\scriptsize 22\end{minipage}%
\begin{minipage}[c][1.24cm][c]{.137\linewidth}\centering$\scalebox{1.55}{ $\OnFour{23}$ }$\par\scriptsize 23\end{minipage}%
\begin{minipage}[c][1.24cm][c]{.137\linewidth}\centering$\scalebox{1.55}{ $\OnFour{24}$ }$\par\scriptsize 24\end{minipage}%
\begin{minipage}[c][1.24cm][c]{.137\linewidth}\centering$\scalebox{1.55}{ $\OnFour{25}$ }$\par\scriptsize 25\end{minipage}%
\begin{minipage}[c][1.24cm][c]{.137\linewidth}\centering$\scalebox{1.55}{ $\OnFour{26}$ }$\par\scriptsize 26\end{minipage}%
\begin{minipage}[c][1.24cm][c]{.137\linewidth}\centering$\scalebox{1.55}{ $\OnFour{27}$ }$\par\scriptsize 27\end{minipage}%
\begin{minipage}[c][1.24cm][c]{.137\linewidth}\centering$\scalebox{1.55}{ $\OnFour{28}$ }$\par\scriptsize 28\end{minipage}%
\par
\begin{minipage}[c][1.24cm][c]{.137\linewidth}\centering$\scalebox{1.55}{ $\OnFour{29}$ }$\par\scriptsize 29\end{minipage}%
\begin{minipage}[c][1.24cm][c]{.137\linewidth}\centering$\scalebox{1.55}{ $\OnFour{30}$ }$\par\scriptsize 30\end{minipage}%
\begin{minipage}[c][1.24cm][c]{.137\linewidth}\centering$\scalebox{1.55}{ $\OnFour{31}$ }$\par\scriptsize 31\end{minipage}%
\begin{minipage}[c][1.24cm][c]{.137\linewidth}\centering$\scalebox{1.55}{ $\OnFour{32}$ }$\par\scriptsize 32\end{minipage}%
\begin{minipage}[c][1.24cm][c]{.137\linewidth}\centering$\scalebox{1.55}{ $\OnFour{33}$ }$\par\scriptsize 33\end{minipage}%
\begin{minipage}[c][1.24cm][c]{.137\linewidth}\centering$\scalebox{1.55}{ $\OnFour{34}$ }$\par\scriptsize 34\end{minipage}%
\begin{minipage}[c][1.24cm][c]{.137\linewidth}\centering$\scalebox{1.55}{ $\OnFour{35}$ }$\par\scriptsize 35\end{minipage}%
\par
\begin{minipage}[c][1.24cm][c]{.137\linewidth}\centering$\scalebox{1.55}{ $\OnFour{36}$ }$\par\scriptsize 36\end{minipage}%
\begin{minipage}[c][1.24cm][c]{.137\linewidth}\centering$\scalebox{1.55}{ $\OnFour{37}$ }$\par\scriptsize 37\end{minipage}%
\begin{minipage}[c][1.24cm][c]{.137\linewidth}\centering$\scalebox{1.55}{ $\OnFour{38}$ }$\par\scriptsize 38\end{minipage}%
\begin{minipage}[c][1.24cm][c]{.137\linewidth}\centering$\scalebox{1.55}{ $\OnFour{39}$ }$\par\scriptsize 39\end{minipage}%
\begin{minipage}[c][1.24cm][c]{.137\linewidth}\centering$\scalebox{1.55}{ $\OnFour{40}$ }$\par\scriptsize 40\end{minipage}%
\begin{minipage}[c][1.24cm][c]{.137\linewidth}\centering$\scalebox{1.55}{ $\OnFour{41}$ }$\par\scriptsize 41\end{minipage}%
\begin{minipage}[c][1.24cm][c]{.137\linewidth}\centering$\scalebox{1.55}{ $\OnFour{42}$ }$\par\scriptsize 42\end{minipage}%
\par
\begin{minipage}[c][1.24cm][c]{.137\linewidth}\centering$\scalebox{1.55}{ $\OnFour{43}$ }$\par\scriptsize 43\end{minipage}%
\begin{minipage}[c][1.24cm][c]{.137\linewidth}\centering$\scalebox{1.55}{ $\OnFour{44}$ }$\par\scriptsize 44\end{minipage}%
\begin{minipage}[c][1.24cm][c]{.137\linewidth}\centering$\scalebox{1.55}{ $\OnFour{45}$ }$\par\scriptsize 45\end{minipage}%
\begin{minipage}[c][1.24cm][c]{.137\linewidth}\centering$\scalebox{1.55}{ $\OnFour{46}$ }$\par\scriptsize 46\end{minipage}%
\begin{minipage}[c][1.24cm][c]{.137\linewidth}\centering$\scalebox{1.55}{ $\OnFour{47}$ }$\par\scriptsize 47\end{minipage}%
\begin{minipage}[c][1.24cm][c]{.137\linewidth}\centering$\scalebox{1.55}{ $\OnFour{48}$ }$\par\scriptsize 48\end{minipage}%
\begin{minipage}[c][1.24cm][c]{.137\linewidth}\centering$\scalebox{1.55}{ $\OnFour{49}$ }$\par\scriptsize 49\end{minipage}%
\par
\begin{minipage}[c][1.24cm][c]{.137\linewidth}\centering$\scalebox{1.55}{ $\OnFour{50}$ }$\par\scriptsize 50\end{minipage}%
\begin{minipage}[c][1.24cm][c]{.137\linewidth}\centering$\scalebox{1.55}{ $\OnFour{51}$ }$\par\scriptsize 51\end{minipage}%
\begin{minipage}[c][1.24cm][c]{.137\linewidth}\centering$\scalebox{1.55}{ $\OnFour{52}$ }$\par\scriptsize 52\end{minipage}%
\begin{minipage}[c][1.24cm][c]{.137\linewidth}\centering$\scalebox{1.55}{ $\OnFour{53}$ }$\par\scriptsize 53\end{minipage}%
\begin{minipage}[c][1.24cm][c]{.137\linewidth}\centering$\scalebox{1.55}{ $\OnFour{54}$ }$\par\scriptsize 54\end{minipage}%
\begin{minipage}[c][1.24cm][c]{.137\linewidth}\centering$\scalebox{1.55}{ $\OnFour{55}$ }$\par\scriptsize 55\end{minipage}%
\begin{minipage}[c][1.24cm][c]{.137\linewidth}\centering$\scalebox{1.55}{ $\OnFour{56}$ }$\par\scriptsize 56\end{minipage}%
\par
\begin{minipage}[c][1.24cm][c]{.137\linewidth}\centering$\scalebox{1.55}{ $\OnFour{57}$ }$\par\scriptsize 57\end{minipage}%
\begin{minipage}[c][1.24cm][c]{.137\linewidth}\centering$\scalebox{1.55}{ $\OnFour{58}$ }$\par\scriptsize 58\end{minipage}%
\begin{minipage}[c][1.24cm][c]{.137\linewidth}\centering$\scalebox{1.55}{ $\OnFour{59}$ }$\par\scriptsize 59\end{minipage}%
\begin{minipage}[c][1.24cm][c]{.137\linewidth}\centering$\scalebox{1.55}{ $\OnFour{60}$ }$\par\scriptsize 60\end{minipage}%
\begin{minipage}[c][1.24cm][c]{.137\linewidth}\centering$\scalebox{1.55}{ $\OnFour{61}$ }$\par\scriptsize 61\end{minipage}%
\begin{minipage}[c][1.24cm][c]{.137\linewidth}\centering$\scalebox{1.55}{ $\OnFour{62}$ }$\par\scriptsize 62\end{minipage}%
\begin{minipage}[c][1.24cm][c]{.137\linewidth}\centering$\scalebox{1.55}{ $\OnFour{63}$ }$\par\scriptsize 63\end{minipage}%
\par
\begin{minipage}[c][1.24cm][c]{.137\linewidth}\centering$\scalebox{1.55}{ $\OnFour{64}$ }$\par\scriptsize 64\end{minipage}%
\begin{minipage}[c][1.24cm][c]{.137\linewidth}\centering$\scalebox{1.55}{ $\OnFour{65}$ }$\par\scriptsize 65\end{minipage}%
\begin{minipage}[c][1.24cm][c]{.137\linewidth}\centering$\scalebox{1.55}{ $\OnFour{66}$ }$\par\scriptsize 66\end{minipage}%
\begin{minipage}[c][1.24cm][c]{.137\linewidth}\centering$\scalebox{1.55}{ $\OnFour{67}$ }$\par\scriptsize 67\end{minipage}%
\begin{minipage}[c][1.24cm][c]{.137\linewidth}\centering$\scalebox{1.55}{ $\OnFour{68}$ }$\par\scriptsize 68\end{minipage}%
\begin{minipage}[c][1.24cm][c]{.137\linewidth}\centering$\scalebox{1.55}{ $\OnFour{69}$ }$\par\scriptsize 69\end{minipage}%
\begin{minipage}[c][1.24cm][c]{.137\linewidth}\centering$\scalebox{1.55}{ $\OnFour{70}$ }$\par\scriptsize 70\end{minipage}%
\par
\begin{minipage}[c][1.24cm][c]{.137\linewidth}\centering$\scalebox{1.55}{ $\OnFour{71}$ }$\par\scriptsize 71\end{minipage}%
\begin{minipage}[c][1.24cm][c]{.137\linewidth}\centering$\scalebox{1.55}{ $\OnFour{72}$ }$\par\scriptsize 72\end{minipage}%
\begin{minipage}[c][1.24cm][c]{.137\linewidth}\centering$\scalebox{1.55}{ $\OnFour{73}$ }$\par\scriptsize 73\end{minipage}%
\begin{minipage}[c][1.24cm][c]{.137\linewidth}\centering$\scalebox{1.55}{ $\OnFour{74}$ }$\par\scriptsize 74\end{minipage}%
\begin{minipage}[c][1.24cm][c]{.137\linewidth}\centering$\scalebox{1.55}{ $\OnFour{75}$ }$\par\scriptsize 75\end{minipage}%
\begin{minipage}[c][1.24cm][c]{.137\linewidth}\centering$\scalebox{1.55}{ $\OnFour{76}$ }$\par\scriptsize 76\end{minipage}%
\begin{minipage}[c][1.24cm][c]{.137\linewidth}\centering$\scalebox{1.55}{ $\OnFour{77}$ }$\par\scriptsize 77\end{minipage}%
\par
\begin{minipage}[c][1.24cm][c]{.137\linewidth}\centering$\scalebox{1.55}{ $\OnFour{78}$ }$\par\scriptsize 78\end{minipage}%
\begin{minipage}[c][1.24cm][c]{.137\linewidth}\centering$\scalebox{1.55}{ $\OnFour{79}$ }$\par\scriptsize 79\end{minipage}%
\begin{minipage}[c][1.24cm][c]{.137\linewidth}\centering$\scalebox{1.55}{ $\OnFour{80}$ }$\par\scriptsize 80\end{minipage}%
\begin{minipage}[c][1.24cm][c]{.137\linewidth}\centering$\scalebox{1.55}{ $\OnFour{81}$ }$\par\scriptsize 81\end{minipage}%
\begin{minipage}[c][1.24cm][c]{.137\linewidth}\centering$\scalebox{1.55}{ $\OnFour{82}$ }$\par\scriptsize 82\end{minipage}%
\begin{minipage}[c][1.24cm][c]{.137\linewidth}\centering$\scalebox{1.55}{ $\OnFour{83}$ }$\par\scriptsize 83\end{minipage}%
\begin{minipage}[c][1.24cm][c]{.137\linewidth}\centering$\scalebox{1.55}{ $\OnFour{84}$ }$\par\scriptsize 84\end{minipage}%
\par
\end{center}

\subsection{Derivation of the two new families}
\label{app:new-family-derivation}

The classification algorithm determines the supports
of the two families in Section~\ref{sec:families}.
For these supports, we give a direct derivation of
the fusion coefficients and character tables.
This calculation explains the explicit formulas
and their parameter restrictions; realizability
is established by the classification and
reconstruction results above.

\paragraph{The family with duality $(12)$.}
Begin with the support of $\mathcal H_{\mathrm{pair}}$
and nonzero quantum dimensions
$d_1=d_2=x$, $d_3=y$, and $d_4=z$.
Commutativity, duality, and the vacuum identities give
\begin{align*}
p_1\star p_1=p_2\star p_2&=A p_3+B p_4,\\
p_1\star p_2&=x p_0+C p_3,\\
p_1\star p_3&=U p_1+V p_2,\\
p_2\star p_3&=V p_1+U p_2,\\
p_1\star p_4&=\gamma p_2,\qquad
p_2\star p_4=\gamma p_1,\\
p_3\star p_3&=y p_0+M p_3+L p_4,\\
p_3\star p_4&=\eta p_3,\qquad
p_4\star p_4=z p_0.
\end{align*}
Here $p_0$ is the convolution unit.

Applying the dimension character to $p_4\star p_4$
gives $z^2=z$. Hence $z=1$, and the dimension equations
for $p_1\star p_4$ and $p_3\star p_4$ give
$\gamma=\eta=1$.
Frobenius reciprocity gives $B=x$ and $L=y$.
Associativity of
$(p_1\star p_1)\star p_4=p_1\star(p_1\star p_4)$
then gives $A=C$.

Frobenius reciprocity and the dimension equation
for $p_1\star p_2$ yield
\[
Ay=xV,\qquad Cy=xU,\qquad x+Cy=x^2.
\]
Thus $U=V=x-1$.
The dimension equation for $p_1\star p_3$ now gives
$xy=2x(x-1)$, so $y=2x-2$.
Consequently,
\[
A=C=\frac{x}{2},\qquad B=x,\qquad
L=2x-2,\qquad M=2x-4,
\]
where the last identity follows from the dimension
equation for $p_3\star p_3$.
Putting $a=x$ recovers
\eqref{eq:pairfive} and \eqref{eq:commonfive}, with
\[
d=(1,a,a,2a-2,1),\qquad \delta^2=4a.
\]
Applying the same calculation on the minus side gives
a parameter $a'$. The common-modulus equation forces
$4a=4a'$, hence $a=a'$.

\paragraph{The self-dual family.}
For the self-dual support, the corresponding calculation
gives the same dimension vector and the same products
in \eqref{eq:commonfive}. The remaining products are
\[
p_1\star p_1=p_2\star p_2
=ap_0+\frac a2p_3,\qquad
p_1\star p_2=\frac a2p_3+ap_4.
\]
These are precisely \eqref{eq:symfive}.
Thus the two families differ in these products and
in the duality on labels $1,2$.

\paragraph{Character tables and Fourier matrices.}
For a character $\chi$ of the paired fusion algebra,
write $z_i=\chi(p_i)$. Since $p_0$ is the convolution
unit, $z_0=1$. The products involving $p_4$ give
\[
z_4^2=1,\qquad z_2=z_1z_4,\qquad
z_3(z_4-1)=0.
\]
If $z_4=-1$, then $z_3=0$, $z_2=-z_1$, and
the equation for $p_1\star p_1$ gives $z_1^2=-a$.
If $z_4=1$, then $z_2=z_1$. When $z_1=0$,
the equation for $p_1\star p_2$ gives $z_3=-2$.
When $z_1\ne0$, the equation for $p_1\star p_3$
gives $z_3=2a-2$, and then $z_1^2=a^2$.

The five characters therefore give the rows
\[
\begin{gathered}
(1,a,a,2a-2,1),\\
(1,r,-r,0,-1),\qquad (1,-r,r,0,-1),\\
(1,0,0,-2,1),\qquad (1,-a,-a,2a-2,1),
\end{gathered}
\qquad r^2=-a.
\]
For the self-dual family, the same calculation gives
these rows with $r^2=a$.
With the row ordering in \eqref{eq:fivefourier},
the resulting matrix $X$ satisfies
\[
X^2=
\begin{cases}
4aI_5,& \mathcal H_{\mathrm{sym}}(a),\\
4aD_{12},& \mathcal H_{\mathrm{pair}}(a).
\end{cases}
\]
Together with $\delta^2=4a$, this gives the
Fourier-square identities for the chosen matching
$S=S'=X/\delta$.

\paragraph{Parameter restrictions.}
Nonzero quantum dimensions and modulus require
$a\ne0,1$. For these values the five characters are
distinct, so the character matrix is invertible.
All prescribed active fusion coefficients remain
nonzero except that $N_{33}^{3}=2(a-2)$ vanishes
at $a=2$. Thus the generic supports have $35$ faces,
and deleting $(3,3,3)$ gives the $34$-face
specializations listed in
Appendix~\ref{app:rank-five-supports}.

\end{document}